\documentclass[11pt,reqno]{amsart}
\usepackage[top=1in,bottom=1in,left=1in,right=1in]{geometry}
\usepackage[T1]{fontenc}
\usepackage[utf8]{inputenc}
\usepackage{amsmath,amssymb,mathrsfs,mathtools}
\usepackage{microtype}
\usepackage[nocompress]{cite}
\usepackage[shortlabels]{enumitem}
\usepackage{aliascnt}
\usepackage{xcolor}
\usepackage[bookmarks=true,pdfstartview=FitH,pdfborder={0 0 0},
  colorlinks=true,citecolor=red,linkcolor=blue,urlcolor=blue]{hyperref}

\theoremstyle{plain}
\newtheorem{thm}{Theorem}[section]
\newaliascnt{lem}{thm}
\newaliascnt{prop}{thm}
\newaliascnt{cor}{thm}
\newtheorem{lem}[lem]{Lemma}
\newtheorem{prop}[prop]{Proposition}
\newtheorem{cor}[cor]{Corollary}
\aliascntresetthe{lem}
\aliascntresetthe{prop}
\aliascntresetthe{cor}

\theoremstyle{definition}
\newaliascnt{defn}{thm}
\newtheorem{defn}[defn]{Definition}
\aliascntresetthe{defn}

\theoremstyle{remark}
\newaliascnt{rem}{thm}
\newaliascnt{exa}{thm}
\newtheorem{rem}[rem]{Remark}
\newtheorem{exa}[exa]{Example}
\aliascntresetthe{rem}
\aliascntresetthe{exa}

\numberwithin{equation}{section}
\setlist[enumerate]{font=\normalfont}
\allowdisplaybreaks[1]
\newcommand{\R}{\mathbb R}
\newcommand{\N}{\mathbb N}
\newcommand{\Z}{\mathbb Z}
\newcommand{\T}{\mathbb T}
\newcommand{\E}{\mathbb E}
\newcommand{\bP}{\mathbb P} 
\newcommand{\PP}{\mathbb P_{\!\sigma}}  
\newcommand{\rH}{\mathrm H}
\newcommand{\rL}{\mathrm L}
\newcommand{\rC}{\mathrm C}
\newcommand{\rW}{\mathrm W}
\newcommand{\cA}{\mathcal A}
\newcommand{\cB}{\mathcal B}
\newcommand{\cF}{\mathcal F}
\newcommand{\cR}{\mathcal R}
\newcommand{\cS}{\mathcal S}
\newcommand{\cQ}{\mathcal Q}
\newcommand{\LL}{\mathscr L}
\newcommand{\bfe}{\boldsymbol e}
\newcommand{\wh}{\widehat}
\newcommand{\eps}{\varepsilon}
\newcommand{\del}{\partial}
\newcommand{\rd}{\,\mathrm d}
\DeclareMathOperator{\mdiv}{div}
\DeclareMathOperator{\Id}{Id}
\DeclareMathOperator{\Ran}{Ran}
\DeclareMathOperator{\tr}{tr}
\DeclarePairedDelimiter{\norm}{\lVert}{\rVert}
\newcommand{\ip}[2]{\langle #1,#2\rangle}
\newcommand{\Hs}[1]{\rH^{#1}_{\sigma,0}(\T^3)}

\title[One-component noise for 3D Navier--Stokes]{Regularizing the Navier--Stokes Equations\\by One Noise Component}

\author{Felix Brandt}
\address{Department of Mathematics, University of Southern California, Los Angeles, CA 90089, USA}
\email{fbrandt@usc.edu}
\author{Igor Kukavica}
\address{Department of Mathematics, University of Southern California, Los Angeles, CA 90089, USA}
\email{kukavica@usc.edu}

\hypersetup{pdftitle={Regularizing the Navier-Stokes Equations by One Noise Component}, pdfauthor={Felix Brandt, Igor Kukavica}}

\begin{document}

\begin{abstract}
We study the three-dimensional incompressible Navier--Stokes equations on the torus with multiplicative It\^o noise which is driven by a single Brownian motion and acts only through the vertical velocity component. Such a noise is strongly degenerate: An infinite-dimensional part of the velocity field is not forced at all. We show that it nevertheless has a regularizing effect. For a linear noise, strong solutions whose initial data have no unforced part are global with probability tending to one as the noise intensity grows. For a suitable nonlinear noise, the same holds for arbitrary $\rH^1$ data, together with smallness of the space-time $\rL^6$-norm of the vertical velocity. As a deterministic counterpart, we give a self-contained proof that a quintic damping of the vertical component yields global regularity, and the same argument shows that, under the linear noise, blow-up can only occur if the vertical velocity leaves~$\rL^6_{t,x}$.
\end{abstract}

\subjclass[2020]{35Q30, 60H15, 76D05, 35R60}
\keywords{Stochastic Navier--Stokes equations, multiplicative noise, degenerate noise, regularization by noise, one-component regularity, stochastic maximal regularity}

\maketitle

\section{Introduction}\label{sec:intro}

Whether smooth solutions of the three-dimensional incompressible Navier--Stokes equations remain smooth for all times is a well-known open problem. A useful way of isolating what is missing is to ask which additional mechanism, as weak as possible, would suffice for global regularity. In this paper we look at regularizing mechanisms acting only on a single velocity component.

Regularization by noise, the phenomenon that a suitable random perturbation improves the well-posedness of an evolution equation, has attracted considerable attention in recent years. For the three-dimensional Euler and Navier--Stokes equations and related models, several noise mechanisms are now known to prevent blow-up with high probability, among them large linear and superlinear multiplicative noise and transport noise~\cite{GHV:14,RZZ:14,FL:21,BMX:25,HLL:24}. In all of these results, the noise acts on the whole velocity field. We show that much less suffices: A noise acting only through the vertical velocity component, which leaves an infinite-dimensional part of the velocity field entirely unforced, still produces global strong solutions with probability close to one when its intensity is large (for arbitrary $\rH^1$ data in the nonlinear case).

\subsection*{A deterministic model}
On the torus $\T^3$ and with viscosity $\nu>0$, consider
\begin{equation}\label{eq:intro-tamed}
 \del_tu+\cB(u)-\nu\Delta u+\alpha\PP(u_3^5\bfe_3)=0,\qquad\alpha>0,
\end{equation}
where $\PP$ is the Leray projection, $\cB(u)=\PP[(u\cdot\nabla)u]$, and $u_3$ is the vertical velocity component. The additional term contains no derivatives. It is a damping and not an extra viscosity, and before the projection it only acts on the vertical component. Testing with $u$ gives
\begin{equation}\label{eq:intro-energy}
 \frac12\frac{\rd}{\rd t}\norm u_2^2+\nu\norm{\nabla u}_2^2+\alpha\norm{u_3}_6^6=0,
 \qquad\text{hence}\qquad
 \int_0^\infty\norm{u_3(t)}_6^6\rd t\le\frac{\norm{u_0}_2^2}{2\alpha}.
\end{equation}
The point of \eqref{eq:intro-energy} is that a space-time norm of a single component is controlled by the initial energy alone, independently of the unknown higher norms. This is precisely the kind of information required by one-component regularity criteria for the Navier--Stokes equations~\cite{KZ:06,CT:08,ZP:10}. For the quintic power, $(6,6)$ is the endpoint pair of exponents in the one-component criterion of Zhou and Pokorn\'y \cite{ZP:10}. Combining a partial damping with component-reduction criteria goes back to Yamazaki \cite{Yam:18}, and global regularity of \eqref{eq:intro-tamed} for arbitrary $\Hs1$ data is in this spirit. We give a self-contained proof in~\autoref{sec:tamed}.

\subsection*{Can noise replace the damping?}
Instead of \eqref{eq:intro-tamed}, we consider
\begin{equation}\label{eq:intro-sns}
 \rd u=[\nu\Delta u-\cB(u)]\rd t+G(u)\rd W_t,
\end{equation}
driven by a single real Brownian motion $W$ in the It\^o sense, where the coefficient $G(u)$ depends on~$u$ only through~$u_3$. At the level of the energy, noise does not help: It\^o's formula gives
\[
 \tfrac12\rd\norm u_2^2+\nu\norm{\nabla u}_2^2\rd t=\tfrac12\norm{G(u)}_2^2\rd t+(u,G(u))_2\rd W_t,
\]
and the correction $\frac12\norm{G(u)}_2^2$ has the wrong sign, so averaging cannot reproduce~\eqref{eq:intro-energy}. The picture changes for concave observables. If $\norm\cdot_H$ is a Hilbert norm and the noise is radial, $G(u)=f(u)u$ with a scalar $f(u)$, then away from $u=0$, we have
\begin{equation}\label{eq:intro-log}
 \rd\log\norm u_H=\biggl(\frac{(b(u),u)_H}{\norm u_H^2}-\frac12f(u)^2\biggr)\rd t+f(u)\rd W_t,
\end{equation}
where $b(u)$ denotes the drift. A large noise coefficient thus produces a large \emph{negative} drift for~$\log\norm u_H$. This mechanism is behind the non-explosion results for large linear multiplicative noise \cite{RZZ:14,GHV:14}, and, for coefficients $f(u)$ growing like a high Sobolev norm, behind the result of Bagnara, Maurelli, and Xu that nonlinear It\^o noise prevents blow-up for the Euler equations~\cite{BMX:25}. For a general coefficient, the It\^o contribution to the drift of $\log\norm u_H$ equals
\begin{equation}\label{eq:intro-geometry}
 \frac{\norm{G(u)}_H^2}{2\norm u_H^2}-\frac{(u,G(u))_H^2}{\norm u_H^4},
\end{equation}
which is negative for radial noise but may be positive if $G(u)$ is not aligned with~$u$. A coefficient acting only through $u_3$ is not aligned with the whole velocity field, so the logarithm of a full norm is not the right observable.

\subsection*{The structure we use}
The projected vertical forcing $Ku=\PP(u_3\bfe_3)$ has a simple form in Fourier variables. On a mode $\xi$ with nonzero horizontal part, it acts on the divergence-free amplitude $\wh u(\xi)$ as $\mu_\xi R_\xi$, where $R_\xi$ is the orthogonal projection onto the single direction $P_\xi\bfe_3$ and $\mu_\xi=|\xi_h|^2/|\xi|^2$. Consequently, every velocity field splits as $u=\cR u+\cS u$, where the \emph{unforced} part $\cS u$ has no vertical component and is never reached by the noise, while the \emph{accessible} part $\cR u$ is generated by~$u_3$. On each Fourier mode, the noise is radial on the accessible part. This restores a favorable It\^o correction for suitably \emph{polarized} observables, although the noise is not radial on the whole field.

\subsection*{Main results}
We consider two coefficients.

\emph{Linear one-component noise}, $G(u)=\beta Ku$ with a constant $\beta>0$. If the unforced part of the initial datum vanishes, $\cS u_0=0$, then the maximal strong solution with $u_0\in\Hs1$ is global with probability tending to one as $\beta\to\infty$ (\autoref{thm:linear}). No smallness of $u_0$ is required. The proof has two parts. First, a polarized version $X_{-1}$ of the Fourier--Wiener norm used by Lei and Lin \cite{LL:11} is a supermartingale as long as it is small, uniformly in $\beta$. Here the radial structure on each mode is crucial, since it excludes a positive It\^o correction. Second, on an initial layer of length $\Theta/\beta^2$, the explicit linear noise flow $\exp(\beta KW_t-\frac12\beta^2tK^2)$ contracts the datum into this small region with high probability, before the nonlinearity has time to act. We also obtain information on space-time integrals of~$u_3$.

\emph{Nonlinear one-component noise}, $G(u)=\beta\norm{\nabla\cR u}_2^k\cR u$ with $k\ge4$. Here, no condition on the initial datum is needed: For every $u_0\in\Hs1$ and $\eta>0$, the solution is global and satisfies $\int_0^\infty\!\int_{\T^3}|u_3|^6\rd x\rd t\le\eta$ with probability tending to one as $\beta\to\infty$ (\autoref{thm:nonlinear}). This is a probabilistic counterpart of~\eqref{eq:intro-energy}. The unforced part is compared with the solution $b$ of a reduced system of 2.5D type, which is globally regular for arbitrary data, and the proof rests on the joint functional
\[
 \norm{\nabla(\cS u-b)}_2^2+\lambda\log\bigl(1+\norm{\nabla\cR u}_2\bigr),
\]
in which viscosity controls small accessible amplitudes and the It\^o correction of the logarithm controls large ones.

\emph{One-component damping and a blow-up criterion.} For the damped system \eqref{eq:intro-tamed}, global well-posedness for arbitrary $\Hs1$ data follows from \eqref{eq:intro-energy} and the mechanism behind one-component regularity criteria, in the spirit of Yamazaki \cite{Yam:18}. We give a self-contained proof (\autoref{thm:tamed}). We formulate the argument as a pathwise a priori bound with an arbitrary constant (\autoref{prop:abstract}), which also applies to the linear noise once the stochastic terms in the energy balances are removed by random integrating factors. For every $\beta>0$, a solution with the linear noise can therefore only blow up in finite time on the event where $\int_0^{\tau_*}\norm{u_3}_6^6\rd t=\infty$ (\autoref{prop:lin-criterion}). Thus the quantity which the damping controls deterministically is the one which governs blow-up under the linear noise, and which \autoref{thm:nonlinear} controls with high probability.

Theorems~\ref{thm:linear} and~\ref{thm:nonlinear} hold with probability close to one for large intensities. We do not obtain almost sure global existence for a fixed~$\beta$. The It\^o interpretation of the noise is essential, see~\autoref{rem:ito}.

\subsection*{Relation to previous work}
As explained above, large linear multiplicative noise \cite{GHV:14,RZZ:14} and superlinear noise \cite{BMX:25} regularize through the It\^o correction of a concave observable, when the noise acts radially on the whole state. For the three-dimensional Navier--Stokes equations, superlinear coefficients such as $c\norm u_{\rH^1}^mu$, driven by a single Brownian motion, are covered by \cite{HLL:24}, and nonlinear Stratonovich noise and stochastic controls are considered in \cite{Bag:25,LC:26}. Transport noise is a different regularizing mechanism, based on enhanced dissipation \cite{FL:21,Luo:23}. It plays no role in this paper. For small $\rL^3$ data and small multiplicative noise, the stochastic Navier--Stokes equations have global solutions with high probability \cite{KX:26}. Pathwise stabilization by degenerate linear multiplicative noise has been studied for linear evolution equations \cite{Bie:11}, and stabilization of the Navier--Stokes equations near an equilibrium by localized finite-dimensional stochastic feedback in \cite{BDP:11}. 

The noise in this paper is degenerate in a strong sense: $K$ has an infinite-dimensional kernel, its nonzero eigenvalues $\mu_\xi$ accumulate at zero, and the unforced part $\cS u$ is not acted upon at all, but only driven by the nonlinearity. The results for radial noise therefore do not apply, and the main difficulty is the control of $\cS u$. For the linear noise, this is the reason for the condition $\cS u_0=0$ in \autoref{thm:linear}(a). For the nonlinear noise, $\cS u$ is compared with a reference flow. We are not aware of non-explosion results for the three-dimensional Navier--Stokes equations under multiplicative noise with a kernel of this type.

\subsection*{Local theory}
Both problems are set up with the theory of stochastic maximal regularity of Agresti and Veraar \cite{AV:22a,AV:22b,AV:24,AV:25}. In the strong setting $(\Hs2,\Hs1,\Hs0)$, both models are subcritical, which gives local well-posedness for $\Hs1$ data together with the blow-up criterion that the $\Hs1$-norm becomes unbounded at a finite explosion time (\autoref{sec:local}). The framework of \cite{DGT:11} could also be used for the local theory.

\subsection*{Organization}
\autoref{sec:prelim} collects notation, the one-component geometry and the main results. \autoref{sec:local} contains the local theory, the energy identities, and, for the linear noise, a continuation criterion and a pathwise description by a random change of variables. Theorems~\ref{thm:linear} and~\ref{thm:nonlinear} are proved in Sections~\ref{sec:linear} and~\ref{sec:nonlinear}, respectively. \autoref{sec:tamed} contains a self-contained proof of global regularity for the damped system~\eqref{eq:intro-tamed}, based on an abstract a priori bound, and uses the same bound to prove a one-component blow-up criterion for the linear noise. \autoref{app:log} justifies the logarithmic It\^o formula at zero.

\section{Preliminaries and main results}\label{sec:prelim}

\subsection{Function spaces and conventions}\label{ss:spaces}
We work on $\T^3=(\R/2\pi\Z)^3$, and spatial integrals are unnormalized. For $f\in \rL^1(\T^3;\mathbb C^d)$ and $\xi\in\Z^3$, let
\[
 \wh f(\xi)=(2\pi)^{-3}\int_{\T^3}f(x)e^{-i\xi\cdot x}\rd x,\qquad\text{so that}\qquad\norm f_2^2=(2\pi)^3\sum_{\xi\in\Z^3}|\wh f(\xi)|^2.
\]
Complex vectors carry their Euclidean norm. We write $x=(x_h,x_3)$, $\xi=(\xi_h,\xi_3)$, $u=(u_h,u_3)$, and~$\nabla_h$, $\mdiv_h$, $\Delta_h$ for the horizontal operators. For $s\in\R$, let
\[
 \Hs s=\Bigl\{u\in \rH^s(\T^3;\R^3):\ \mdiv u=0,\ \int_{\T^3}u\rd x=0\Bigr\},\qquad\norm u_{\Hs s}=\norm{(-\Delta)^{s/2}u}_2 .
\]
Thus, $\norm u_{\Hs1}=\norm{\nabla u}_2$ and $\norm u_{\Hs2}=\norm{\Delta u}_2$, and we write $\rL^2_\sigma=\Hs0$. On mean-zero functions these norms are equivalent to the usual Sobolev norms. We use the Sobolev and Gagliardo--Nirenberg inequalities for mean-zero functions without further comment, for instance
\[
 \norm f_6\le C\norm{\nabla f}_2,\qquad\norm f_3\le C\norm f_2^{1/2}\norm{\nabla f}_2^{1/2},\qquad\norm f_4\le C\norm f_2^{1/4}\norm{\nabla f}_2^{3/4},
\]
as well as the Poincar\'e inequality $\norm f_2\le\norm{\nabla f}_2$.

The Leray projection $\PP$ is the Fourier multiplier with symbol $P_\xi=\Id-\xi\otimes\xi/|\xi|^2$ for $\xi\ne0$ and~$P_0=0$; thus $\PP$ also removes the spatial mean. For the convection term this is irrelevant, since $(u\cdot\nabla)u=\mdiv(u\otimes u)$ has mean zero for divergence-free $u$, but it matters for the damping term in~\autoref{sec:tamed}.
For a Fourier multiplier $M$, we denote its symbol by $M_\xi$, i.e., $\wh{Mu}(\xi)=M_\xi\wh u(\xi)$.
The Stokes operator is $A=-\PP\Delta=-\Delta$ on~$\Hs2$. We write
\[
 \cB(a,b)=\PP[(a\cdot\nabla)b],\qquad\cB(u)=\cB(u,u),
\]
and recall that $(\cB(a,b),b)_2=0$ for divergence-free~$a$.

\begin{lem}[Convection estimates]\label{lem:convection}
There is a constant $C$ such that, for divergence-free mean-zero fields $a,b$ and whenever the right-hand sides are finite,
\begin{align}
 \norm{\cB(a,b)}_{\rH^{-1/2}}&\le C\norm{\nabla a}_2\norm{\nabla b}_2,\label{eq:conv-weak}\\
 \norm{\cB(a,b)}_2&\le C\norm a_{\Hs{5/4}}\norm b_{\Hs{5/4}},\label{eq:conv-sub}\\
 \norm{\cB(a,b)}_2&\le C\norm{\nabla a}_2\norm{\nabla b}_2^{1/2}\norm{\Delta b}_2^{1/2}.\label{eq:conv-strong}
\end{align}
Moreover, for $u\in\Hs2$,
\begin{equation}\label{eq:conv-commutator}
 (\cB(u),\Delta u)_2=-\sum_{j=1}^3\int_{\T^3}(\del_ju\cdot\nabla)u\cdot\del_ju\rd x,
 \qquad\text{hence}\qquad
 |(\cB(u),\Delta u)_2|\le\norm{\nabla u}_\infty\norm{\nabla u}_2^2.
\end{equation}
\end{lem}
\begin{proof}
Since $(a\cdot\nabla)b=\mdiv(a\otimes b)$ has mean zero and $\PP$ is bounded on $\rH^{-1/2}$, the inequality \eqref{eq:conv-weak} follows from the dual Sobolev embedding $\rL^{3/2}\hookrightarrow \rH^{-1/2}$ and $\norm{(a\cdot\nabla)b}_{3/2}\le\norm a_6\norm{\nabla b}_2$. For \eqref{eq:conv-sub}, we use H\"older's inequality with exponents $12$ and $\frac{12}5$ together with $\rH^{5/4}\hookrightarrow \rL^{12}$ and $\rH^{1/4}\hookrightarrow \rL^{12/5}$, and for \eqref{eq:conv-strong} the bound $\norm{(a\cdot\nabla)b}_2\le\norm a_6\norm{\nabla b}_3$. Since $\Delta u$ is divergence-free with mean zero, $(\cB(u),\Delta u)_2=((u\cdot\nabla)u,\Delta u)_2$. For smooth $u$, the identity in \eqref{eq:conv-commutator} then follows by integration by parts, because $((u\cdot\nabla)\del_ju,\del_ju)_2=0$. Both sides are continuous on $\Hs2$ (the right-hand side is bounded by $\norm{\nabla u}_3^3$), so the identity extends to $\Hs2$, and the inequality is immediate.
\end{proof}

\subsection{The one-component geometry}\label{ss:geometry}
For $u\in \rL^2_\sigma$, let $Ku=\PP(u_3\bfe_3)$, and for $\xi\ne0$, set
\[
 b_\xi=P_\xi\bfe_3=\bfe_3-\frac{\xi_3}{|\xi|^2}\,\xi,\qquad\mu_\xi=|b_\xi|^2=\frac{|\xi_h|^2}{|\xi|^2}\in[0,1].
\]

\begin{lem}[Polarization]\label{lem:geometry}
\begin{enumerate}[(i)]
 \item For $u\in \rL^2_\sigma$ and $\xi\ne0$, $\wh{Ku}(\xi)=(b_\xi\otimes b_\xi)\wh u(\xi)=\mu_\xi R_\xi\wh u(\xi)$, where $R_\xi=\mu_\xi^{-1}\,b_\xi\otimes b_\xi$ is the orthogonal projection of $\mathbb C^3$ onto $\mathbb Cb_\xi$ if $\xi_h\ne0$, and $R_\xi=0$ if $\xi_h=0$. We set $S_\xi=P_\xi-R_\xi$.
 \item $K$ is self-adjoint on $\rL^2_\sigma$, $0\le K\le\Id$, and $\ker K=\{u\in \rL^2_\sigma:u_3=0\}$.
 \item Let $\cS$ be the orthogonal projection of $\rL^2_\sigma$ onto $\ker K$, and $\cR=\Id-\cS$. Then $\wh{\cR u}(\xi)=R_\xi\wh u(\xi)$, $\wh{\cS u}(\xi)=S_\xi\wh u(\xi)$, and for $\xi_h\ne0$, it holds that
 \begin{equation}\label{eq:R-symbol}
  R_\xi\wh u(\xi)=\wh u_3(\xi)\begin{pmatrix}-\xi_3\xi_h|\xi_h|^{-2}\\ 1\end{pmatrix}.
 \end{equation}
 In physical variables, $\cR u=(-\nabla_h\Delta_h^{-1}\del_3u_3,\ u_3)$ and $\cS u=(u_h+\nabla_h\Delta_h^{-1}\del_3u_3,\ 0)$, where $\Delta_h^{-1}$ acts on the modes with $\xi_h\ne0$.
 \item $K$, $\cR$, and $\cS$ are Fourier multipliers. They commute with each other and with derivatives, and they have norm at most one on every $\Hs s$ and every space $\cA^r$ introduced below. Moreover, $K\cR=\cR K=K$, $K\cS=\cS K=0$, and
 \[
  (u,Ku)_2=\norm{u_3}_2^2,\qquad(\nabla u,\nabla Ku)_2=\norm{\nabla u_3}_2^2,\qquad(\nabla_hu,\nabla_hKu)_2=\norm{\nabla_hu_3}_2^2.
 \]
\end{enumerate}
\end{lem}
\begin{proof}
(i) Since $\wh u(\xi)\perp\xi$ and $b_\xi$ is real, $\wh u_3(\xi)=\bfe_3\cdot\wh u(\xi)=(P_\xi\bfe_3)\cdot\wh u(\xi)=b_\xi\cdot\wh u(\xi)$, with the bilinear dot product. Hence, $\wh{Ku}(\xi)=P_\xi(\wh u_3(\xi)\bfe_3)=\wh u_3(\xi)\,b_\xi=(b_\xi\otimes b_\xi)\wh u(\xi)$. If $\xi_h=0$, then $b_\xi=0$.

(ii) The symbol $b_\xi\otimes b_\xi$ is symmetric with eigenvalues $0$ and $\mu_\xi\le1$, which gives self-adjointness and $0\le K\le\Id$. Moreover, $Ku=0$ if and only if $\wh u_3(\xi)=0$ for all $\xi$ with $\xi_h\ne0$; on the modes with $\xi_h=0$, incompressibility gives $\xi_3\wh u_3(\xi)=0$, so $\wh u_3$ vanishes there anyway.

(iii) On $\xi^\perp$ the kernel of $b_\xi\otimes b_\xi$ is the orthogonal complement of $b_\xi$, so $\cS$ acts as $S_\xi$ on the mode~$\xi$. Formula \eqref{eq:R-symbol} follows from $R_\xi\wh u(\xi)=\mu_\xi^{-1}(b_\xi\cdot\wh u(\xi))\,b_\xi=\wh u_3(\xi)\,b_\xi/\mu_\xi$ and $b_\xi/\mu_\xi=(-\xi_3\xi_h|\xi_h|^{-2},1)$; the symbol of $-\nabla_h\Delta_h^{-1}\del_3$ is $-\xi_3\xi_h|\xi_h|^{-2}$.

(iv) The first assertions are clear from (i)--(iii). The last three identities hold because $\PP$ is self-adjoint and $\PP u=u$, $\PP\Delta u=\Delta u$, $\PP\Delta_hu=\Delta_hu$.
\end{proof}

\begin{rem}\label{rem:geometry}
The splitting $u=\cR u+\cS u$ is the classical poloidal-toroidal decomposition:
With $\phi=-\Delta_h^{-1}u_3$, one has $\cR u=\nabla\times\nabla\times(\phi\bfe_3)$,
while $\cS u=(w_h,0)$ is horizontal with $\mdiv_hw_h=0$. The point for us is how the
noise interacts with it: By (i), the projected vertical forcing acts \emph{radially}
on each poloidal mode, with strength $\mu_\xi$, and since $\cS K=\cS\cR=0$, the noise in
the models below never acts on the toroidal part. The first fact is what makes
polarized observables such as $X_r$ behave well under It\^o's formula. The second is
the main difficulty, and it is addressed in different ways in
Theorems~\ref{thm:linear} and~\ref{thm:nonlinear}. Note also that the noise acts on
the third coordinate before the incompressibility projection, so the stochastic
increment $Ku\rd W$ generally has horizontal components.
\end{rem}

\begin{exa}[Large accessible data]\label{exa:data}
For $a\in\R$, the field
\[
 u_0^{(a)}=a\begin{pmatrix}-\sin(x_1+x_3)\\-2\sin(x_2+2x_3)\\ \sin(x_1+x_3)+\sin(x_2+2x_3)+\sin(x_1+x_2)\end{pmatrix}
\]
is smooth, mean-zero and divergence-free, and $\cS u_0^{(a)}=0$ by~\eqref{eq:R-symbol}. Its wavevectors $(1,0,1)$, $(0,1,2)$ and $(1,1,0)$ span $\R^3$, and $a$ is arbitrary; so the condition $\cS u_0=0$ forces neither smallness nor two-dimensionality. It is not preserved by the evolution.
\end{exa}

\subsection{Fourier--Wiener norms}\label{ss:wiener}
For $r\in\R$ and $f$ with $\wh f(0)=0$, let $\norm f_{\cA^r}=\sum_{\xi\ne0}|\xi|^r|\wh f(\xi)|$, and denote by $\cA^r$ the space of such $f$ with finite norm. For $u\in \rL^2_\sigma$, we use the polarized version
\begin{equation}\label{eq:Xr}
 X_r(u)=\sum_{\xi\ne0}|\xi|^r\bigl(|R_\xi\wh u(\xi)|+|S_\xi\wh u(\xi)|\bigr)\in[0,\infty].
\end{equation}
Clearly, $X_r(u)=\norm{\cR u}_{\cA^r}+\norm{\cS u}_{\cA^r}$, $\norm u_{\cA^r}\le X_r(u)\le\sqrt2\norm u_{\cA^r}$, $\norm u_\infty\le X_0(u)$, and $\norm{\nabla u}_\infty\le X_1(u)$. We stress that $\Hs1\not\hookrightarrow\cA^0$. The quantities $X_r$ with $r\ge0$ are only used where their finiteness is known.

\begin{lem}\label{lem:wiener}
\begin{enumerate}[(i)]
 \item For $r_0<r<r_1$, $X_r\le X_{r_0}^{(r_1-r)/(r_1-r_0)}X_{r_1}^{(r-r_0)/(r_1-r_0)}$. Moreover, $X_r\le X_{r'}$ for $r\le r'$.
 \item If $s>r+\frac32$, then $X_r(u)\le C_{s,r}\norm u_{\Hs s}$. In particular, there is a constant $C_E$ such that
 \begin{equation}\label{eq:H1-embeddings}
  X_{-1}(u)\le C_E\norm{\nabla u}_2,\qquad X_{-3/4}(u)\le C\norm{\nabla u}_2,\qquad X_0(u)\le C\norm{\Delta u}_2.
 \end{equation}
 \item For $u\in\Hs2$ and $\xi\ne0$, it holds that
 \begin{equation}\label{eq:convolution}
  |\wh{\cB(u)}(\xi)|\le\min\Bigl\{|\xi|\sum_{\eta+\zeta=\xi}|\wh u(\eta)||\wh u(\zeta)|,\ \sum_{\eta+\zeta=\xi}|\wh u(\eta)||\zeta||\wh u(\zeta)|\Bigr\},
 \end{equation}
 and consequently
 \begin{equation}\label{eq:wiener-products}
    \begin{aligned}
   X_{-1}(\cB(u))&\le\sqrt2X_0(u)^2\le\sqrt2X_{-1}(u)X_1(u),\\
   X_0(\cB(u))&\le\sqrt2X_0(u)X_1(u)\le\sqrt2X_{-1}(u)X_2(u).
  \end{aligned}
 \end{equation}
 The first bound in \eqref{eq:convolution} also holds for $\cB(a,b)$ with divergence-free $a$, with $\wh a(\eta)$, $\wh b(\zeta)$ in place of $\wh u(\eta)$, $\wh u(\zeta)$. In particular, $\norm{\cB(a,b)}_{\cA^{-1}}\le\norm a_{\cA^0}\norm b_{\cA^0}$.
\end{enumerate}
\end{lem}
\begin{proof}
(i) is H\"older's inequality for weighted sums. The monotonicity follows from $|\xi|\ge1$. (ii) follows from the Cauchy--Schwarz inequality, $\sum_{\xi\ne0}|\xi|^{2(r-s)}<\infty$, and $|R_\xi z|+|S_\xi z|\le\sqrt2|z|$. For (iii), $\wh{(a\cdot\nabla)b}(\xi)=i\sum_{\eta+\zeta=\xi}(\wh a(\eta)\cdot\zeta)\,\wh b(\zeta)$, and $\wh a(\eta)\cdot\zeta=\wh a(\eta)\cdot\xi$ because $\wh a(\eta)\perp\eta$. Together with $|P_\xi z|\le|z|$ this gives~\eqref{eq:convolution}. Note that the sums converge absolutely for $u\in\Hs2$, since $X_0(u)<\infty$. Summing the first bound with weight $|\xi|^{-1}$ gives $\sum_\xi|\xi|^{-1}|\wh{\cB(u)}(\xi)|\le(\sum_\eta|\wh u(\eta)|)^2\le X_0(u)^2$, and summing the second one gives $\sum_\xi|\wh{\cB(u)}(\xi)|\le X_0(u)X_1(u)$. The factor $\sqrt2$ accounts for the polarization, and the remaining inequalities follow from (i), which gives $X_0^2\le X_{-1}X_1$ and $X_0X_1\le X_{-1}X_2$.
\end{proof}

\subsection{The stochastic models}\label{ss:models}
We fix a stochastic basis $(\Omega,\cF,(\cF_t)_{t\ge0},\bP)$ satisfying the usual conditions, a real $(\cF_t)$-Brownian motion $W$, and $\nu>0$. For $\beta>0$ and $k\ge1$, we consider
\begin{align}
 \rd u&=[\nu\Delta u-\cB(u)]\rd t+\beta Ku\rd W_t,\label{eq:lin}\\
 \rd u&=[\nu\Delta u-\cB(u)]\rd t+\beta\norm{\nabla\cR u}_2^k\,\cR u\rd W_t,\label{eq:nl}
\end{align}
with deterministic initial data $u(0)=u_0\in\Hs1$. Both are of the form \eqref{eq:intro-sns}, with $G_{\rm lin}(u)=\beta Ku$ and $G_{\rm nl}(u)=\beta\norm{\nabla\cR u}_2^k\cR u$. Both coefficients map $\Hs1$ into itself. In both cases, the intensity~$\beta$ is constant in time: There is no cut-off and no switching off of the noise.

\begin{defn}[Maximal solutions]\label{def:solution}
Let $G\in\{G_{\rm lin},G_{\rm nl}\}$ and $u_0\in\Hs1$. A pair $(u,\tau_*)$ is a \emph{local solution} if $\tau_*$ is a stopping time with values in $(0,\infty]$, $u\colon [0,\tau_*)\times\Omega\to\Hs1$ is progressively measurable, and there are stopping times $\tau_n\uparrow\tau_*$ such that a.s.\ $u\in \rC([0,\tau_n];\Hs1)\cap \rL^2(0,\tau_n;\Hs2)$ and, for all $t\in[0,\tau_n]$,
\[
 u(t)=u_0+\int_0^t[\nu\Delta u-\cB(u)]\rd s+\int_0^t\mathbf 1_{[0,\tau_n]}G(u)\rd W_s\qquad\text{in }\rL^2_\sigma.
\]
A local solution $(u,\tau_*)$ is \emph{maximal} if for every local solution $(v,\tau)$ one has a.s.\ $\tau\le\tau_*$ and $u=v$ on~$[0,\tau)$. We call $\tau_*=\tau_*^{(\beta)}$ the lifetime.
\end{defn}

\subsection{Main results}\label{ss:results}
By \autoref{prop:local} below, \eqref{eq:lin} and \eqref{eq:nl} have maximal solutions $(u^{(\beta)},\tau_*^{(\beta)})$, and $\sup_{t<\tau_*}\norm{\nabla u(t)}_2=\infty$ a.s.\ on $\{\tau_*<\infty\}$. Throughout we use the radius
\begin{equation}\label{eq:rho}
 \rho=\frac{\nu}{2\sqrt2}.
\end{equation}

\begin{thm}[Linear one-component noise]\label{thm:linear}
Let $u_0\in\Hs1$ be deterministic, and let $(u^{(\beta)},\tau_*^{(\beta)})$ be the maximal solution of~\eqref{eq:lin}.
\begin{enumerate}[(a)]
 \item If $\cS u_0=0$, then $\displaystyle\lim_{\beta\to\infty}\bP\bigl(\tau_*^{(\beta)}=\infty\bigr)=1$.
 \item In general, $\displaystyle\liminf_{\beta\to\infty}\bP\bigl(\tau_*^{(\beta)}=\infty\bigr)\ge\Bigl(1-\frac{\norm{\cS u_0}_{\cA^{-1}}}\rho\Bigr)_+$.
 \item If $\cS u_0=0$, then $\displaystyle\lim_{\beta\to\infty}\bP\Bigl(\tau_*^{(\beta)}=\infty,\ \int_0^\infty\!\!\int_{\T^3}|u_3^{(\beta)}|^6\rd x\rd t<\infty\Bigr)=1$.
 \item If $\cS u_0=0$ and in addition $u_0\in\cA^0$, then for every integer $m\ge1$ and every $\eta>0$,
 \[
  \lim_{\beta\to\infty}\bP\Bigl(\tau_*^{(\beta)}=\infty,\ \int_0^\infty\!\!\int_{\T^3}|u_3^{(\beta)}|^{2m}\rd x\rd t\le\eta\Bigr)=1.
 \]
\end{enumerate}
\end{thm}

\begin{rem}\label{rem:linear}
(1) By \autoref{exa:data}, part (a) covers arbitrarily large, genuinely three-dimensional data. The intensity $\beta_0(\eps)$ beyond which $\bP(\tau_*^{(\beta)}=\infty)\ge1-\eps$ depends on $u_0$ through the distribution of~$\cR u_0$ over the values of $\mu_\xi$ (see \autoref{rem:contraction}). 

(2) The restriction $\cS u_0=0$ in (a) is due to the fact that the linear noise never acts on $\cS u$, whereas the argument needs the solution to enter a small ball after a short time; part (b) quantifies this. We do not know whether (a) holds for arbitrary data. For every $\beta>0$, however, finite-time blow-up can only occur on the event where $\int_0^{\tau_*}\norm{u_3}_6^6\rd t=\infty$ (\autoref{prop:lin-criterion}). For the nonlinear noise, this restriction is removed by means of a reference flow, see~\autoref{thm:nonlinear}.

(3) Part (c) asserts finiteness only. The smallness in (d) uses the additional assumption $u_0\in\cA^0$, which serves to control the integrals near $t=0$; we do not know whether (d) holds for arbitrary $u_0\in\Hs1$ with $\cS u_0=0$ (cf.~\autoref{rem:powers}).
\end{rem}

\begin{thm}[Nonlinear one-component noise]\label{thm:nonlinear}
Let $k\ge4$, let $u_0\in\Hs1$ be deterministic, and let $(u^{(\beta)},\tau_*^{(\beta)})$ be the maximal solution of~\eqref{eq:nl}. Then for every $\eta>0$, it holds that
\begin{equation}\label{eq:main-nl}
 \lim_{\beta\to\infty}\bP\Bigl(\tau_*^{(\beta)}=\infty,\ \int_0^\infty\!\!\int_{\T^3}|u_3^{(\beta)}|^6\rd x\rd t\le\eta\Bigr)=1.
\end{equation}
Moreover, there are constants $C,c>0$, depending on $u_0$, $\nu$ and $k$, such that
\[
 \bP\bigl(\tau_*^{(\beta)}<\infty\bigr)\le C\exp\bigl(-c\,\beta^{2/(2k+1)}\bigr)\qquad\text{for all }\beta\ge1.
\]
\end{thm}

\begin{rem}\label{rem:nonlinear}
(1) The coefficient in \eqref{eq:nl} is radial on the accessible part: It points along $\cR u$, and its size only depends on $\norm{\nabla\cR u}_2$. Since $\cR u$ is determined by $u_3$, the noise is driven by the vertical component alone. For $k=4$, it is quintic in the amplitude, but it is not the local coefficient $u_3^5\bfe_3$, and $\cR$ cannot be replaced by $K$ in our argument.

(2) No smallness or structural condition on $u_0$ is imposed. The thresholds depend on the datum through the choice of a reference flow (\autoref{ss:choice}). 

(3) For the scalar model $\rd X=cX^3\rd W$, It\^o's formula gives $\rd\log|X|=-\frac{c^2}2X^4\rd t+cX^2\rd W$, while the quintic damping $\dot X=-\alpha X^5$ gives $\frac{\rd}{\rd t}\log|X|=-\alpha X^4$. In this sense, cubic It\^o noise is the degree-matched analogue of quintic damping. In \autoref{thm:nonlinear} the noise has the higher degree $k+1\ge5$: The logarithmic correction $\frac{\beta^2}2h_k(r)\sim\frac{\beta^2}2r^{2k}$, for large $r=\norm{\nabla\cR u}_2$, has to absorb a residual of order $r^8$ produced by the convection terms (\autoref{lem:absorb}).
\end{rem}

\begin{rem}\label{rem:ito}
Both theorems concern It\^o noise, and the It\^o correction drives both proofs: The Fourier multipliers of the linear noise flow are $\exp(\beta\mu_\xi W_t-\frac12\beta^2\mu_\xi^2t)$, which tend to zero almost surely as $t\to\infty$ (for $\mu_\xi>0$) although their expectation equals one, and the drift in \eqref{eq:intro-log} carries the term~$-\frac12f(u)^2$. Rewriting the noise in Stratonovich form changes the drift. For the scalar model equation, $cX^3\rd W=cX^3\circ\rd W-\frac{3c^2}2X^5\rd t$. Thus, the Stratonovich form of the same equation contains an explicit damping term, and the Stratonovich equation without this term is a different equation, see \cite[Section~6]{BMX:25} for a discussion in the Euler setting.
\end{rem}

\section{Local well-posedness via stochastic maximal regularity}\label{sec:local}

\subsection{The abstract result}\label{ss:AV}
We use the critical variational setting of Agresti and Veraar \cite{AV:24} in the form presented in \cite[Sections~4--6]{AV:25}; it corresponds to the case $p=2$, $\kappa=0$ of the $\rL^p$-theory developed in~\cite{AV:22a,AV:22b}. Let
\[
 V=\Hs2,\qquad H=\Hs1,\qquad V^*\cong\Hs0=\rL^2_\sigma,
\]
with the norms $\norm v_V=\norm{\Delta v}_2$, $\norm f_H=\norm{\nabla f}_2$ and $\norm f_{V^*}=\norm f_2$, and with the duality pairing $\ip vf=(-\Delta v,f)_2$, which identifies $\rL^2_\sigma$ isometrically with the dual of~$V$, since $-\Delta:\Hs2\to \rL^2_\sigma$ is an isometric isomorphism. As $\ip vf=(v,f)_H$ for $f\in H$, $(V,H,V^*)$ is a Gelfand triple. Moreover, $[V^*,V]_\theta=\Hs{2\theta}$, and $V_{\theta,1}=(V^*,V)_{\theta,1}\hookrightarrow\Hs{2\theta}$ for $\theta\in(0,1)$. For a real Brownian motion~$W$, we consider
\begin{equation}\label{eq:abstract}
 \rd u+A_0u\rd t=F(u)\rd t+G(u)\rd W_t,\qquad u(0)=u_0.
\end{equation}

\begin{thm}[{\cite[Theorems~4.7, 5.1, 5.2 and~6.2]{AV:25}}]\label{thm:AV}
Assume the following.
\begin{enumerate}[(H1)]
 \item $A_0\in\mathscr L(V,V^*)$, and there are $\theta>0$ and $M\ge0$ such that $\ip v{A_0v}\ge\theta\norm v_V^2-M\norm v_H^2$ for all $v\in V$.
 \item $F\colon V\to V^*$ and $G\colon V\to H$ are measurable, and for every $n\ge1$, there is $L_n$ such that for all $u,v\in V$ with $\norm u_H,\norm v_H\le n$,
 \[
  \norm{F(u)-F(v)}_{V^*}+\norm{G(u)-G(v)}_H\le L_n\sum_{j=1}^m\bigl(1+\norm u_{V_{\beta_j,1}}^{\rho_j}+\norm v_{V_{\beta_j,1}}^{\rho_j}\bigr)\norm{u-v}_{V_{\beta_j,1}},
 \]
 where $\beta_j\in(\frac12,1)$, $\rho_j\ge0$ and $(2\beta_j-1)(\rho_j+1)\le1$ for $j=1,\dots,m$. 
\end{enumerate}
Then for every $u_0\in \rL^0_{\cF_0}(\Omega;H)$ there is a unique maximal solution $(u,\tau_*)$ of \eqref{eq:abstract} in the sense of \autoref{def:solution}, where the identity holds in $V^*$. Moreover, $\tau_*>0$ a.s., and the following assertions hold.
\begin{enumerate}[(a)]
 \item $\bP\bigl(\tau_*<\infty,\ \sup_{t<\tau_*}\norm{u(t)}_H+\int_0^{\tau_*}\norm{u(t)}_V^2\rd t<\infty\bigr)=0$.
 \item If $(2\beta_j-1)(\rho_j+1)<1$ for all $j$, then $\bP\bigl(\tau_*<\infty,\ \sup_{t<\tau_*}\norm{u(t)}_H<\infty\bigr)=0$.
 \item If $(v,\tau)$ is the maximal solution with initial value $v_0\in \rL^0_{\cF_0}(\Omega;H)$, then a.s.\ on $\{u_0=v_0\}$, one has $\tau=\tau_*$ and $u=v$ on $[0,\tau_*)$.
\end{enumerate}
\end{thm}

The conditions of \autoref{thm:AV} are the case $p=2$, $\kappa=0$, $B=0$ of \cite[Assumption~4.1]{AV:25}, written in the variational form of \cite[Assumption~6.1]{AV:25}. By \cite[Theorem~3.17]{AV:25}, (H1) implies stochastic maximal~$\rL^2$-regularity of $(A_0,0)$. 
The theorem includes deterministic equations ($G=0$), which we use in \autoref{sec:tamed}; in that case one may alternatively use the deterministic theory of critical spaces~\cite{PSW:18}.

\subsection{Application to the two models}\label{ss:application}
\begin{prop}[Local well-posedness]\label{prop:local}
Let $\beta>0$, $k\ge1$ and $u_0\in\Hs1$. Then \eqref{eq:lin} and~\eqref{eq:nl} are of the form \eqref{eq:abstract} with
\[
 A_0=\nu A,\qquad F(u)=-\cB(u),\qquad G\in\{G_{\rm lin},G_{\rm nl}\},
\]
where $A$ and $\cB$ are as in \autoref{ss:spaces}, and $G_{\rm lin}$ and $G_{\rm nl}$ are as in \autoref{ss:models}.
Both models satisfy (H1)--(H2) with $m=2$, $(\beta_1,\rho_1)=(\frac58,1)$ and $(\beta_2,\rho_2)=(\frac58,0)$.
In particular, the setting is subcritical. Consequently, each of them has a unique maximal solution $(u,\tau_*)$ in the sense of \autoref{def:solution}, and
\begin{equation}\label{eq:blowup}
 \bP\Bigl(\tau_*<\infty,\ \sup_{t<\tau_*}\norm{\nabla u(t)}_2<\infty\Bigr)=0.
\end{equation}
\end{prop}
\begin{proof}
(H1) holds with $\theta=\nu$ and $M=0$, since $\ip v{\nu Av}=\nu\norm{\Delta v}_2^2$. 
For the drift, bilinearity and~\eqref{eq:conv-sub} give
\[
 \norm{\cB(u)-\cB(v)}_2\le\norm{\cB(u-v,u)}_2+\norm{\cB(v,u-v)}_2\le C\bigl(\norm u_{\Hs{5/4}}+\norm v_{\Hs{5/4}}\bigr)\norm{u-v}_{\Hs{5/4}},
\]
and $V_{5/8,1}\hookrightarrow\Hs{5/4}$. This is the first term in (H2), and $(2\beta_1-1)(\rho_1+1)=\frac12<1$. For the noise, \autoref{lem:geometry}(iv) gives $\norm{G_{\rm lin}(u)-G_{\rm lin}(v)}_H\le\beta\norm{\nabla(u-v)}_2$. 
For $G_{\rm nl}(u)=\beta\varphi(\cR u)$ with $\varphi(x)=\norm x_H^kx$, the identity $\varphi(a)-\varphi(b)=\norm a_H^k(a-b)+(\norm a_H^k-\norm b_H^k)\,b$ and $|s^k-t^k|\le k\max\{s,t\}^{k-1}|s-t|$ give $\norm{\varphi(a)-\varphi(b)}_H\le(k+1)\max\{\norm a_H,\norm b_H\}^k\norm{a-b}_H$. Since $\cR$ is a contraction on $H$,
\[
 \norm{G_{\rm nl}(u)-G_{\rm nl}(v)}_H\le\beta(k+1)n^k\norm{\nabla(u-v)}_2
 \qquad\text{if }\norm u_H,\norm v_H\le n.
\]
Since $\norm\cdot_H\le C\norm\cdot_{V_{5/8,1}}$, both noise coefficients satisfy the second term in (H2), and here $(2\beta_2-1)(\rho_2+1)=\frac14<1$. \autoref{thm:AV} provides the maximal solution, and part (b) gives~\eqref{eq:blowup}.
\end{proof}

\subsection{Energy identities}\label{ss:energy}
\begin{lem}[$\Hs1$-balances]\label{lem:energy}
Let $(u,\tau_*)$ be the maximal solution of \eqref{eq:lin} or~\eqref{eq:nl}. Then,
\begin{equation}\label{eq:energy-H1}
\begin{split}
 \norm{\nabla u(t)}_2^2+2\nu\int_0^t\norm{\Delta u}_2^2\rd s
 &=\norm{\nabla u_0}_2^2+\int_0^t\bigl[2(\cB(u),\Delta u)_2+\norm{\nabla G(u)}_2^2\bigr]\rd s\\
 &\quad+2\int_0^t(\nabla G(u),\nabla u)_2\rd W_s
\end{split}
\end{equation}
a.s., for all $t<\tau_*$.
For \eqref{eq:lin}, $\norm{\nabla G(u)}_2^2=\beta^2\norm{\nabla Ku}_2^2\le\beta^2\norm{\nabla u}_2^2$ and $(\nabla G(u),\nabla u)_2=\beta\norm{\nabla u_3}_2^2$. For \eqref{eq:nl}, the quantities $v=\cR u$ and $r=\norm{\nabla v}_2$ satisfy
\begin{equation}\label{eq:energy-v}
 r(t)^2+2\nu\int_0^t\norm{\Delta v}_2^2\rd s=r(0)^2+\int_0^t\bigl[2(\cB(u),\Delta v)_2+\beta^2r^{2k+2}\bigr]\rd s+2\beta\int_0^tr^{k+2}\rd W_s.
\end{equation}
\end{lem}
\begin{proof}
Let $\tau_n\uparrow\tau_*$ be localizing times such that $\norm{\nabla u}_2\le n$ and
$\int_0^{\tau_n}\norm{\Delta u}_2^2\rd s\le n$ on $[0,\tau_n]$, and~$\tau_n\le n$. We use
the It\^o formula in the Gelfand triple $(V,H,V^*)$, \cite[Lemma~A.5]{AV:25}. Since
$U=\R$, we identify $\rL_2(\R,H)$ with $H$ via $\Psi\mapsto\Psi(1)$; then
$\Psi^*u=(\Psi,u)_H=(\nabla\Psi,\nabla u)_2$, and the It\^o correction is
$\norm{\Psi}_{\rL_2(\R,H)}^2=\norm{\nabla\Psi}_2^2$.

\emph{Integrability.} Set $\Phi=\nu\Delta u-\cB(u)$ and $\Psi=G(u)$. On $[0,\tau_n]$, we
have $\Delta u\in \rL^2(0,\tau_n;\rL^2_\sigma)$, and by~\eqref{eq:conv-strong},
\[
 \int_0^{\tau_n}\norm{\cB(u)}_2^2\rd s\le Cn^3\int_0^{\tau_n}\norm{\Delta u}_2\rd s
 \le Cn^{4},
\]
so $\Phi\in \rL^2(0,\tau_n;V^*)$. Moreover, $\Psi$ is continuous in time with values in $H$,
since $u\in \rC([0,\tau_n];H)$ and~$G$ is locally Lipschitz on $H$
(\autoref{prop:local}), and it is bounded there:
In fact, $\norm{\nabla G_{\rm lin}(u)}_2\le\beta n$ and
$\norm{\nabla G_{\rm nl}(u)}_2=\beta\norm{\nabla\cR u}_2^{k+1}\le\beta n^{k+1}$.
Both $\Phi$ and $\Psi$ are progressively measurable.

\emph{Localization.} Note that \cite[Lemma~A.5]{AV:25} is formulated on deterministic intervals for processes in~$\rL^2(0,T;V)$. We therefore extend $u$ beyond $\tau_n$ by the heat flow: Let $\bar u=u$ on
$[0,\tau_n]$, and let $\bar u$ solve $\del_t\bar u=\nu\Delta\bar u$ on $(\tau_n,T]$ with
$\bar u(\tau_n)=u(\tau_n)\in H$. Then $\bar u\in \rC([0,T];H)\cap \rL^2(0,T;V)$ satisfies the
assumptions of \cite[Lemma~A.5]{AV:25} with $\bar\Phi=\mathbf 1_{[0,\tau_n]}\Phi
+\mathbf 1_{(\tau_n,T]}\nu\Delta\bar u$ and $\bar\Psi=\mathbf 1_{[0,\tau_n]}\Psi$.
Restricting the resulting identity to $t\le\tau_n$ and using
$\ip u\Phi=(-\Delta u,\nu\Delta u-\cB(u))_2=-\nu\norm{\Delta u}_2^2+(\cB(u),\Delta u)_2$
gives \eqref{eq:energy-H1} on $[0,\tau_n]$; letting $n\to\infty$ yields it for all
$t<\tau_*$.

\emph{The noise terms.} For \eqref{eq:lin}, the formulas in the statement follow from \autoref{lem:geometry}(iv). For~\eqref{eq:nl}, since $\cR$ is a
self-adjoint idempotent commuting with $\nabla$, we have
$(\nabla\cR u,\nabla u)_2=r^2$, hence $(\nabla G_{\rm nl}(u),\nabla u)_2=\beta r^{k+2}$
and $\norm{\nabla G_{\rm nl}(u)}_2^2=\beta^2r^{2k+2}$.

\emph{The accessible part.} Applying $\cR$ to the equation, and using that $\cR$ is
bounded on $V$, $H$ and $V^*$, commutes with $\Delta$ and satisfies
$\cR G_{\rm nl}(u)=G_{\rm nl}(u)=\beta r^kv$, we see that $v=\cR u$ solves
\[
 v(t)=\cR u_0+\int_0^t\bigl(\nu\Delta v-\cR\cB(u)\bigr)\rd s+\int_0^t\beta r^kv\rd W_s
\]
in the same class. The same argument, together with
$(-\Delta v,\cR\cB(u))_2=(\cB(u),-\Delta v)_2$, which is in turn due to the self-adjointness of $\cR$ and
$\cR\Delta v=\Delta v$, gives~\eqref{eq:energy-v}.
\end{proof}

\subsection{A continuation criterion for the linear noise}\label{ss:continuation}
\begin{lem}\label{lem:continuation}
Let $(u,\tau_*)$ be the maximal solution of~\eqref{eq:lin}. Then
\begin{equation}\label{eq:continuation}
 \bP\Bigl(\tau_*<\infty,\ \int_0^{\tau_*}X_1(u(t))\rd t<\infty\Bigr)=0.
\end{equation}
\end{lem}
\begin{proof}
Let $E(t)=\norm{\nabla u(t)}_2^2$, and for $t<\tau_*$ define
\[
 a(t)=\frac{2(\cB(u),\Delta u)_2+\beta^2\norm{\nabla Ku}_2^2}{E(t)},\qquad b(t)=\frac{2\beta\norm{\nabla u_3}_2^2}{E(t)}
\]
on $\{E(t)>0\}$, and $a(t)=b(t)=0$ otherwise. By \eqref{eq:conv-commutator}, $\norm{\nabla u}_\infty\le X_1(u)$ and \autoref{lem:geometry}(iv),
\[
 |a|\le2X_1(u)+\beta^2,\qquad|b|\le2\beta,
\]
and \eqref{eq:energy-H1} reads $\rd E=(aE-2\nu\norm{\Delta u}_2^2)\rd t+bE\rd W$.
Since $X_1(u)$ need not be integrable in time a priori, we localize with
\[
 \vartheta_m=\inf\Bigl\{t<\tau_*:\int_0^tX_1(u)\rd s\ge m\Bigr\},\qquad m\in\N,
\]
and work on $[0,\vartheta_m\wedge\tau_n]$, where $\int|a|\rd s\le2m+\beta^2n$. There
the integrating factor
\[
 \Gamma_t=\exp\Bigl(-\int_0^ta\rd s+\frac12\int_0^tb^2\rd s-\int_0^tb\rd W_s\Bigr)
\]
is well defined. It satisfies $\rd\Gamma=\Gamma[(b^2-a)\rd t-b\rd W]$, so that It\^o's
product rule gives $\rd(\Gamma E)=-2\nu\Gamma\norm{\Delta u}_2^2\rd t$, i.e.,
\begin{equation}\label{eq:integrating-factor}
 \Gamma_tE(t)+2\nu\int_0^t\Gamma_s\norm{\Delta u(s)}_2^2\rd s=E(0),
 \qquad t\le\vartheta_m\wedge\tau_n .
\end{equation}
On the event in \eqref{eq:continuation} we have $\vartheta_m\ge\tau_*$ for all large
$m$, so \eqref{eq:integrating-factor} holds for all $t<\tau_*$ there. Moreover,
$\int_0^{\tau_*}|a|\rd s<\infty$ and $\int_0^{\tau_*}b^2\rd s\le4\beta^2\tau_*<\infty$
on this event. 
Since a continuous local martingale converges on the event where its quadratic variation is finite, $\int_0^tb\rd W$ has a finite limit as $t\uparrow\tau_*$, a.s.\ on this event. Thus, $\inf_{t<\tau_*}\Gamma_t>0$ there, and \eqref{eq:integrating-factor} yields $\sup_{t<\tau_*}E(t)<\infty$. By \eqref{eq:blowup}, the event in~\eqref{eq:continuation} is a null set.
\end{proof}

The criterion does not require $X_1(u_0)<\infty$.
Indeed, in \autoref{sec:linear}, the integral $\int X_1(u)\rd t$ is produced by the viscous term in the evolution of~$X_{-1}(u)$.

\subsection{The linear noise as a random change of variables}\label{ss:random-pde}
For the linear model, we use the solution operator of $\rd z=\beta Kz\rd W$,
\begin{equation}\label{eq:Qbeta}
 Q_\beta(t)=\exp\bigl(\beta KW_t-\tfrac12\beta^2tK^2\bigr),\qquad(Q_\beta(t))_\xi=S_\xi+q_\xi(t)R_\xi,\quad q_\xi(t)=\exp\bigl(\beta\mu_\xi W_t-\tfrac12\beta^2\mu_\xi^2t\bigr),
\end{equation}
where we used $K_\xi=\mu_\xi R_\xi$ and $K_\xi^2=\mu_\xi^2R_\xi$. By It\^o's formula, $\rd q_\xi=\beta\mu_\xi q_\xi\rd W$ and $\rd(q_\xi^{-1})=q_\xi^{-1}(-\beta\mu_\xi\rd W+\beta^2\mu_\xi^2\rd t)$. The deterministic part of the argument is conveniently formulated for general operator paths of this form.

\begin{defn}\label{def:paths}
For $L\ge1$ and $T_0>0$, let $\cQ_L(T_0)$ be the set of families $(Q(t))_{t\in[0,T_0]}$ of Fourier multipliers $Q(t)_\xi=S_\xi+q_\xi(t)R_\xi$ with continuous functions $q_\xi\colon[0,T_0]\to[L^{-1},L]$.
\end{defn}

For $Q\in\cQ_L(T_0)$ the operators $Q(t)^{\pm1}$ commute with $\Delta$, $\PP$, $\cR$ and~$\cS$. They have norm at most $L$ on every $\Hs s$, $s\in\R$, and on the divergence-free part of every $\cA^r$, since $|S_\xi z+qR_\xi z|\le\max\{1,q\}|z|$ for $z\in\xi^\perp$, and $X_r(Q(t)^{\pm1}u)\le L\,X_r(u)$. For $Q\in\cQ_L(T_0)$ and $y_0\in\Hs1$, we consider the integral equation
\begin{equation}\label{eq:mild}
 y(t)=e^{\nu t\Delta}y_0-\int_0^te^{\nu(t-s)\Delta}N_Q(s,y(s))\rd s,\qquad N_Q(s,y)=Q(s)^{-1}\cB(Q(s)y).
\end{equation}

\begin{lem}[Pathwise mild solutions]\label{lem:mild}
Let $L\ge1$, $T_0>0$, $Q\in\cQ_L(T_0)$, $y_0\in\Hs1$, $M_1=\max\{1,\norm{\nabla y_0}_2\}$, and
\begin{equation}\label{eq:tL}
 t_L=\min\bigl\{1,\ c_*\,\nu^3L^{-12}M_1^{-4}\bigr\},
\end{equation}
where $c_*>0$ is a numerical constant.
\begin{enumerate}[(a)]
 \item For $T\le\min\{T_0,t_L\}$, there is a solution $y\in\rC([0,T];\Hs1)$ of \eqref{eq:mild}, and it satisfies $\sup_{t\le T}\norm{\nabla y(t)}_2\le2M_1$.
 \item For every $T'\le T_0$, \eqref{eq:mild} has at most one solution in~$\rC([0,T'];\Hs1)$.
 \item The solution in (a) satisfies, for $t\le T$,
 \begin{equation}\label{eq:displacement}
  \norm{\nabla(y(t)-y_0)}_2\le\omega_{y_0}(t)+C_1L^3\nu^{-3/4}t^{1/4}M_1^2,\qquad\omega_{y_0}(t)=\sup_{s\le t}\norm{\nabla(e^{\nu s\Delta}-\Id)y_0}_2.
 \end{equation}
 \item If, in addition, $y_0\in\cA^0$ and $M_0=\max\{M_1,\norm{y_0}_{\cA^0}\}$, let $t_L'=\min\{t_L,c_*'\nu L^{-6}M_0^{-2}\}$. For $T\le\min\{T_0,t_L'\}$ the solution in (a) belongs to $\rC([0,T];\cA^0)$, $\sup_{t\le T}\norm{y(t)}_{\cA^0}\le2M_0$, and
 \begin{equation}\label{eq:displacement-A0}
  \norm{y(t)-y_0}_{\cA^0}\le\sup_{s\le t}\norm{(e^{\nu s\Delta}-\Id)y_0}_{\cA^0}+C_1'L^3\nu^{-1/2}t^{1/2}M_0^2,\qquad t\le T.
 \end{equation}
\end{enumerate}
The constants $c_*,c_*',C_1,C_1'$ do not depend on $Q$, $L$, $T_0$, $\nu$, or~$y_0$.
\end{lem}
\begin{proof}
We use the heat estimates for mean-zero fields
\[
 \norm{\nabla e^{\nu t\Delta}f}_2\le C_h(\nu t)^{-3/4}\norm f_{\rH^{-1/2}},\qquad\norm{e^{\nu t\Delta}f}_{\cA^0}\le(2e\nu t)^{-1/2}\norm f_{\cA^{-1}},
\]
which follow from $\sup_{\lambda\ge1}\lambda^{3/2}e^{-\nu t\lambda^2}\le C_h(\nu t)^{-3/4}$ and $\sup_{\lambda>0}\lambda e^{-\nu t\lambda^2}=(2e\nu t)^{-1/2}$. By \eqref{eq:conv-weak}, bilinearity and the bounds on $Q^{\pm1}$,
\[
 \norm{N_Q(s,y)-N_Q(s,z)}_{\rH^{-1/2}}\le C_0L^3\bigl(\norm{\nabla y}_2+\norm{\nabla z}_2\bigr)\norm{\nabla(y-z)}_2,\qquad N_Q(s,0)=0.
\]
Denote by $\Phi(y)$ the right-hand side of~\eqref{eq:mild}. Since $\int_0^t(t-s)^{-3/4}\rd s=4t^{1/4}$ and $\norm{\nabla e^{\nu t\Delta}y_0}_2\le\norm{\nabla y_0}_2$, on the set $Z_T=\{y\in \rC([0,T];\Hs1):\sup_t\norm{\nabla y(t)}_2\le2M_1\}$, we have
\begin{align*}
 \sup_{t\le T}\norm{\nabla\Phi(y)(t)}_2&\le M_1+16C_hC_0L^3\nu^{-3/4}T^{1/4}M_1^2,\\
 \sup_{t\le T}\norm{\nabla(\Phi(y)-\Phi(z))(t)}_2&\le16C_hC_0L^3\nu^{-3/4}T^{1/4}M_1\sup_{t\le T}\norm{\nabla(y-z)(t)}_2.
\end{align*}
If $16C_hC_0L^3\nu^{-3/4}T^{1/4}M_1\le\frac12$, which holds for $T\le t_L$ with $c_*=(32C_hC_0)^{-4}$, then $\Phi$ is a contraction of $Z_T$ (continuity in time of $\Phi(y)$ is standard), and (a) follows.

For (b), let $y,z\in \rC([0,T'];\Hs1)$ be two solutions with $\norm{\nabla y}_2,\norm{\nabla z}_2\le R$. Then $\delta(t)=\norm{\nabla(y-z)(t)}_2$ satisfies $\delta(t)\le2C_hC_0L^3R\nu^{-3/4}\int_0^t(t-s)^{-3/4}\delta(s)\rd s$. On an interval $[0,h]$ with $8C_hC_0L^3R\nu^{-3/4}h^{1/4}<1$ this gives $\delta=0$, and iterating over consecutive intervals of length $h$ yields $\delta\equiv0$ and thus shows~(b).

The estimate \eqref{eq:displacement} follows directly from \eqref{eq:mild}, the heat estimate and $\sup\norm{\nabla y}_2\le2M_1$, with $C_1=16C_hC_0$.

For (d) we use $\norm{\cB(a,b)}_{\cA^{-1}}\le\norm a_{\cA^0}\norm b_{\cA^0}$ from \autoref{lem:wiener}(iii), which gives
\[
 \norm{N_Q(s,y)-N_Q(s,z)}_{\cA^{-1}}\le L^3\bigl(\norm y_{\cA^0}+\norm z_{\cA^0}\bigr)\norm{y-z}_{\cA^0}.
\]
Together with $\int_0^t(t-s)^{-1/2}\rd s=2t^{1/2}$, the map $\Phi$ is a contraction of
\[
 \bigl\{y\in Z_T\cap \rC([0,T];\cA^0):\sup_t\norm{y(t)}_{\cA^0}\le2M_0\bigr\},\qquad\text{with the metric }\sup_t\bigl(\norm{\nabla(y-z)}_2+\norm{y-z}_{\cA^0}\bigr),
\]
provided $T\le t_L'$ with a suitable~$c_*'$. By (b), its fixed point is the solution from (a). The displacement estimate follows as before, and $(e^{\nu s\Delta}-\Id)y_0\to0$ in $\cA^0$ as $s\downarrow0$ by dominated convergence.
\end{proof}

\begin{lem}[Random change of variables]\label{lem:transformation}
Let $(u,\tau_*)$ be the maximal solution of \eqref{eq:lin}, and let $y=Q_\beta^{-1}u$ on $[0,\tau_*)$. Then a.s.\ $y\in \rC([0,\tau_*);\Hs1)$, and for every $T<\tau_*$, the restriction $y|_{[0,T]}$ solves \eqref{eq:mild} with $Q=Q_\beta$ and $y_0=u_0$.
\end{lem}
\begin{proof}
Fix $\xi\ne0$. Testing the equation in \autoref{def:solution} with $e^{-i\xi\cdot x}$ shows that the Fourier coefficients satisfy, in integral form on $[0,\tau_*)$,
\[
 \rd\wh u(\xi)=-\bigl[\nu|\xi|^2\wh u(\xi)+\wh{\cB(u)}(\xi)\bigr]\rd t+\beta\mu_\xi R_\xi\wh u(\xi)\rd W_t.
\]
By It\^o's product rule and the formula for $\rd(q_\xi^{-1})$,
\[
 \rd\bigl(q_\xi^{-1}R_\xi\wh u(\xi)\bigr)=-q_\xi^{-1}R_\xi\bigl[\nu|\xi|^2\wh u(\xi)+\wh{\cB(u)}(\xi)\bigr]\rd t,
\]
because the two stochastic integrals cancel, and so do the It\^o correction $\beta^2\mu_\xi^2q_\xi^{-1}R_\xi\wh u(\xi)\rd t$ and the covariation term $-\beta^2\mu_\xi^2q_\xi^{-1}R_\xi\wh u(\xi)\rd t$. Moreover, $\rd(S_\xi\wh u(\xi))=-S_\xi[\nu|\xi|^2\wh u(\xi)+\wh{\cB(u)}(\xi)]\rd t$, since $S_\xi K_\xi=0$. Hence, $\wh y(\xi)$ is absolutely continuous on compact subintervals of $[0,\tau_*)$, and, as $\cB(u)=\cB(Q_\beta y)$,
\[
 \frac{\rd}{\rd t}\wh y(\xi)=-\nu|\xi|^2\wh y(\xi)-\wh{N_{Q_\beta}(t,y)}(\xi).
\]
Variation of constants on each mode gives \eqref{eq:mild} mode by mode, hence in~$\Hs1$. Finally, $t\mapsto Q_\beta(t)^{-1}$ is continuous in operator norm (uniformly in $\xi$, since $\mu_\xi\in[0,1]$), so $y\in \rC([0,\tau_*);\Hs1)$.
\end{proof}

\begin{cor}[The initial layer]\label{cor:layer}
Let $(u,\tau_*)$ be the maximal solution of \eqref{eq:lin}, let $L\ge1$, put $M_1=\max\{1,\norm{\nabla u_0}_2\}$, and let $T\le t_L$ be deterministic, where $t_L$ is given by~\eqref{eq:tL}. On the event
\[
 \Omega_{L,T}=\Bigl\{\sup_{t\le T}\sup_{\xi\ne0}\max\{q_\xi(t),q_\xi(t)^{-1}\}\le L\Bigr\}\in\cF_T,
\]
we have a.s.\ $\tau_*>T$, $\sup_{t\le T}\norm{\nabla u(t)}_2\le2LM_1$, and
\begin{equation}\label{eq:comparison}
 \norm{\nabla(u(t)-Q_\beta(t)u_0)}_2\le L\,\omega_{u_0}(t)+C_1L^4\nu^{-3/4}t^{1/4}M_1^2,\qquad t\le T.
\end{equation}
If $u_0\in\cA^0$ and $T\le t_L'$, then in addition, with $M_0=\max\{M_1,\norm{u_0}_{\cA^0}\}$,
\begin{equation}\label{eq:comparison-A0}
  \begin{aligned}
  &\sup_{t\le T}X_0(u(t))\le4LM_0,\\
  &\norm{u(t)-Q_\beta(t)u_0}_{\cA^0}\le L\sup_{s\le t}\norm{(e^{\nu s\Delta}-\Id)u_0}_{\cA^0}+C_1'L^4\nu^{-1/2}t^{1/2}M_0^2 .
 \end{aligned}
\end{equation}
\end{cor}
\begin{proof}
On $\Omega_{L,T}$ the restriction of $Q_\beta$ to $[0,T]$ belongs to $\cQ_L(T)$; let $\tilde y$ be the solution from \autoref{lem:mild}(a). By \autoref{lem:transformation} and \autoref{lem:mild}(b), $Q_\beta^{-1}u=\tilde y$ on $[0,\tau_*\wedge T)$. If $\tau_*\le T$, then $\sup_{t<\tau_*}\norm{\nabla u(t)}_2\le2LM_1$, and by \eqref{eq:blowup} this happens only on a null set. Thus $\tau_*>T$ a.s.\ on $\Omega_{L,T}$ and $u=Q_\beta\tilde y$ on $[0,T]$, and \eqref{eq:comparison} follows from \eqref{eq:displacement} because $u-Q_\beta u_0=Q_\beta(\tilde y-u_0)$. The statements in $\cA^0$ follow in the same way from \autoref{lem:mild}(d), using $X_0(Q_\beta\tilde y)\le L\,X_0(\tilde y)\le\sqrt2L\norm{\tilde y}_{\cA^0}$.
\end{proof}

\subsection{Persistence of Wiener-algebra regularity}\label{ss:persistence}
Note that $t\mapsto X_r(u(t))$ is lower semicontinuous on $[0,\tau_*)$ for every $r$, as a supremum of finite sums of continuous functions. In particular, for an adapted lower semicontinuous process, the first time it exceeds a given level is a stopping time, and the process does not exceed this level before (and at) that time.

\begin{lem}\label{lem:persistence}
Let $(u,\tau_*)$ be the maximal solution of \eqref{eq:lin}, let $t_0\ge0$ be deterministic, and let~$\Gamma\in\cF_{t_0}$ with $\Gamma\subset\{t_0<\tau_*,\ X_0(u(t_0))<\infty\}$. Then a.s.\ on $\Gamma$ there is $\tau_A\in(t_0,\tau_*]$ such that $u\in \rC([t_0,\tau_A);\cA^0)$, and $\limsup_{t\uparrow\tau_A}X_0(u(t))=\infty$ if $\tau_A<\tau_*$.
\end{lem}
\begin{proof}
For $t\ge t_0$, set $\tilde Q(t)=Q_\beta(t)Q_\beta(t_0)^{-1}$, whose multipliers
\[
 \tilde q_\xi(t)=\exp\Bigl(\beta\mu_\xi(W_t-W_{t_0})-\frac12\beta^2\mu_\xi^2(t-t_0)\Bigr)
\]
are continuous in $t$, uniformly in~$\xi$. In particular, $t\mapsto\tilde Q(t)^{\pm1}$ is
continuous in operator norm on $\Hs1$ and on $\cA^0$, and on every bounded interval
$\tilde Q$ belongs to some class~$\cQ_L$. Fix $\omega\in\Gamma$ outside a null set.
As in \autoref{lem:transformation}, $\tilde y=\tilde Q^{-1}u\in \rC([t_0,\tau_*);\Hs1)$
solves the time-shifted version of~\eqref{eq:mild} started at $t_0$ with initial value~$u(t_0)$. By the semigroup property, it then also solves, for every
$s\in[t_0,\tau_*)$,
\begin{equation}\label{eq:mild-restart}
 \tilde y(t)=e^{\nu(t-s)\Delta}\tilde y(s)
 -\int_s^te^{\nu(t-r)\Delta}N_{\tilde Q}(r,\tilde y(r))\rd r,\qquad s\le t<\tau_* .
\end{equation}

\emph{Restarting.} Let $s\in[t_0,\tau_*)$ with $\tilde y(s)\in\cA^0$, let $L\ge1$ bound
$\tilde Q^{\pm1}$ on $[s,s+1]$, and put $M_1(s)=\max\{1,\norm{\nabla\tilde y(s)}_2\}$
and $M_0(s)=\max\{M_1(s),\norm{\tilde y(s)}_{\cA^0}\}$. After shifting time by $s$,
\autoref{lem:mild}(a) and~(d) provide a solution $z$ of \eqref{eq:mild-restart} on
$[s,s+t_L'(s)]$ which is continuous with values in $\Hs1\cap\cA^0$, where
\[
 t_L'(s)=\min\bigl\{1,\ c_*\nu^3L^{-12}M_1(s)^{-4},\ c_*'\nu L^{-6}M_0(s)^{-2}\bigr\}.
\]
Since $z$ and $\tilde y$ both solve \eqref{eq:mild-restart} in $\rC(\Hs1)$,
\autoref{lem:mild}(b) gives $z=\tilde y$ on $[s,(s+t_L'(s))\wedge\tau_*)$. Hence,
$\tilde y$ is continuous with values in $\cA^0$ on this interval. Since
$\tilde Q^{\pm1}$ is continuous in operator norm on $\cA^0$, the same holds for
$u=\tilde Q\tilde y$, and $\norm{\tilde y}_{\cA^0}\le L\norm u_{\cA^0}\le L\,X_0(u)$.

\emph{Continuation.} Let
\[
 \tau_A=\sup\bigl\{T\in(t_0,\tau_*]:\ u\in \rC([t_0,T);\cA^0)\bigr\}.
\]
Since $\norm{u(t_0)}_{\cA^0}\le X_0(u(t_0))<\infty$ and $\tilde y(t_0)=u(t_0)$, the
restarting step with $s=t_0$ shows~$\tau_A>t_0$. Suppose that $\tau_A<\tau_*$ and
$\Lambda:=\sup_{t_0\le t<\tau_A}X_0(u(t))<\infty$. Since $[t_0,\tau_A]$ is a compact
subset of~$[0,\tau_*)$ and $u\in \rC([0,\tau_*);\Hs1)$, the quantity
$\norm{\nabla u}_2$ is bounded on $[t_0,\tau_A]$. Moreover,~$\tilde Q^{\pm1}$ is
bounded by some $L$ on $[t_0,\tau_A+1]$, and $\norm{\tilde y(s)}_{\cA^0}\le L\Lambda$
for $s\in[t_0,\tau_A)$. Hence, $M_1(s)$ and~$M_0(s)$ are bounded on $[t_0,\tau_A)$, and
there is $\delta>0$ with $t_L'(s)\ge\delta$ for all $s\in[t_0,\tau_A)$. Choosing
$s\in(\tau_A-\delta,\tau_A)$ and applying the restarting step at $s$, we find that $u$
is continuous with values in~$\cA^0$ on $[t_0,\min\{s+\delta,\tau_*\})$, which
contradicts the definition of $\tau_A$ because $s+\delta>\tau_A$. Consequently, if
$\tau_A<\tau_*$, then $\sup_{t_0\le t<\tau_A}X_0(u(t))=\infty$, and since $X_0(u)$ is
continuous on $[t_0,\tau_A)$, this means $\limsup_{t\uparrow\tau_A}X_0(u(t))=\infty$.
\end{proof}

\section{Linear one-component noise}\label{sec:linear}
Throughout this section, $(u,\tau_*)=(u^{(\beta)},\tau_*^{(\beta)})$ is the maximal solution of \eqref{eq:lin}, and for $\xi\ne0$ and~$t<\tau_*$, we write
\[
 v_\xi(t)=R_\xi\wh u(\xi,t),\qquad w_\xi(t)=S_\xi\wh u(\xi,t),\qquad\mathbf m_r(t)=\sum_{\xi\ne0}|\xi|^r\mu_\xi|v_\xi(t)|\le X_r(u(t)).
\]
Let us first explain the strategy. On each Fourier mode, the noise acts radially on $v_\xi$ and not at all on~$w_\xi$. Hence, the moduli $|v_\xi|$ and $|w_\xi|$ satisfy inequalities without It\^o correction (\autoref{lem:modes}), and the polarized norm $X_{-1}(u)$, which is just a weighted sum of these moduli, behaves like a supermartingale as long as it stays below the threshold $\rho$, with bounds which do not depend on~$\beta$ (\autoref{prop:small}). It therefore suffices to bring the solution into the region~$\{X_{-1}<\rho\}$. This is done by the noise itself on the initial layer $[0,\Theta/\beta^2]$: There, the solution is close to the linear noise flow $Q_\beta(t)u_0$ (\autoref{cor:layer}), and this flow contracts the accessible datum with high probability (\autoref{lem:contract}).

\subsection{Mode-wise moduli}\label{ss:modes}
\begin{lem}\label{lem:modes}
A.s., for all $\xi\ne0$ and all $0\le s\le t<\tau_*$, we have
\begin{align}
 |v_\xi(t)|+\nu|\xi|^2\int_s^t|v_\xi|\rd r&\le|v_\xi(s)|+\int_s^t|R_\xi\wh{\cB(u)}(\xi)|\rd r+\beta\int_s^t\mu_\xi|v_\xi|\rd W_r,\label{eq:mode-v}\\
 |w_\xi(t)|+\nu|\xi|^2\int_s^t|w_\xi|\rd r&\le|w_\xi(s)|+\int_s^t|S_\xi\wh{\cB(u)}(\xi)|\rd r.\label{eq:mode-w}
\end{align}
\end{lem}
\begin{proof}
As in the proof of \autoref{lem:transformation}, $\rd v_\xi=f_\xi\rd t+\beta\mu_\xi v_\xi\rd W$ with $f_\xi=-\nu|\xi|^2v_\xi-R_\xi\wh{\cB(u)}(\xi)$, and $z_\xi=q_\xi^{-1}v_\xi$ is absolutely continuous with $\dot z_\xi=q_\xi^{-1}f_\xi$. We regard $\mathbb C^3$ as~$\R^6$. Then $|z_\xi|$ is absolutely continuous, and for a.e.\ $t$, it holds that
\[
 \frac{\rd}{\rd t}|z_\xi|\le-\nu|\xi|^2|z_\xi|+q_\xi^{-1}|R_\xi\wh{\cB(u)}(\xi)|.
\]
For $z_\xi\ne0$, this follows from $\frac{\rd}{\rd t}|z_\xi|=\mathrm{Re}(\overline{z_\xi}\cdot\dot z_\xi)/|z_\xi|$, and for a.e.\ $t$ in $\{z_\xi=0\}$ the derivative of~$|z_\xi|$ vanishes while the right-hand side is nonnegative. Since $|z_\xi|$ has finite variation, It\^o's product rule gives $\rd|v_\xi|=\rd(q_\xi|z_\xi|)=q_\xi\rd|z_\xi|+|z_\xi|\rd q_\xi$, and $\rd q_\xi=\beta\mu_\xi q_\xi\rd W$; integrating yields~\eqref{eq:mode-v}. The proof of \eqref{eq:mode-w} is the same with $q_\xi$ replaced by $1$, since $S_\xi K_\xi=0$.
\end{proof}

\begin{rem}\label{rem:radial}
No It\^o correction appears in \eqref{eq:mode-v}, because the mode noise $\beta\mu_\xi v_\xi\rd W$ is radial for~$v_\xi$ and the Hessian of the Euclidean norm vanishes in the radial direction. For the unpolarized modulus~$|\wh u(\xi)|$, the noise $\beta\mu_\xi R_\xi\wh u(\xi)\rd W$ is not radial unless $S_\xi\wh u(\xi)=0$, and the positive correction
\[
 \frac{\beta^2\mu_\xi^2}{2|\wh u(\xi)|}\biggl(|v_\xi|^2-\frac{|v_\xi|^4}{|\wh u(\xi)|^2}\biggr)
\]
would appear. This is the mode-wise version of \eqref{eq:intro-geometry}, and the reason for working with $X_r$ instead of~$\norm\cdot_{\cA^r}$.
\end{rem}

\subsection{Summation over the modes}\label{ss:summation}
\begin{lem}\label{lem:Xminus}
Let $t_0\ge0$ be deterministic. A.s.\ on $\{t_0<\tau_*\}$, for all $t\in[t_0,\tau_*)$, we have
\begin{equation}\label{eq:Xminus}
 X_{-1}(u(t))+\nu\int_{t_0}^tX_1(u)\rd s\le X_{-1}(u(t_0))+\sqrt2\int_{t_0}^tX_0(u)^2\rd s+\beta\int_{t_0}^t\mathbf m_{-1}\rd W_s.
\end{equation}
In particular, $\int_{t_0}^tX_1(u)\rd s<\infty$ for all $t<\tau_*$.
\end{lem}
\begin{proof}
Let $\tau_n=\inf\{t\ge0:\norm{\nabla u(t)}_2\ge n\text{ or }\int_0^t\norm{\Delta u}_2^2\rd s\ge n\}\wedge n$. Then $\tau_n\uparrow\tau_*$, and on $[0,\tau_n]$, we have $X_{-1}(u)\le C_En$ and $\int_0^{\tau_n}X_0(u)^2\rd s\le Cn$ by~\eqref{eq:H1-embeddings}. Let $F_j$ be finite sets increasing to $\Z^3\setminus\{0\}$. We add \eqref{eq:mode-v} and \eqref{eq:mode-w} with weights $|\xi|^{-1}$ over $\xi\in F_j$, with $s=t_0$ and $t\le\tau_n$ (on $\{t_0<\tau_n\}$), and bound the drift by the full sum $X_{-1}(\cB(u))\le\sqrt2X_0(u)^2$ from~\eqref{eq:wiener-products}. This gives \eqref{eq:Xminus} with $X_{\pm1}$ and~$\mathbf m_{-1}$ replaced by the corresponding sums over~$F_j$. As $j\to\infty$, the left-hand side converges pointwise, the viscous term monotonically. 
The martingale integrands $\mathbf m_{-1}^{F_j}=\sum_{\xi\in F_j}|\xi|^{-1}\mu_\xi|v_\xi|$ increase to $\mathbf m_{-1}\le X_{-1}(u)\le C_En$ on $[0,\tau_n]$, and $\tau_n\le n$. By dominated convergence, the It\^o isometry and Doob's inequality, the stochastic integrals, stopped at $\tau_n$, converge in $\rL^2(\Omega;\rC([t_0,n]))$, hence a.s.\ uniformly along a subsequence. This proves \eqref{eq:Xminus} on $[t_0,\tau_n]$, and $n\to\infty$ completes the proof.
\end{proof}

\subsection{Small-state continuation}\label{ss:small}
\begin{prop}[Small-state continuation, uniformly in $\beta$]\label{prop:small}
Let $\rho=\nu/2\sqrt2$ be as in \eqref{eq:rho}, let~$t_0\ge0$ be deterministic,
$x\in(0,\rho)$, and $\Gamma\in\cF_{t_0}$ with $\Gamma\subset\{t_0<\tau_*,\ X_{-1}(u(t_0))\le x\}$. Let
\[
 \sigma=\inf\{t\in[t_0,\tau_*):X_{-1}(u(t))\ge\rho\},
\]
where $\inf\emptyset=\infty$.
Then,
\begin{equation}\label{eq:small-prob}
    \begin{aligned}
        \bP(\Gamma\cap\{\sigma<\infty\})
        &\le\frac1\rho\,\E\bigl[\mathbf 1_\Gamma X_{-1}(u(t_0))\bigr]\le\frac x\rho\,\bP(\Gamma),\\
        \E\Bigl[\mathbf 1_\Gamma\int_{t_0}^{\sigma\wedge\tau_*}X_1(u)\rd t\Bigr]
        &\le\frac2\nu\,\E\bigl[\mathbf 1_\Gamma X_{-1}(u(t_0))\bigr],
    \end{aligned}
\end{equation}
and $\tau_*=\infty$ a.s.\ on $\Gamma\cap\{\sigma=\infty\}$. These bounds do not depend on~$\beta$.
\end{prop}
\begin{proof}
For $t\in[t_0,\sigma\wedge\tau_*)$ we have $X_{-1}(u)<\rho$, and hence $\sqrt2X_0(u)^2\le\sqrt2X_{-1}(u)X_1(u)\le\frac\nu2X_1(u)$ by \eqref{eq:wiener-products} and the choice of $\rho$ in~\eqref{eq:rho}. Thus,
\eqref{eq:Xminus} implies
\begin{equation}\label{eq:small-energy}
 X_{-1}(u(t))+\frac\nu2\int_{t_0}^tX_1(u)\rd s\le X_{-1}(u(t_0))+\beta\int_{t_0}^t
 \mathbf m_{-1}\rd W_s,\qquad t_0\le t\le\sigma,\ t<\tau_* ,
\end{equation}
where the case $t=\sigma$ follows by continuity. Recall that $X_{-1}(u(\cdot))$ is
continuous on $[t_0,\tau_*)$, since $X_{-1}\le C_E\norm{\nabla\cdot}_2$. Hence, on
$\Gamma\cap\{\sigma<\infty\}$, the infimum defining $\sigma$ is attained, $\sigma<\tau_*$,
and, as $X_{-1}(u(t_0))\le x<\rho$, we have $X_{-1}(u(\sigma))=\rho$.

Let $\tau_n$ be as in the proof of \autoref{lem:Xminus}, $\Gamma_n=\Gamma\cap\{t_0<\tau_n\}
\in\cF_{t_0}$ and $\theta_n=\sigma\wedge\tau_n$. On~$[t_0,\theta_n]$, we have
$0\le\mathbf m_{-1}\le X_{-1}(u)\le\rho$. Since $\mathbf 1_{\Gamma_n}$ is
$\cF_{t_0}$-measurable, it can be moved inside the stochastic integral.
Then $\int_{t_0}^{t\wedge\theta_n}\mathbf 1_{\Gamma_n}\mathbf m_{-1}\rd W_s$ is a square-integrable martingale due to the boundedness of the integrand and~$\theta_n\le n$, so its expectation
vanishes. Multiplying \eqref{eq:small-energy} at the time $t\wedge\theta_n$ by
$\mathbf 1_{\Gamma_n}$ and taking expectations, we obtain, with
$x_\Gamma=\E[\mathbf 1_\Gamma X_{-1}(u(t_0))]$,
\begin{equation}\label{eq:small-expectation}
 \E\bigl[\mathbf 1_{\Gamma_n}X_{-1}(u(t\wedge\theta_n))\bigr]+\frac\nu2\,
 \E\biggl[\mathbf 1_{\Gamma_n}\int_{t_0}^{t\wedge\theta_n}X_1(u)\rd s\biggr]\le x_\Gamma,
 \qquad t\ge t_0 .
\end{equation}

\emph{The integral bound.} Since $\tau_n\uparrow\tau_*$, we have $\Gamma_n\uparrow\Gamma$ and
$\theta_n\uparrow\sigma\wedge\tau_*$. Monotone convergence in
\eqref{eq:small-expectation}, first as $n\to\infty$ and then as $t\to\infty$, gives the
second bound in~\eqref{eq:small-prob}.

\emph{Global existence on $\Gamma\cap\{\sigma=\infty\}$.} By the integral bound,
$\int_{t_0}^{\tau_*}X_1(u)\rd t<\infty$ a.s.\ on $\Gamma\cap\{\sigma=\infty\}\cap
\{\tau_*<\infty\}$, and $\int_0^{t_0}X_1(u)\rd t<\infty$ there by \autoref{lem:Xminus}
with $t_0=0$. By \autoref{lem:continuation}, this event is a null set, so $\tau_*=\infty$
a.s.\ on $\Gamma\cap\{\sigma=\infty\}$.

\emph{The exit probability.} Fix $t\ge t_0$ and put
$f_n=\mathbf 1_{\Gamma_n}X_{-1}(u(t\wedge\theta_n))\ge0$. On $\Gamma\cap\{\sigma\le t\}$,
we have $\sigma<\tau_*$, so $\tau_n>\sigma\ge t_0$ for all $n\ge N(\omega)$; for such $n$,
$\omega\in\Gamma_n$ and $t\wedge\theta_n=\sigma$, whence $f_n=X_{-1}(u(\sigma))=\rho$.
Thus, $\liminf_nf_n\ge\rho\,\mathbf 1_{\Gamma\cap\{\sigma\le t\}}$, and Fatou's lemma
together with \eqref{eq:small-expectation} yields
\[
 \rho\,\bP\bigl(\Gamma\cap\{\sigma\le t\}\bigr)\le\E\bigl[\liminf_nf_n\bigr]
 \le\liminf_n\E[f_n]\le x_\Gamma .
\]
Letting $t\to\infty$ and using $x_\Gamma\le x\,\bP(\Gamma)$ gives the first bound in~\eqref{eq:small-prob}.
\end{proof}

The bounds are uniform in $\beta$ although $\beta$ appears in the martingale term: Before the exit, the martingale has a bounded integrand, and its expectation vanishes regardless of the size of $\beta$. In this sense the small region $\{X_{-1}<\rho\}$ is stable for every intensity.

\subsection{Contraction by the noise flow}\label{ss:contraction}
In this subsection, $B$ denotes a standard real Brownian motion. For $\tau\ge0$, let
\begin{equation}\label{eq:noise-flow}
 \cQ_\tau=\exp\bigl(KB_\tau-\tfrac12\tau K^2\bigr),\qquad(\cQ_\tau)_\xi=S_\xi+e^{\mu_\xi B_\tau-\mu_\xi^2\tau/2}R_\xi.
\end{equation}
Thus, $\cQ_\tau$ is the operator $Q_\beta$ from \eqref{eq:Qbeta} with $\beta=1$.
In fact, by Brownian
scaling, $Q_\beta(t)=\cQ_{\beta^2t}$ if $\cQ$ is computed with the Brownian motion
$B_\tau=\beta W_{\tau/\beta^2}$, see also \eqref{eq:rescaled} below.
Since $\mu_\xi\in[0,1]$, we have $|\mu_\xi B_\tau-\frac12\mu_\xi^2\tau|\le|B_\tau|+\frac\tau2$, and therefore, for $\Theta<\infty$,
\begin{equation}\label{eq:path-bound}
 \sup_{\tau\le\Theta}\sup_{\xi\ne0}\max\bigl\{e^{\mu_\xi B_\tau-\mu_\xi^2\tau/2},e^{-\mu_\xi B_\tau+\mu_\xi^2\tau/2}\bigr\}\le\exp\Bigl(\sup_{\tau\le\Theta}|B_\tau|+\frac\Theta2\Bigr).
\end{equation}

\begin{lem}[Contraction of a fixed accessible datum]\label{lem:contract}
Let $a\in\Ran\cR\cap\Hs1$, $d>0$, and $\gamma\in(0,1)$. There exist deterministic $\Theta,L<\infty$ and an event $\Omega_*\in\sigma(B_\tau:\tau\le\Theta)$ with $\bP(\Omega_*)\ge1-\gamma$ such that, on $\Omega_*$,
\begin{equation}\label{eq:contract}
 \norm{\nabla\cQ_\Theta a}_2\le d,\qquad\sup_{\tau\le\Theta}\sup_{\xi\ne0}\max\bigl\{e^{\mu_\xi B_\tau-\mu_\xi^2\tau/2},e^{-\mu_\xi B_\tau+\mu_\xi^2\tau/2}\bigr\}\le L.
\end{equation}
If moreover $a\in\cA^0$, then $\Theta$, $L$ and $\Omega_*$ can be chosen such that, in addition, $\norm{\cQ_\Theta a}_{\cA^0}\le d$ on $\Omega_*$.
\end{lem}
\begin{proof}
Choose $z>0$ with $2e^{-z^2/2}\le\gamma/2$. For $\kappa>0$, denote by $R_{<\kappa}$ the Fourier projection onto the modes with $0<\mu_\xi<\kappa$.
Since $a\in\Ran\cR$ and $R_\xi=0$ for $\xi_h=0$, the Fourier coefficients $\wh a(\xi)$
vanish on the modes with $\mu_\xi=0$. Hence, the indicators of $\{0<\mu_\xi<\kappa\}$
tend to zero as $\kappa\downarrow0$ on the support of $\wh a$, and dominated convergence
gives $\norm{\nabla R_{<\kappa}a}_2\to0$ as $\kappa\downarrow0$.
We fix $\kappa\in(0,1]$ with $e^{z^2/2}\norm{\nabla R_{<\kappa}a}_2\le d/2$. Next, we choose $\Theta$ so large that
\[
 \kappa\sqrt\Theta\ge4z\qquad\text{and}\qquad e^{-\kappa^2\Theta/4}\norm{\nabla a}_2\le d/2 .
\]
On $\Omega_1=\{|B_\Theta|\le z\sqrt\Theta\}$, completing the square gives, for every $\mu\ge0$,
\[
 \mu B_\Theta-\tfrac12\mu^2\Theta\le\mu z\sqrt\Theta-\tfrac12\mu^2\Theta\le\tfrac12z^2,
\]
and for $\mu\ge\kappa$, since then $z\le\frac14\mu\sqrt\Theta$, the same expression is at most
\[
 \mu\sqrt\Theta\Bigl(\frac14\mu\sqrt\Theta-\frac12\mu\sqrt\Theta\Bigr)=-\frac14\mu^2\Theta\le-\frac14\kappa^2\Theta .
\] Splitting $a=R_{<\kappa}a+(a-R_{<\kappa}a)$ into parts with disjoint Fourier supports, on $\Omega_1$, we obtain 
\[
 \norm{\nabla\cQ_\Theta a}_2\le e^{z^2/2}\norm{\nabla R_{<\kappa}a}_2+e^{-\kappa^2\Theta/4}\norm{\nabla a}_2\le d.
\]
Finally, choose $b>0$ with $4e^{-b^2/(2\Theta)}\le\gamma/2$, let $\Omega_2=\{\sup_{\tau\le\Theta}|B_\tau|\le b\}$ and $L=e^{b+\Theta/2}$. 
By the reflection principle, applied to $B$ and $-B$, and the Gaussian tail bound
$\bP(N\ge y)\le e^{-y^2/2}$ for a standard normal variable $N$ and~$y\ge0$,
we get
\[
 \bP(\Omega_2^c)\le4\,\bP(B_\Theta\ge b)\le4e^{-b^2/(2\Theta)}\le\frac\gamma2,
 \qquad
 \bP(\Omega_1^c)=\bP(|N|>z)\le2e^{-z^2/2}\le\frac\gamma2 .
\]
Thus, $\Omega_*=\Omega_1\cap\Omega_2$ has the required properties, by~\eqref{eq:path-bound}. If $a\in\cA^0$, we choose $\kappa$ and $\Theta$ such that the same two smallness conditions also hold for $\norm{\cdot}_{\cA^0}$ in place of $\norm{\nabla\cdot}_2$; the argument is unchanged.
\end{proof}

\begin{rem}\label{rem:contraction}
(1) Each multiplier $e^{\mu_\xi B_\tau-\mu_\xi^2\tau/2}$ has expectation one. The contraction holds on an event of high probability and is compensated by rare large realizations. It is not an estimate in expectation, and it relies on the It\^o correction $-\frac12\mu_\xi^2\tau$, cf.~\autoref{rem:ito}.

(2) Since $\mu_\xi$ has no positive lower bound on accessible modes, the contraction time $\Theta$ depends on how much of $a$ lives on modes with small~$\mu_\xi$. This is why the parameters in \autoref{thm:linear} depend on the Fourier tail of $u_0$ and are not uniform on bounded sets. No spectral gap is used.
\end{rem}

\subsection{Proof of parts (a) and (b)}\label{ss:proof-ab}
\begin{proof}[Proof of \autoref{thm:linear}(a)]
Fix $\eps\in(0,1)$ and choose $x\in(0,\rho)$ with $x/\rho\le\eps/2$. Since $\cS u_0=0$, we can apply \autoref{lem:contract} to $a=u_0$, with $d=x/(2C_E)$ and $\gamma=\eps/2$.
We keep the resulting $\Theta$ and $L$ fixed. For $\beta>0$ and $\tau\ge0$, put
\begin{equation}\label{eq:rescaled}
 t_\beta=\Theta/\beta^2,\qquad B^{(\beta)}_\tau=\beta W_{\tau/\beta^2}.
\end{equation}
By Brownian scaling, $B^{(\beta)}$ is a standard Brownian motion, and $Q_\beta(t)=\cQ_{\beta^2t}$ when $\cQ$ is computed with $B^{(\beta)}$. In particular, $Q_\beta(t_\beta)=\cQ_\Theta$, and the multipliers $q_\xi(t)$, $t\le t_\beta$, are those in \eqref{eq:contract} with $\tau=\beta^2t\le\Theta$. Let $\Omega_\beta$ be the event of \autoref{lem:contract} for~$B^{(\beta)}$. 
Then $\Omega_\beta\in\cF_{t_\beta}$, $\bP(\Omega_\beta)\ge1-\eps/2$ for every $\beta$, and
$\Omega_\beta\subset\Omega_{L,t_\beta}$, where $\Omega_{L,T}$ is the event from
\autoref{cor:layer}. Indeed, the multiplier $q_\xi(t)$ coincides with
$e^{\mu_\xi B^{(\beta)}_\tau-\mu_\xi^2\tau/2}$ at $\tau=\beta^2t$, and $\tau\le\Theta$
if and only if $t\le t_\beta$.

Let $M_1=\max\{1,\norm{\nabla u_0}_2\}$, and let $t_L$ be as in \eqref{eq:tL}, i.e.,
$t_L=\min\{1,c_*\nu^3L^{-12}M_1^{-4}\}$. Since $L$ and $M_1$ are fixed, $t_L$ does not
depend on~$\beta$. Choose $\beta_0$ such that for all $\beta\ge\beta_0$, it holds that
\begin{equation}\label{eq:beta-choice}
 t_\beta\le t_L,\qquad L\,\omega_{u_0}(t_\beta)+C_1L^4\nu^{-3/4}t_\beta^{1/4}M_1^2\le\frac x{2C_E}.
\end{equation}
This is possible since $\omega_{u_0}(t)\to0$ as $t\downarrow0$ by strong continuity of the heat semigroup. By \autoref{cor:layer}, a.s.\ on $\Omega_\beta$, we have $\tau_*>t_\beta$ and
\[
 \norm{\nabla u(t_\beta)}_2\le\norm{\nabla\cQ_\Theta u_0}_2+\frac x{2C_E}\le\frac x{C_E},\qquad\text{hence}\qquad X_{-1}(u(t_\beta))\le x
\]
by~\eqref{eq:H1-embeddings}.
We now apply \autoref{prop:small} with the deterministic time $t_0=t_\beta$ and
$\Gamma=\Omega_\beta\cap\{t_\beta<\tau_*\}$. Here $\Gamma\in\cF_{t_\beta}$, since
$\Omega_\beta\in\cF_{t_\beta}$ and $\tau_*$ is a stopping time, and
$X_{-1}(u(t_\beta))\le x$ on $\Gamma$ after modification on a null set, which
changes neither the measurability (the filtration satisfies the usual conditions)
nor~$\bP(\Gamma)$. Then
the exit time $\sigma$ satisfies $\bP(\Gamma\cap\{\sigma<\infty\})\le x/\rho\le\eps/2$, and $\tau_*=\infty$ a.s.\ on $\Gamma\cap\{\sigma=\infty\}$. Therefore,
\[
 \bP(\tau_*<\infty)\le\bP(\Omega_\beta^c)+\bP(\Gamma\cap\{\sigma<\infty\})\le\eps\qquad\text{for all }\beta\ge\beta_0.\qedhere
\]
\end{proof}

\begin{rem}\label{rem:layer}
The noise coefficient is $\beta K$ at all times. The initial layer $[0,t_\beta]$ is only a device in the proof. The full nonlinear term, including the pressure and the generation of an unforced component, is contained in the comparison error, which by \eqref{eq:mild} equals
\[
 u(t)-Q_\beta(t)u_0=Q_\beta(t)\Bigl[(e^{\nu t\Delta}-\Id)u_0-\int_0^te^{\nu(t-s)\Delta}Q_\beta(s)^{-1}\cB(u(s))\rd s\Bigr].
\]
It is small because the layer is short in physical time, while the noise acts on the time scale $\beta^{-2}$ and completes its contraction within the layer. No deterministic term is dropped or tamed.
\end{rem}

\begin{proof}[Proof of \autoref{thm:linear}(b)]
Let $s_0=\norm{\cS u_0}_{\cA^{-1}}$. If $s_0\ge\rho$ there is nothing to prove, so let $s_0<\rho$, fix~$\delta>0$ with $s_0+\delta<\rho$, and $\gamma\in(0,1)$. 
We apply \autoref{lem:contract} to $a=\cR u_0\in\Ran\cR\cap\Hs1$, with $d=\delta/(2C_E)$
and failure tolerance $\gamma$, and keep the resulting $\Theta$ and $L$ fixed. With
$t_\beta$ and $B^{(\beta)}$ as in~\eqref{eq:rescaled}, let $\Omega_\beta$ be the event of
\autoref{lem:contract} for $B^{(\beta)}$. As in the proof of (a), $\Omega_\beta\in
\cF_{t_\beta}$, $\bP(\Omega_\beta)\ge1-\gamma$, and $\Omega_\beta\subset
\Omega_{L,t_\beta}$. We choose $\beta_0=\beta_0(\gamma,\delta)$ such that
\eqref{eq:beta-choice} holds with $\delta$ in place of $x$ for all $\beta\ge\beta_0$.
Since $Q_\beta(t)$ acts as the identity on $\Ran\cS$ and by multiplication with $q_\xi(t)$
on each accessible mode, we have $Q_\beta(t_\beta)u_0=\cQ_\Theta\cR u_0+\cS u_0$. Moreover,
$X_{-1}(f+g)=X_{-1}(f)+X_{-1}(g)$ for $f\in\Ran\cR$ and $g\in\Ran\cS$, so that
\[
 X_{-1}(Q_\beta(t_\beta)u_0)=X_{-1}(\cQ_\Theta\cR u_0)+\norm{\cS u_0}_{\cA^{-1}}\le C_E\norm{\nabla\cQ_\Theta\cR u_0}_2+s_0 .
\]
By \autoref{cor:layer}, a.s.\ on $\Omega_\beta$, we have $\tau_*>t_\beta$, so $u(t_\beta)$
is defined, and we split
\[
 u(t_\beta)=Q_\beta(t_\beta)u_0+\bigl(u(t_\beta)-Q_\beta(t_\beta)u_0\bigr)
\]
into the linear noise flow applied to the datum and the comparison error. On
$\Omega_\beta$, the choice $d=\delta/(2C_E)$ gives
$C_E\norm{\nabla\cQ_\Theta\cR u_0}_2\le\delta/2$, and \eqref{eq:comparison} together with
\eqref{eq:beta-choice} (with $\delta$ in place of $x$) gives
$C_E\norm{\nabla(u(t_\beta)-Q_\beta(t_\beta)u_0)}_2\le\delta/2$. Using the subadditivity
of $X_{-1}$ and \eqref{eq:H1-embeddings}, we obtain a.s.\ on $\Omega_\beta$ the estimate
\[
 X_{-1}(u(t_\beta))\le X_{-1}(Q_\beta(t_\beta)u_0)
 +C_E\norm{\nabla(u(t_\beta)-Q_\beta(t_\beta)u_0)}_2
 \le\Bigl(s_0+\frac\delta2\Bigr)+\frac\delta2=s_0+\delta .
\]
\autoref{prop:small}, applied with $t_0=t_\beta$, $x=s_0+\delta$ and
$\Gamma=\Omega_\beta\cap\{t_\beta<\tau_*\}$ as in the proof of (a), gives
$\bP(\Gamma\cap\{\sigma<\infty\})\le(s_0+\delta)/\rho$ and $\tau_*=\infty$ a.s.\ on
$\Gamma\cap\{\sigma=\infty\}$. Since $\Gamma$ and $\Omega_\beta$ differ only by a null set,
\[
 \bP(\tau_*<\infty)\le\bP(\Omega_\beta^c)+\bP(\Gamma\cap\{\sigma<\infty\})
 \le\gamma+\frac{s_0+\delta}\rho\qquad\text{for all }\beta\ge\beta_0(\gamma,\delta).
\]
Hence, $\liminf_{\beta\to\infty}\bP(\tau_*^{(\beta)}=\infty)\ge1-\gamma-(s_0+\delta)/\rho$.
Since $\Theta$, $L$, and $\beta_0$ depend on $\gamma$ and $\delta$, but the left-hand side
does not, we may now let $\gamma,\delta\to0$, which proves (b).
\end{proof}

\subsection{Space-time integrals of the vertical component}\label{ss:integrals}
For the integrals of $u_3$, we also need the norm $X_0$, which is not controlled by the $\Hs1$-norm. Its evolution is described by the analogue of~\autoref{lem:Xminus}.

\begin{lem}\label{lem:Xzero}
Let $t_0$, $\Gamma$, and $\tau_A$ be as in~\autoref{lem:persistence}. Then a.s.\ on $\Gamma$, for all $t\in[t_0,\tau_A)$,
\begin{equation}\label{eq:Xzero}
 X_0(u(t))+\int_{t_0}^t\bigl(\nu-\sqrt2X_{-1}(u)\bigr)X_2(u)\rd s\le X_0(u(t_0))+\beta\int_{t_0}^t\mathbf m_0\rd W_s .
\end{equation}
\end{lem}
\begin{proof}
We argue as in the proof of \autoref{lem:Xminus}, now with weights $|\xi|^0$, and localize with $\theta_\Lambda=\inf\{t\ge t_0:X_0(u(t))>\Lambda\}$ in addition to $\tau_n$.
Then $\theta_\Lambda$ is a stopping time by lower semicontinuity, and $\theta_\Lambda<\tau_A$ on $\{\tau_A<\tau_*\}$ by~\autoref{lem:persistence}. Before $\theta_\Lambda$, by \eqref{eq:wiener-products}, the drift satisfies $X_0(\cB(u))\le\sqrt2X_0(u)X_1(u)\le\sqrt2\Lambda X_1(u)$, which is integrable by \autoref{lem:Xminus}, and the martingale integrand satisfies $\mathbf m_0\le X_0(u)\le\Lambda$. Passing to the limit along finite sets of modes as before, we obtain \eqref{eq:Xzero} with the drift $\sqrt2\int X_0X_1$ on the right-hand side and a finite viscous term $\nu\int X_2$ on the left-hand side. Since $X_0X_1\le X_{-1}X_2$ by \eqref{eq:wiener-products}, \eqref{eq:Xzero} follows on $[t_0,\theta_\Lambda\wedge\tau_n]$, and we let $\Lambda,n\to\infty$.
\end{proof}

Before the exit time of \autoref{prop:small}, we have $\sqrt2X_{-1}(u)<\sqrt2\rho=\nu/2$, so the coefficient of $X_2$ in \eqref{eq:Xzero} is at least $\nu/2$, and $X_0(u)$ is a nonnegative local supermartingale there.

\begin{proof}[Proof of \autoref{thm:linear}(c)]
\emph{Strategy.} On a good event $G_\beta$ constructed as in the proof of (a), the
solution is global and $\int_{t_\beta}^\infty X_1(u)\rd t<\infty$. Since
$\norm{u_3}_\infty\le X_0(u)$ and $X_0(u)^2\le X_{-1}(u)X_1(u)$, the sixth power of $u_3$
is controlled by $X_0(u)^4X_1(u)$, so it suffices to show that $X_0(u)$ stays bounded on~$G_\beta$. As $X_0(u)$ need not be finite at time $t_\beta$, we restart at a time $r$ with
$u(r)\in\Hs2$, propagate~$\cA^0$-regularity by \autoref{lem:persistence}, and use the
supermartingale property of $X_0(u)$ from \autoref{lem:Xzero} to rule out that $X_0(u)$
becomes unbounded.

\emph{Step 1: The good event.} Fix $\eps\in(0,1)$, and let $x$, $\Theta$, $L$,
$\beta_0$, $t_\beta$, and $\Omega_\beta$ be as in the proof of (a). For $\beta\ge\beta_0$
let $\Gamma=\Omega_\beta\cap\{t_\beta<\tau_*\}$, let $\sigma$ be the exit time from
\autoref{prop:small} with $t_0=t_\beta$, and set $G_\beta=\Gamma\cap\{\sigma=\infty\}$.
By the proof of (a), $\bP(G_\beta)\ge1-\eps$. Moreover, a.s.\ on $G_\beta$, we have
$\tau_*=\infty$, $X_{-1}(u)<\rho$ on $[t_\beta,\infty)$, and, by \eqref{eq:small-prob},
$\int_{t_\beta}^\infty X_1(u)\rd t<\infty$.

\emph{Step 2: Restart at a time of $\Hs2$-regularity.} Since
$u\in \rL^2_{\rm loc}([0,\tau_*);\Hs2)$ a.s., Fubini's theorem provides a deterministic
$r>t_\beta$ such that $u(r)\in\Hs2$ a.s.\ on $\{r<\tau_*\}$. By \eqref{eq:H1-embeddings},
$X_0(u(r))\le C\norm{\Delta u(r)}_2<\infty$ a.s.\ on $\{r<\tau_*\}$. This finiteness is
the input required in \autoref{lem:persistence} and \autoref{lem:Xzero}, which we now
apply with $t_0=r$ and $\Gamma=\{r<\tau_*\}$ (modified on a null set).
Let $\tau_A$ be as
in~\autoref{lem:persistence}.

\emph{Step 3: $X_0(u)$ stays bounded on~$G_\beta$.} On $[r,\sigma\wedge\tau_A)$, we have
$\nu-\sqrt2X_{-1}(u)\ge\nu/2$, so \eqref{eq:Xzero} yields
\begin{equation}\label{eq:Xzero-super}
 X_0(u(t))\le X_0(u(r))+\beta\int_r^t\mathbf m_0\rd W_s,\qquad r\le t<\sigma\wedge\tau_A .
\end{equation}
For $\Lambda>0$, let $\theta_\Lambda$ be the first time after $r$ at which $X_0(u)$ exceeds
$\Lambda$, and put $\vartheta_n=r\vee(\theta_\Lambda\wedge\sigma\wedge\tau_n)$, with
$\tau_n$ as in the proof of~\autoref{lem:Xminus}. By \autoref{lem:persistence},
$\theta_\Lambda<\tau_A$ on $\{\tau_A<\tau_*\}$, and $\tau_n<\tau_*$. 
Hence~\eqref{eq:Xzero-super} applies up to time $\vartheta_n$, and $\vartheta_n\le n\vee r$.

Fix $0<K<\Lambda$, and let $\Gamma_K=\{r<\tau_*,\ X_0(u(r))\le K\}\in\cF_r$. On $\Gamma_K$ and on $[r,\vartheta_n]$, we have $0\le\mathbf m_0\le X_0(u)\le\Lambda$. Exactly as in the proof of \autoref{prop:small}, multiplying \eqref{eq:Xzero-super} at time $t\wedge\vartheta_n$ by the $\cF_r$-measurable indicator $\mathbf 1_{\Gamma_K}$ and taking expectations gives $\E[\mathbf 1_{\Gamma_K}X_0(u(t\wedge\vartheta_n))]\le K$ for $t\ge r$. On $\Gamma_K\cap\{\theta_\Lambda\le t,\ \theta_\Lambda<\sigma\wedge\tau_A\}$, we have $\theta_\Lambda<\tau_*$, hence $\vartheta_n=\theta_\Lambda$ for all large $n$, and $X_0(u(\theta_\Lambda))=\Lambda$, because $X_0(u)$ is continuous on $[r,\tau_A)$ and $X_0(u(r))\le K<\Lambda$. Fatou's lemma and the limit $t\to\infty$ give
\[
 \Lambda\,\bP\bigl(\Gamma_K\cap\{\theta_\Lambda<\sigma\wedge\tau_A\}\bigr)\le K .
\]
On $G_\beta$, we have $\tau_*=\infty>r$ and $\sigma=\infty$, and $\theta_\Lambda<\infty$
implies $\theta_\Lambda<\tau_A$: Either $\tau_A=\infty$, or $\tau_A<\infty=\tau_*$, in
which case $X_0(u)$ is unbounded before $\tau_A$ by~\autoref{lem:persistence}. Therefore,
\[
 \bP\bigl(G_\beta\cap\{\theta_\Lambda<\infty\}\bigr)
 \le\bP\bigl(r<\tau_*,\ X_0(u(r))>K\bigr)+\frac K\Lambda .
\]
Letting first $\Lambda\to\infty$ and then $K\to\infty$, the right-hand side tends to zero,
since $X_0(u(r))<\infty$ a.s.\ on $\{r<\tau_*\}$ by Step~2. As
$\{\sup_{t\ge r}X_0(u(t))=\infty\}\subset\{\theta_\Lambda<\infty\}$ for every $\Lambda$, we
conclude that, a.s.\ on $G_\beta$, $\tau_A=\infty$ and
$\Lambda_*:=\sup_{t\ge r}X_0(u(t))<\infty$.

\emph{Step 4: The space-time integral.} By Step~1, $X_0(u)^2\le X_{-1}(u)X_1(u)\le\rho X_1(u)$ on $[r,\infty)$. Together with $\norm{u_3}_\infty\le X_0(u)$, we find a.s.\ on~$G_\beta$ the estimate
\[
 \int_r^\infty\!\!\int_{\T^3}|u_3|^6\rd x\rd t\le(2\pi)^3\int_r^\infty X_0(u)^6\rd t
 \le(2\pi)^3\rho\,\Lambda_*^4\int_r^\infty X_1(u)\rd t<\infty .
\]
On $[0,r]$, the integral is finite, since $u\in \rC([0,r];\Hs1)$ and
$\norm{u_3}_6\le C\norm{\nabla u}_2$. Hence,
\[
 \bP\Bigl(\tau_*=\infty,\ \int_0^\infty\!\!\int_{\T^3}|u_3|^6\rd x\rd t<\infty\Bigr)
 \ge\bP(G_\beta)\ge1-\eps\qquad\text{for all }\beta\ge\beta_0 .\qedhere
\]
\end{proof}

\begin{rem}\label{rem:powers}
On compact subintervals of $[0,\tau_*)$ the solution belongs to $\rL^{10}_{t,x}$: The embedding $\rH^{6/5}\hookrightarrow \rL^{10}$ and interpolation give $\norm u_{10}^{10}\le C\norm{\nabla u}_2^8\norm{\Delta u}_2^2$. Combined with the tail argument above, which controls every power, part (c) holds simultaneously for the even powers $2,4,\dots,10$. For higher powers the tail is still controlled, but integrability near $t=0$ is not automatic for $\Hs1$ data.
This is the reason for the additional assumption in (d). The bound on $X_0$ obtained in the proof of~(c) is not uniform in $\beta$, which is why (c) asserts finiteness and not smallness.
\end{rem}

\begin{proof}[Proof of \autoref{thm:linear}(d)]
Fix $m\ge1$, $\eta>0$, $\eps\in(0,1)$, and a level $\Lambda>0$, and choose $x\in(0,\min\{\rho,\Lambda\})$ so small that
\begin{equation}\label{eq:tail-probability}
 \frac x\rho+\frac x\Lambda+\frac{4(2\pi)^3\rho\Lambda^{2m-2}x}{\nu\eta}\le\frac\eps2 .
\end{equation}
We apply \autoref{lem:contract} to $a=u_0\in\Ran\cR\cap\cA^0$ with $d=x/4$ and
$\gamma=\eps/2$, so that, on the resulting event, both $\norm{\nabla\cQ_\Theta u_0}_2$ and
$\norm{\cQ_\Theta u_0}_{\cA^0}$ are at most $x/4$ (only the second bound is used below), and define $t_\beta$, $\Omega_\beta$ as in the proof of (a). Let $M_0=\max\{1,\norm{\nabla u_0}_2,\norm{u_0}_{\cA^0}\}$. 
By \autoref{cor:layer}, for all sufficiently large $\beta$, we have $t_\beta\le t_L'$ and
\[
 L\sup_{s\le t_\beta}\norm{(e^{\nu s\Delta}-\Id)u_0}_{\cA^0}
 +C_1'L^4\nu^{-1/2}t_\beta^{1/2}M_0^2\le\frac x4 ,
\]
the first term tending to zero by strong continuity of the heat semigroup on~$\cA^0$. Hence,
a.s.\ on~$\Omega_\beta$,~$\tau_*>t_\beta$ and
\[
 X_0(u(t_\beta))\le\sqrt2\bigl(\norm{\cQ_\Theta u_0}_{\cA^0}+\norm{u(t_\beta)-Q_\beta(t_\beta)u_0}_{\cA^0}\bigr)\le\sqrt2\Bigl(\frac x4+\frac x4\Bigr)<x,
\]
so that also $X_{-1}(u(t_\beta))\le X_0(u(t_\beta))<x$. Moreover, $\sup_{t\le t_\beta}X_0(u(t))\le4LM_0$ by \eqref{eq:comparison-A0}, so that, using $\norm{u_3}_\infty\le X_0(u)$,
\[
 \int_0^{t_\beta}\!\!\int_{\T^3}|u_3|^{2m}\rd x\rd t\le(2\pi)^3t_\beta(4LM_0)^{2m}\le\frac\eta2
\]
for large $\beta$, since $L$ and $M_0$ are fixed.

After $t_\beta$, we apply \autoref{lem:persistence} and \autoref{lem:Xzero} with
$t_0=t_\beta$ and $\Gamma=\Omega_\beta\cap\{t_\beta<\tau_*\}\in\cF_{t_\beta}$ (modified on a
null set), on which $X_{-1}(u(t_\beta))\le X_0(u(t_\beta))<x$. Let $\sigma$ be the exit
time from \autoref{prop:small} with $t_0=t_\beta$, let $\theta_\Lambda$ be the first time
after $t_\beta$ at which $X_0(u)$ exceeds $\Lambda$, and put
$\vartheta=\sigma\wedge\theta_\Lambda$. By \autoref{prop:small},
$\bP(\Gamma\cap\{\sigma<\infty\})\le x/\rho$, and $\tau_*=\infty$ a.s.\ on
$\Gamma\cap\{\sigma=\infty\}$. Moreover, the argument in the proof of (c), with $r$
replaced by $t_\beta$, $\Gamma_K$ replaced by $\Gamma$, and $K$ replaced by $x$ (recall
$X_0(u(t_\beta))<x<\Lambda$ on $\Gamma$), gives
$\bP(\Gamma\cap\{\theta_\Lambda<\sigma\wedge\tau_A\})\le x/\Lambda$. If $\vartheta<\infty$
and $\sigma=\infty$, then a.s.~$\tau_*=\infty$ and hence $\theta_\Lambda<\tau_A$ by~\autoref{lem:persistence}. Therefore,
\[
 \bP(\Gamma\cap\{\vartheta<\infty\})\le\frac x\rho+\frac x\Lambda,\qquad
 \E\Bigl[\mathbf 1_\Gamma\int_{t_\beta}^{\vartheta\wedge\tau_*}X_1(u)\rd t\Bigr]\le\frac{2x}\nu ,
\]
where the second bound follows from \eqref{eq:small-prob}, since $\vartheta\le\sigma$. On
$\Gamma\cap\{\vartheta=\infty\}$, the solution is global, and $X_0(u)\le\Lambda$ on
$[t_\beta,\infty)$.
Before $\vartheta$, we have
\[
 \norm{u_3}_\infty^{2m}\le X_0(u)^{2m}\le\Lambda^{2m-2}X_0(u)^2\le\Lambda^{2m-2}X_{-1}(u)X_1(u)\le\rho\Lambda^{2m-2}X_1(u),
\]
hence $\E[\mathbf 1_\Gamma\int_{t_\beta}^{\vartheta\wedge\tau_*}\!\int|u_3|^{2m}]\le2(2\pi)^3\rho\Lambda^{2m-2}x/\nu$, and by Markov's inequality, the probability that this integral exceeds $\eta/2$ on $\Gamma$ is at most $4(2\pi)^3\rho\Lambda^{2m-2}x/(\nu\eta)$. Collecting the exceptional events and using \eqref{eq:tail-probability}, the event in (d) has probability at least $1-\eps$ for all sufficiently large~$\beta$.
\end{proof}

\section{Nonlinear one-component noise}\label{sec:nonlinear}
Throughout this section, $k\ge4$, and $(u,\tau_*)=(u^{(\beta)},\tau_*^{(\beta)})$ is the maximal solution of~\eqref{eq:nl}. We write $C$ for numerical constants and $C_\nu$ for constants depending only on $\nu$.
They may change from line to line, but they never depend on $\beta$, $\lambda$, $k$, the reference flow $b$ introduced below, or the solution.

The noise acts radially on $v=\cR u$, so the It\^o correction of $\log(1+\norm{\nabla v}_2)$ is negative, and large if $\beta$ is large (\autoref{lem:log}). This controls the accessible part, but it does not say anything about the unforced part $\cS u$, which is not small in general. We therefore compare $\cS u$ with a reference flow $b$: The solution of the reduced system which $\cS u$ would follow in the absence of the accessible part. This reduced system is globally regular for arbitrary data (\autoref{lem:reference}). For the error $e=\cS u-b$ we only need smallness of $\norm{\nabla e}_2$, and at time zero this is achieved by choosing $b(0)$ as a Fourier truncation of~$\cS u_0$. Both pieces of information are combined in the functional $\LL_\lambda=\norm{\nabla e}_2^2+\lambda\log(1+\norm{\nabla v}_2)$: Viscosity controls the coupling at small accessible amplitudes, and the It\^o correction at large ones.

\subsection{The reference flow}\label{ss:reference}
For $b_0\in\Ran\cS\cap\Hs2$, consider
\begin{equation}\label{eq:reference}
 \del_tb+\nu Ab+\cS\cB(b,b)=0,\qquad b(0)=b_0.
\end{equation}
Since $\cS\PP$ is the orthogonal projection of the mean-zero $\rL^2$-fields onto $\Ran\cS=\{(w_h,0):\mdiv_hw_h=0\}$, writing $b=(b_h,0)$, the system \eqref{eq:reference} reads
\[
 \del_tb_h-\nu\Delta b_h+\mathbb P_h[(b_h\cdot\nabla_h)b_h]=0,\qquad\mdiv_hb_h=0,
\]
where $\mathbb P_h$ is the horizontal Leray projection, acting as the identity on the modes with $\xi_h=0$. Thus,~$b$ is a family of two-dimensional Navier--Stokes flows, parametrized by $x_3$ and coupled only through the vertical diffusion $\nu\del_3^2b_h$. This is a reduced system of 2.5D type, not the three-dimensional Navier--Stokes system with vanishing third component. Related reduced systems appear in~\cite{CZ}. We set
\[
 Y_b^h=\norm{\nabla_hb}_2^2,\qquad Y_b^v=\norm{\del_3b}_2^2,\qquad Y_b=Y_b^h+Y_b^v=\norm{\nabla b}_2^2,\qquad D_b=\norm{\Delta b}_2^2 .
\]

\begin{lem}[Reference bounds]\label{lem:reference}
For every $b_0\in\Ran\cS\cap\Hs2$, \eqref{eq:reference} has a unique global solution $b\in \rC([0,\infty);\Hs2)\cap \rL^2_{\rm loc}([0,\infty);\Hs3)$ with values in~$\Ran\cS$. For every $M<\infty$, there is $F_\nu(M)<\infty$ such that $Y_b(0)\le M$ implies
\begin{equation}\label{eq:base-H1}
 \sup_{t\ge0}e^{\nu t}Y_b(t)+\int_0^\infty e^{\nu t/2}D_b(t)\rd t\le F_\nu(M).
\end{equation}
Moreover, there is $\Phi_\nu(M)<\infty$, depending only on $\nu$ and $M$, such that
\begin{equation}\label{eq:base-coefficients}
 \mathbf A_b:=\int_0^\infty\norm{\nabla b}_3^2\rd t\le\Phi_\nu(M),\qquad\mathbf B_b:=\int_0^\infty\bigl(D_b+D_b^{2/3}\bigr)\rd t\le\Phi_\nu(M),
\end{equation}
independently of $\norm{\Delta b_0}_2$. Finally, $\sup_{t\ge0}\norm{\nabla b(t)}_3<\infty$.
\end{lem}
\begin{proof}
Local existence and uniqueness in $\rC([0,T];\Hs2)\cap \rL^2(0,T;\Hs3)$ are standard.
For instance, \autoref{thm:AV} applies with $G=0$ in the shifted triple $(\Hs3,\Hs2,\Hs1)$. The solution stays in $\Ran\cS$, since $(\Id-\cS)b$ solves the heat equation with zero initial value. It suffices to prove the a priori bounds.

\emph{Horizontal enstrophy.} Testing with $-\Delta_hb$ and using that $\Delta_hb_h$ is horizontally divergence-free, the nonlinear term becomes $\sum_{j=1}^2\int(\del_jb_h\cdot\nabla_h)b_h\cdot\del_jb_h\rd x$, after the transport term has been removed by $\mdiv_hb_h=0$. 
With $J=\nabla_hb_h$, i.e., $J_{ki}=\del_ib_k$, the integrand equals~$\tr(J^TJ^TJ)$. By
the Cayley--Hamilton theorem for $2\times2$ matrices,
$(J^T)^2=\tr(J)\,J^T-\det(J)\,I$, so that
$\tr(J^TJ^TJ)=\tr(J)\tr(J^TJ)-\det(J)\tr(J)$, which gives the pointwise identity
\[
 \sum_{i,j,k=1}^2\del_jb_i\,\del_ib_k\,\del_jb_k=\tr(J)\bigl(|J|^2-\det J\bigr)=0,
\]
because $\tr J=\mdiv_hb_h=0$. This is the two-dimensional enstrophy cancellation on each slice. Since $(\Delta b,\Delta_hb)_2=\norm{\nabla\nabla_hb}_2^2$, we obtain
\begin{equation}\label{eq:base-horizontal}
 \frac{\rd}{\rd t}Y_b^h+2\nu\norm{\nabla\nabla_hb}_2^2=0 .
\end{equation}

\emph{Vertical enstrophy.} Testing with $-\del_3^2b$, the transport term vanishes again, and the remaining nonlinear term is bounded by
\[
 \Bigl|\int(\del_3b_h\cdot\nabla_h)b_h\cdot\del_3b_h\rd x\Bigr|\le\norm{\nabla_hb}_2\norm{\del_3b}_4^2\le C(Y_b^h)^{1/2}(Y_b^v)^{1/4}\norm{\nabla\del_3b}_2^{3/2},
\]
where we used that $\del_3b$ has mean zero. Young's inequality gives
\begin{equation}\label{eq:base-vertical}
 \frac{\rd}{\rd t}Y_b^v+\nu\norm{\nabla\del_3b}_2^2\le C_\nu(Y_b^h)^2Y_b^v .
\end{equation}

\emph{$\Hs1$-bounds.} The Poincar\'e inequality gives $\norm{\nabla\nabla_hb}_2\ge\norm{\nabla_hb}_2$, so that \eqref{eq:base-horizontal} implies $Y_b^h(t)\le Y_b^h(0)e^{-2\nu t}$. Adding \eqref{eq:base-horizontal} and \eqref{eq:base-vertical}, and using $\norm{\nabla\nabla_hb}_2^2+\norm{\nabla\del_3b}_2^2=D_b\ge Y_b$, we find
\[
 \frac{\rd}{\rd t}Y_b+\nu D_b\le a_0(t)Y_b,\qquad a_0(t)=C_\nu Y_b^h(0)^2e^{-4\nu t}.
\]
Since $\int_0^\infty a_0\le C_\nu M^2/(4\nu)$, Gronwall's lemma and $D_b\ge Y_b$ give $Y_b(t)\le Me^{C_\nu M^2/(4\nu)}e^{-\nu t}$.
To obtain a weighted bound for the dissipation, put $A(t)=\int_0^ta_0\rd s$. Then, using
again $Y_b\le D_b$,
\[
 \frac{\rd}{\rd t}\bigl[e^{\nu t/2-A(t)}Y_b\bigr]
 =e^{\nu t/2-A(t)}\Bigl(Y_b'+\bigl(\tfrac\nu2-a_0\bigr)Y_b\Bigr)
 \le e^{\nu t/2-A(t)}\Bigl(-\nu D_b+\tfrac\nu2Y_b\Bigr)
 \le-\frac\nu2e^{\nu t/2-A(t)}D_b ,
\]
and integration together with $A(t)\le\int_0^\infty a_0\rd s\le C_\nu M^2/(4\nu)$ gives
$\int_0^\infty e^{\nu t/2}D_b\rd t\le\frac2\nu Me^{C_\nu M^2/(4\nu)}$. The weight is used in the last step of the proof.
This proves~\eqref{eq:base-H1}.

\emph{$\Hs2$-bounds.} Testing with $\Delta^2b$ and using the transport cancellation once more, the remaining terms of $\Delta[(b\cdot\nabla)b]$ are bounded by $C\int|\nabla b||\nabla^2b|^2\le C\norm{\nabla b}_2\norm{\nabla^2b}_4^2\le CY_b^{1/2}D_b^{1/4}\norm{\nabla\Delta b}_2^{3/2}$. Hence,
\[
 \frac{\rd}{\rd t}D_b+\nu\norm{\nabla\Delta b}_2^2\le C_\nu Y_b^2D_b,\qquad\text{so}\qquad D_b(t)\le D_b(0)\exp\Bigl(C_\nu\int_0^\infty Y_b^2\rd s\Bigr)e^{-\nu t},
\]
and $\int_0^\infty Y_b^2\rd s\le F_\nu(M)^2\int_0^\infty e^{-2\nu s}\rd s=F_\nu(M)^2/(2\nu)$,
since $Y_b(s)\le F_\nu(M)e^{-\nu s}$ by \eqref{eq:base-H1}, which has been proved in the previous step.
This gives the global solution.

\emph{Coefficients.} Since $\norm{\nabla b}_3^2\le C\norm{\nabla b}_2\norm{\nabla b}_6\le CY_b^{1/2}D_b^{1/2}$, \eqref{eq:base-H1} and the Cauchy--Schwarz inequality give
\[
 \int_0^\infty\norm{\nabla b}_3^2\rd t\le C\Bigl(\int_0^\infty e^{-\nu t/2}Y_b\rd t\Bigr)^{1/2}\Bigl(\int_0^\infty e^{\nu t/2}D_b\rd t\Bigr)^{1/2}\le C\nu^{-1/2}F_\nu(M).
\]
Moreover, $\int_0^\infty D_b\le F_\nu(M)$ and, by H\"older's inequality,
\[
 \int_0^\infty D_b^{2/3}\rd t\le\Bigl(\int_0^\infty e^{\nu t/2}D_b\rd t\Bigr)^{2/3}\Bigl(\int_0^\infty e^{-\nu t}\rd t\Bigr)^{1/3}\le\nu^{-1/3}F_\nu(M)^{2/3}.
\]
Finally, $\sup_t\norm{\nabla b}_3^2\le C\sup_t(Y_bD_b)^{1/2}<\infty$ by the $\Hs2$-bound.
\end{proof}

\subsection{Decomposition and the radial logarithm}\label{ss:decomposition}
Fix a reference flow $b$ as in \autoref{lem:reference}. It will be chosen in~\autoref{ss:choice}. Let
\begin{equation}\label{eq:quantities}
 v=\cR u,\quad e=\cS u-b,\quad r=\norm{\nabla v}_2,\quad Y=\norm{\nabla e}_2^2,\quad D_v=\norm{\Delta v}_2^2,\quad D_e=\norm{\Delta e}_2^2,
\end{equation}
so that $u=v+e+b$. Applying $\cR$ and $\cS$ to \eqref{eq:nl} and subtracting \eqref{eq:reference}, we obtain
\begin{align}
 \rd v+[\nu Av+\cR\cB(u)]\rd t&=\beta r^kv\rd W_t,\label{eq:v}\\
 \del_te+\nu Ae+\cS[\cB(u)-\cB(b,b)]&=0.\label{eq:e}
\end{align}
The term $\cR\cB(b,b)$ is kept in \eqref{eq:v}: The reference flow is only used in the proof, and the equation is not modified. Since $\cS u$ solves an equation without stochastic integral, on every interval $[0,T]$ with~$T<\tau_*$, we have $e\in \rW^{1,2}(0,T;\rL^2_\sigma)\cap \rL^2(0,T;\Hs2)$ a.s., and $Y$ is absolutely continuous with
\begin{equation}\label{eq:Y-identity}
 \frac{\rd}{\rd t}Y+2\nu D_e=-2\bigl(\cB(u)-\cB(b,b),Ae\bigr)_2,
\end{equation}
where we used that $Ae\in\Ran\cS$. For the accessible part we write
\[
 T_v=(\cB(u),Av)_2=-(\cB(u),\Delta v)_2,
\]
so that \eqref{eq:energy-v} reads $\rd(r^2)=(-2\nu D_v-2T_v+\beta^2r^{2k+2})\rd t+2\beta r^{k+2}\rd W$. Let
\begin{equation}\label{eq:weights}
 \psi(r)=\log(1+r),\qquad a(r)=\frac1{r(1+r)},\qquad h_k(r)=\frac{r^{2k+2}}{(1+r)^2}.
\end{equation}

\begin{lem}[Radial logarithm]\label{lem:log}
A.s., for all $t<\tau_*$, it holds that
\begin{equation}\label{eq:log}
 \psi(r(t))+\int_0^t\Bigl[\nu a(r)D_v+\frac{\beta^2}2h_k(r)\Bigr]\rd s=\psi(r(0))-\int_0^ta(r)T_v\rd s+\beta\int_0^t\frac{r^{k+1}}{1+r}\rd W_s,
\end{equation}
where $a(r)D_v$ and $a(r)T_v$ are set to zero at times with $v=0$, and all integrals converge absolutely.
\end{lem}
\begin{proof}
Where $r>0$, this is It\^o's formula for $E\mapsto\log(1+\sqrt E)$ applied to $E=r^2$: The first two derivatives are $\frac1{2r(1+r)}$ and $-\frac{1+2r}{4r^3(1+r)^2}$, and the noise contributes the drift
\[
 \frac{\beta^2r^{2k+2}}{2r(1+r)}-\frac{(1+2r)\,4\beta^2r^{2k+4}}{8r^3(1+r)^2}=-\frac{\beta^2}2\,\frac{r^{2k+2}}{(1+r)^2}.
\]
The function $\log(1+\sqrt E)$ is not $\rC^2$ at $E=0$; the justification through $r=0$ by regularization is given in~\autoref{app:log}.
\end{proof}

The choice of $\log(1+r)$ instead of $\log(1+r^2)$ matters: The It\^o drift of $\log(1+r^2)$ is $\beta^2r^{2k+2}(1-r^2)/(1+r^2)^2$, which is positive for $r<1$, whereas $\log(1+r)$ has a negative noise drift at all amplitudes.

\subsection{Convection estimates}\label{ss:convection-nl}
\begin{lem}\label{lem:convection-nl}
There are a numerical constant $C_0$ and constants $C_\nu$, $Y_*=Y_*(\nu)>0$ such that
\begin{equation}\label{eq:Tv}
 |T_v|\le C_0\Bigl[r^{3/2}D_v^{3/4}+\bigl(D_e^{1/2}+D_b^{1/2}\bigr)r^{3/2}D_v^{1/4}+(D_e+D_b)\,r\Bigr],
\end{equation}
and, at all times with $Y<Y_*$,
\begin{equation}\label{eq:error-energy}
 \frac{\rd}{\rd t}Y+\nu D_e\le C_\nu\norm{\nabla b}_3^2\,(Y+r^2)+C_\nu\Bigl[r^4Y+\bigl((Y+Y_b)r+r^3\bigr)D_v^{1/2}\Bigr].
\end{equation}
\end{lem}
\begin{proof}
Write $w=e+b=\cS u$, so that $u=v+w$, and recall $\norm v_6\le Cr$, $\norm{\nabla v}_3\le Cr^{1/2}D_v^{1/4}$ and $\norm{\nabla v}_{12/5}\le Cr^{3/4}D_v^{1/8}$, which follow from interpolation between $\rL^2$ and $\rL^6$, with analogous bounds for $e$ and $b$, and $\norm{\Delta w}_2\le D_e^{1/2}+D_b^{1/2}$.

\emph{Proof of~\eqref{eq:Tv}.} Since $Av$ is divergence-free, $T_v=((u\cdot\nabla)u,-\Delta v)_2$, and we split $(u\cdot\nabla)u$ into four terms. First,
\[
 |((v\cdot\nabla)v,\Delta v)_2|\le\norm v_6\norm{\nabla v}_3\norm{\Delta v}_2\le Cr^{3/2}D_v^{3/4}.
\]
Second, integrating by parts and using $\mdiv w=0$, $((w\cdot\nabla)v,-\Delta v)_2=\sum_j((\del_jw\cdot\nabla)v,\del_jv)_2$. The transported derivative cancels, and the term is bounded by $\norm{\nabla w}_6\norm{\nabla v}_{12/5}^2\le C\norm{\Delta w}_2r^{3/2}D_v^{1/4}$. Third, $((v\cdot\nabla)w,-\Delta v)_2=\sum_j[((\del_jv\cdot\nabla)w,\del_jv)_2+((v\cdot\nabla)\del_jw,\del_jv)_2]$, which is bounded by the same expression plus $\norm v_6\norm{\nabla^2w}_2\norm{\nabla v}_3\le C\norm{\Delta w}_2r^{3/2}D_v^{1/4}$. Finally, $((w\cdot\nabla)w,-\Delta v)_2=(\nabla[(w\cdot\nabla)w],\nabla v)_2$, and $\norm{\nabla[(w\cdot\nabla)w]}_2\le\norm{\nabla w}_4^2+\norm w_\infty\norm{\nabla^2w}_2\le C\norm{\Delta w}_2^2$, so that
\begin{equation}\label{eq:regeneration}
 |((w\cdot\nabla)w,\Delta v)_2|\le C\norm{\Delta w}_2^2\,r\le C(D_e+D_b)\,r .
\end{equation}

\emph{Proof of~\eqref{eq:error-energy}.} By \eqref{eq:Y-identity}, $\frac12\frac{\rd}{\rd t}Y+\nu D_e=-(\cB(u)-\cB(b,b),Ae)_2$, and
\[
 \cB(u)-\cB(b,b)=\cB(e,e)+\bigl[\cB(v,e)+\cB(e,v)+\cB(v,v)\bigr]+\bigl[\cB(e,b)+\cB(b,e)\bigr]+\bigl[\cB(b,v)+\cB(v,b)\bigr].
\]
(i) By \eqref{eq:conv-commutator} (with $e$ in place of $u$), $|(\cB(e,e),Ae)_2|\le\norm{\nabla e}_3^3\le CY^{3/4}D_e^{3/4}\le CY^{1/2}D_e$, where we used $Y\le D_e$ (Poincar\'e). For $Y<Y_*:=(\nu/(4C))^2$ this is at most $\frac\nu4D_e$.

(ii) The remaining terms which are quadratic in $(v,e)$ satisfy
\begin{gather*}
 |(\cB(v,e),Ae)_2|\le CrY^{1/4}D_e^{3/4},\qquad|(\cB(e,v),Ae)_2|\le CY^{1/2}r^{1/2}D_v^{1/4}D_e^{1/2},\\
 |(\cB(v,v),Ae)_2|\le Cr^{3/2}D_v^{1/4}D_e^{1/2},
\end{gather*}
by H\"older's inequality with $\rL^6\times \rL^3\times \rL^2$. By Young's inequality, their sum is at most $\frac\nu8D_e+C_\nu[r^4Y+(Yr+r^3)D_v^{1/2}]$.

(iii) For the terms which are linear in $e$, the transport cancellation and H\"older's inequality give
\begin{align*}
 |(\cB(b,e),Ae)_2|&=\Bigl|\sum_j((\del_jb\cdot\nabla)e,\del_je)_2\Bigr|\le\norm{\nabla b}_3\norm{\nabla e}_2\norm{\nabla e}_6\le C\norm{\nabla b}_3Y^{1/2}D_e^{1/2},\\
 |(\cB(e,b),Ae)_2|&\le\norm e_6\norm{\nabla b}_3\norm{\Delta e}_2\le C\norm{\nabla b}_3Y^{1/2}D_e^{1/2}.
\end{align*}
Their sum is at most $\frac\nu{16}D_e+C_\nu\norm{\nabla b}_3^2Y$.

(iv) For the terms which are linear in $v$,
\begin{align*}
 |(\cB(b,v),Ae)_2|&\le\norm b_6\norm{\nabla v}_3\norm{\Delta e}_2\le CY_b^{1/2}r^{1/2}D_v^{1/4}D_e^{1/2},\\
 |(\cB(v,b),Ae)_2|&\le\norm v_6\norm{\nabla b}_3\norm{\Delta e}_2\le Cr\norm{\nabla b}_3D_e^{1/2}.
\end{align*}
Their sum is at most $\frac\nu{16}D_e+C_\nu[Y_brD_v^{1/2}+\norm{\nabla b}_3^2r^2]$.

Adding (i)--(iv), and noting that the fractions of $D_e$ used there sum to
$\frac\nu4+\frac\nu8+\frac\nu{16}+\frac\nu{16}=\frac\nu2$, for~$Y<Y_*$, we obtain
\[
 \frac12\frac{\rd}{\rd t}Y+\frac\nu2D_e\le C_\nu\norm{\nabla b}_3^2\,(Y+r^2)
 +C_\nu\Bigl[r^4Y+\bigl((Y+Y_b)r+r^3\bigr)D_v^{1/2}\Bigr],
\]
which gives \eqref{eq:error-energy} after multiplication by~$2$ (and enlarging $C_\nu$).
\end{proof}

The bound \eqref{eq:regeneration} is the reason why the accessible forcing regenerated by the unforced part can be paid for: After multiplication by $a(r)$ its contribution is at most $C(D_e+D_b)$, independently of $r$, and $D_e$ is controlled by the viscous dissipation of the error.

From now on, we fix $C_\nu\ge1$ such that \autoref{lem:convection-nl} holds, and put
\begin{gather*}
 a_b(t)=C_\nu\norm{\nabla b(t)}_3^2,\qquad f_b(t)=C_\nu\bigl(D_b(t)+D_b(t)^{2/3}\bigr),\\
 A_b=\int_0^\infty a_b\rd t=C_\nu\mathbf A_b,\qquad B_b=\int_0^\infty f_b\rd t=C_\nu\mathbf B_b .
\end{gather*}
By \autoref{lem:reference}, $A_b+B_b\le2C_\nu\Phi_\nu(M)$ whenever $Y_b(0)\le M$, and $\sup_ta_b<\infty$ for each fixed reference.

\subsection{The joint inequality}\label{ss:joint}
For $\lambda\in(0,1]$, let
\begin{equation}\label{eq:L}
 \LL_\lambda=Y+\lambda\psi(r),\qquad\rd M_\lambda=\lambda\beta\frac{r^{k+1}}{1+r}\rd W_t,\qquad\rd\langle M_\lambda\rangle=\lambda^2\beta^2h_k(r)\rd t.
\end{equation}
Below, inequalities between differentials are understood in integrated form on the time interval on which they are asserted. We first record the Young inequalities which produce the powers of~$\lambda$.

\begin{lem}[Weighted Young inequalities]\label{lem:young}
For every $\eps\in(0,1]$, there is $C_{\nu,\eps}$ such that for all~$r>0$, $\lambda\in(0,1]$ and $D_v,D_e,D_b,Y,Y_b\ge0$, we have
\begin{align}
 C_0\lambda a(r)r^{3/2}D_v^{3/4}&\le\eps\lambda\nu a(r)D_v+C_{\nu,\eps}\lambda\frac{r^5}{1+r},\label{eq:young-self}\\
 C_0\lambda a(r)D_e^{1/2}r^{3/2}D_v^{1/4}&\le\eps\lambda\nu a(r)D_v+\eps\nu D_e+C_{\nu,\eps}\lambda^3\Bigl(\frac r{1+r}\Bigr)^3,\label{eq:young-error}\\
 C_0\lambda a(r)D_b^{1/2}r^{3/2}D_v^{1/4}&\le\eps\lambda\nu a(r)D_v+C_{\nu,\eps}\lambda D_b^{2/3},\label{eq:young-base}\\
 C_0\lambda a(r)(D_e+D_b)\,r&\le C_0\lambda(D_e+D_b),\label{eq:young-regeneration}\\
 C_\nu\bigl((Y+Y_b)r+r^3\bigr)D_v^{1/2}&\le\eps\lambda\nu a(r)D_v+\frac{C_{\nu,\eps}}\lambda\bigl((Y+Y_b)^2r^3+r^7\bigr)(1+r).\label{eq:young-mixed}
\end{align}
\end{lem}
\begin{proof}
For \eqref{eq:young-self}, write the left-hand side as $(\lambda aD_v)^{3/4}\cdot C_0(\lambda a)^{1/4}r^{3/2}$ and use Young's inequality with exponents $\frac43$ and $4$. The remainder is $C_{\nu,\eps}\lambda ar^6=C_{\nu,\eps}\lambda r^5/(1+r)$. For \eqref{eq:young-error}, Young's inequality with exponents $4$ and $\frac43$ in $D_v$ gives first $\eps\lambda\nu aD_v+C_{\nu,\eps}\lambda ar^2D_e^{2/3}=\eps\lambda\nu aD_v+C_{\nu,\eps}\lambda\frac r{1+r}D_e^{2/3}$, and Young's inequality with exponents $\frac32$ and $3$ then gives the claim. The same first step and $r/(1+r)\le1$ give~\eqref{eq:young-base}. Inequality \eqref{eq:young-regeneration} is $a(r)r=(1+r)^{-1}\le1$. Finally, \eqref{eq:young-mixed} follows from Young's inequality with exponents $2$ and $2$, $1/a(r)=r(1+r)$, and $((Y+Y_b)r+r^3)^2\le2((Y+Y_b)^2r^2+r^6)$.
\end{proof}

\begin{prop}[Joint inequality]\label{prop:joint}
There is $\lambda_0=\lambda_0(\nu)\in(0,1]$ such that for $0<\lambda\le\lambda_0$, as long as $Y<Y_*$, it holds that
\begin{equation}\label{eq:preclosed}
 \rd\LL_\lambda+\Bigl[\frac\nu2D_e+\frac\lambda2\nu a(r)D_v+\frac\lambda2\beta^2h_k(r)\Bigr]\rd t\le\bigl[a_b\LL_\lambda+\lambda f_b+R_{\lambda,b}(r)\bigr]\rd t+\rd M_\lambda,
\end{equation}
where
\begin{equation}\label{eq:residual}
 0\le R_{\lambda,b}(r)\le C_{\nu,b}\bigl[\lambda^{-1}(r^3+r^8)+r^2\bigr],
\end{equation}
and $C_{\nu,b}$ depends only on $\nu$, $\sup_tY_b$ and~$\sup_ta_b$. The coefficients $a_b$ and $f_b$ do not depend on~$\lambda$.
\end{prop}
\begin{proof}
By \autoref{lem:log} and \eqref{eq:error-energy}, as long as $Y<Y_*$, we get
\[
 \rd\LL_\lambda\le\Bigl[-\nu D_e+a_b(Y+r^2)+C_\nu\bigl(r^4Y+((Y+Y_b)r+r^3)D_v^{1/2}\bigr)-\lambda\nu aD_v-\frac\lambda2\beta^2h_k+\lambda a|T_v|\Bigr]\rd t+\rd M_\lambda .
\]
We bound $\lambda a|T_v|$ by \eqref{eq:Tv} and \eqref{eq:young-self}--\eqref{eq:young-regeneration}, and the mixed term by \eqref{eq:young-mixed}, all with $\eps=\frac18$. This uses up at most $\frac12\lambda\nu aD_v$ and $\frac18\nu D_e$, and produces the additional term $C_0\lambda D_e$. We choose $\lambda_0\le\nu/(8C_0)$, so that the total coefficient of $D_e$ remains at most~$-\frac\nu2$. The terms $C_0\lambda D_b+C_{\nu,1/8}\lambda D_b^{2/3}$ are bounded by $\lambda f_b$ after enlarging $C_\nu$ if necessary. Since $\psi\ge0$, we have $a_bY\le a_b\LL_\lambda$. The remaining terms are
\[
 R_{\lambda,b}(r)=a_br^2+C_\nu r^4Y+\frac{C_{\nu}}\lambda\bigl((Y+Y_b)^2r^3+r^7\bigr)(1+r)+C_\nu\lambda\frac{r^5}{1+r}+C_\nu\lambda^3\Bigl(\frac r{1+r}\Bigr)^3 .
\]
Using $Y<Y_*$, $\lambda\le1$ as well as the elementary inequalities $r^4,r^5,r^7\le r^3+r^8$ for $r\ge0$, we obtain~\eqref{eq:residual}.
\end{proof}

\begin{lem}[Scalar absorption]\label{lem:absorb}
There is $C_{\nu,k,b}$ such that, if $0<\lambda\le\lambda_0$ and
\begin{equation}\label{eq:beta-lambda}
 \beta^2\lambda^{2k+1}\ge C_{\nu,k,b},
\end{equation}
then $R_{\lambda,b}(r)\le\frac14\lambda\nu\frac r{1+r}+\frac14\lambda\beta^2h_k(r)$ for all $r\ge0$.
\end{lem}
\begin{proof}
Let $C_b=C_{\nu,b}$ be the constant in \eqref{eq:residual} and $d_0=\min\{1,\nu/(24C_b)\}$. We split the range of $r$ at $d_0\lambda$ and at~$1$.

If $r\le d_0\lambda\le1$, then $r/(1+r)\ge r/2$, and
\[
 C_b\bigl[\lambda^{-1}(r^3+r^8)+r^2\bigr]\le C_b\bigl[2\lambda^{-1}r^3+r^2\bigr]\le C_b\bigl(2d_0^2+d_0\bigr)\lambda r\le\frac{\lambda\nu}8r\le\frac{\lambda\nu}4\,\frac r{1+r}.
\]

If $d_0\lambda\le r\le1$, then $h_k(r)\ge r^{2k+2}/4$, and $C_b[\lambda^{-1}(r^3+r^8)+r^2]\le2C_b\lambda^{-1}r^3+C_br^2$. The two ratios with respect to $\frac1{16}\lambda\beta^2r^{2k+2}$ are
\[
 \frac{32C_b\lambda^{-1}r^3}{\lambda\beta^2r^{2k+2}}=\frac{32C_b}{\lambda^2\beta^2r^{2k-1}}\le\frac{32C_bd_0^{-(2k-1)}}{\beta^2\lambda^{2k+1}},\qquad
 \frac{16C_br^2}{\lambda\beta^2r^{2k+2}}=\frac{16C_b}{\lambda\beta^2r^{2k}}\le\frac{16C_bd_0^{-2k}}{\beta^2\lambda^{2k+1}},
\]
and both are at most $\frac12$ if $\beta^2\lambda^{2k+1}\ge64C_bd_0^{-2k}$.

If $r\ge1$, then $h_k(r)\ge r^{2k}/4$ and, since $2k\ge8$ and $\lambda\le1$,
\[
 C_b\bigl[\lambda^{-1}(r^3+r^8)+r^2\bigr]\le3C_b\lambda^{-1}r^8\le3C_b\lambda^{-1}r^{2k}.
\]
This is at most $\frac1{16}\lambda\beta^2r^{2k}$ if $\lambda^2\beta^2\ge48C_b$, which follows from $\beta^2\lambda^{2k+1}\ge48C_b$.

Thus \eqref{eq:beta-lambda} with $C_{\nu,k,b}=64C_bd_0^{-2k}$ suffices.
\end{proof}

This is the only place where $k\ge4$ is used: The It\^o correction must dominate the residual $r^8$ at large accessible amplitudes.

\begin{cor}[Closed inequality]\label{cor:closed}
Under \eqref{eq:beta-lambda} and $0<\lambda\le\lambda_0$, as long as $Y<Y_*$, we have
\begin{equation}\label{eq:closed}
 \rd\LL_\lambda+\Bigl[\frac\nu2D_e+\frac\lambda4\nu a(r)D_v+\frac\lambda4\beta^2h_k(r)\Bigr]\rd t\le\bigl[a_b\LL_\lambda+\lambda f_b\bigr]\rd t+\rd M_\lambda .
\end{equation}
\end{cor}
\begin{proof}
By the Poincar\'e inequality, $D_v\ge r^2$, and therefore $a(r)D_v\ge r/(1+r)$. Thus, the bound of \autoref{lem:absorb} is absorbed by the dissipative terms in~\eqref{eq:preclosed}.
\end{proof}

The constants in \eqref{eq:closed} do not depend on any stochastic localization of the strong norms. The reference $b$ is fixed before $\lambda$ and $\beta$ are chosen.

\subsection{Stopping and the exit estimate}\label{ss:exit}
Define the deterministic weights
\[
 \varrho_b(t)=\exp\Bigl(-\int_0^ta_b(s)\rd s\Bigr),\qquad\varrho_*=e^{-A_b}>0,\qquad\varrho_*\le\varrho_b\le1,
\]
the exit time $\sigma=\inf\{t\in[0,\tau_*):Y(t)\ge Y_*\}$, where $\inf\emptyset=\infty$, and
\[
 r_0=\norm{\nabla\cR u_0}_2,\qquad Y_0=\norm{\nabla(\cS u_0-b(0))}_2^2,\qquad\ell_\lambda=Y_0+\lambda\bigl[\log(1+r_0)+B_b\bigr].
\]

\begin{prop}\label{prop:exit}
Assume \eqref{eq:beta-lambda}, $0<\lambda\le\lambda_0$ and $Y_0<Y_*$. Then we have
\begin{enumerate}[(a)]
 \item $\tau_*=\infty$ a.s.\ on $\{\sigma=\infty\}$;
 \item $\displaystyle\E\int_0^{\sigma\wedge\tau_*}\varrho_b\Bigl[\frac\nu2D_e+\frac\lambda4\nu a(r)D_v+\frac\lambda4\beta^2h_k(r)\Bigr]\rd t\le\ell_\lambda$;
 \item $\displaystyle\bP(\sigma<\infty)\le\min\Bigl\{1,\ \sqrt{1+r_0}\,e^{B_b/2}\exp\Bigl[-\frac{\varrho_*Y_*-Y_0}{2\lambda}\Bigr]\Bigr\}$.
\end{enumerate}
\end{prop}
\begin{proof}
\emph{Weighted inequality.} Since $\varrho_b$ is deterministic and absolutely continuous
with $\varrho_b'=-a_b\varrho_b$, the product rule gives
$\rd(\varrho_b\LL_\lambda)=\varrho_b\rd\LL_\lambda-a_b\varrho_b\LL_\lambda\rd t$, without a
covariation term. Multiplying~\eqref{eq:closed} by~$\varrho_b>0$, the terms
$\pm a_b\varrho_b\LL_\lambda$ cancel, and integration over $[0,t]$ yields the bound with
$\LL_\lambda(0)+\lambda\int_0^t\varrho_bf_b\rd s$ on the right-hand side. Since
$\LL_\lambda(0)=Y_0+\lambda\log(1+r_0)$ and $\lambda\int_0^t\varrho_bf_b\rd s\le\lambda B_b$, as $0<\varrho_b\le1$ and~$f_b\ge0$, for $t<\sigma\wedge\tau_*$, we obtain
\begin{equation}\label{eq:weighted}
 \varrho_b(t)\LL_\lambda(t)+\int_0^t\varrho_b\Bigl[\frac\nu2D_e+\frac\lambda4\nu aD_v+\frac\lambda4\beta^2h_k\Bigr]\rd s\le\ell_\lambda+\widetilde M_t,\qquad\widetilde M_t=\int_0^t\varrho_b\rd M_\lambda .
\end{equation}

\emph{Proof of (a) and (b).} Let $\tau_n$ be as in the proof of~\autoref{lem:Xminus}. On $[0,\tau_n]$, we have $r\le n$, so the integrand of $\widetilde M$ is bounded, and taking expectations in \eqref{eq:weighted} at $t\wedge\sigma\wedge\tau_n$ gives
\[
 \E\bigl[\varrho_b\LL_\lambda(t\wedge\sigma\wedge\tau_n)\bigr]+\E\int_0^{t\wedge\sigma\wedge\tau_n}\varrho_b\Bigl[\frac\nu2D_e+\frac\lambda4\nu aD_v+\frac\lambda4\beta^2h_k\Bigr]\rd s\le\ell_\lambda .
\]
Monotone convergence ($n\to\infty$, then $t\to\infty$) yields (b). 
For (a), fix $T<\infty$, let $N>r_0$, $\theta_N=\inf\{t<\tau_*:r(t)\ge N\}$ and
$\vartheta_n=\theta_N\wedge T\wedge\sigma\wedge\tau_n$. Evaluating \eqref{eq:weighted} at
$\vartheta_n$ and taking expectations as above gives
$\E[\varrho_b(\vartheta_n)\LL_\lambda(\vartheta_n)]\le\ell_\lambda$. On
$E_N^n=\{\theta_N<\sigma\wedge\tau_*\wedge T,\ \theta_N<\tau_n\}$, we have $\vartheta_n=\theta_N$
and $r(\theta_N)=N$, hence
$\varrho_b(\vartheta_n)\LL_\lambda(\vartheta_n)\ge\varrho_*\lambda\log(1+N)$. Since the left-hand
side is nonnegative, this yields
$\varrho_*\lambda\log(1+N)\,\bP(E_N^n)\le\ell_\lambda$, and letting $n\to\infty$, so that
$E_N^n\uparrow\{\theta_N<\sigma\wedge\tau_*\wedge T\}$, we obtain
\[
 \bP\bigl(\theta_N<\sigma\wedge\tau_*\wedge T\bigr)\le\frac{\ell_\lambda}{\varrho_*\lambda\log(1+N)}\longrightarrow0\qquad(N\to\infty).
\]
Hence, a.s.\ on $\{\sigma=\infty,\ \tau_*\le T\}$, $r$ is bounded on $[0,\tau_*)$. Moreover, $Y<Y_*$ there and $Y_b$ is bounded, so $\sup_{t<\tau_*}\norm{\nabla u(t)}_2\le\sup_{t<\tau_*}(r+Y^{1/2}+Y_b^{1/2})<\infty$. By \eqref{eq:blowup} this event is a null set, and $T\to\infty$ gives (a).

\emph{Proof of (c).} By \eqref{eq:L}, $\rd\langle\widetilde M\rangle=\varrho_b^2\lambda^2\beta^2h_k\rd t\le\varrho_b\lambda^2\beta^2h_k\rd t$, so $\frac\lambda4\beta^2\int_0^t\varrho_bh_k\rd s\ge\frac1{4\lambda}\langle\widetilde M\rangle_t$. Dropping the other nonnegative terms in \eqref{eq:weighted},
\begin{equation}\label{eq:mart-bound}
 \varrho_b(t)\LL_\lambda(t)+\frac1{4\lambda}\langle\widetilde M\rangle_t\le\ell_\lambda+\widetilde M_t,\qquad t<\sigma\wedge\tau_* .
\end{equation}
For every continuous local martingale $N$ with $N_0=0$ and every $h>0$,
\begin{equation}\label{eq:exp-martingale}
 \bP\Bigl(\sup_{t\ge0}\Bigl[N_t-\frac1{4\lambda}\langle N\rangle_t\Bigr]\ge h\Bigr)\le e^{-h/(2\lambda)}.
\end{equation}
Indeed, with $\alpha=1/(2\lambda)$, the process $Z=\exp(\alpha N-\frac{\alpha^2}2\langle N\rangle)$ is a nonnegative continuous local martingale, hence a supermartingale with $Z_0=1$. Since $\frac{\alpha^2}2=\frac\alpha{4\lambda}$, the event in \eqref{eq:exp-martingale} is $\{\sup_tZ_t\ge e^{\alpha h}\}$, whose probability is at most $e^{-\alpha h}$ by the maximal inequality for nonnegative supermartingales. No Novikov condition is needed. 
On $\{\sigma<\infty\}$, we have $\sigma<\tau_*$ and, by continuity of $Y$, $Y(\sigma)=Y_*$,
so $\varrho_b(\sigma)\LL_\lambda(\sigma)\ge\varrho_*Y(\sigma)=\varrho_*Y_*$. Since all terms
in \eqref{eq:mart-bound} are continuous in $t$, the inequality extends to $t=\sigma$, and we
obtain on $\{\sigma<\infty\}$
\[
 \widetilde M_\sigma-\frac1{4\lambda}\langle\widetilde M\rangle_\sigma\ge
 \varrho_*Y_*-\ell_\lambda=:h .
\]
Let $N_t=\widetilde M_{t\wedge\sigma}$. By (a), $\sigma\wedge\tau_*=\sigma$ a.s., so $N$ is a
continuous local martingale on $[0,\infty)$ with~$N_0=0$, localized by $(\tau_n)$, and
$\{\sigma<\infty\}\subset\{\sup_{t\ge0}[N_t-\frac1{4\lambda}\langle N\rangle_t]\ge h\}$.
If $h>0$, \eqref{eq:exp-martingale} and
$\ell_\lambda=Y_0+\lambda\log(1+r_0)+\lambda B_b$ yield
\[
 \bP(\sigma<\infty)\le e^{-h/(2\lambda)}
 =\exp\Bigl[-\frac{\varrho_*Y_*-Y_0}{2\lambda}\Bigr]\,
 e^{\frac12\log(1+r_0)+\frac12B_b}
 =\sqrt{1+r_0}\,e^{B_b/2}\exp\Bigl[-\frac{\varrho_*Y_*-Y_0}{2\lambda}\Bigr].
\]
If $h\le0$, the right-hand side is at least one, and (c) holds trivially.
\end{proof}

\subsection{Choice of the reference flow}\label{ss:choice}
Let $M=\norm{\nabla\cS u_0}_2^2$, and for $N\in\N$, let $b_0^{(N)}=P_{\le N}\cS u_0$ be the Fourier truncation of $\cS u_0$ to $|\xi|\le N$. It belongs to $\Ran\cS\cap\Hs2$, and $Y_b(0)\le M$ for every~$N$.
By \autoref{lem:reference}, $A_b=C_\nu\mathbf A_b\le C_\nu\Phi_\nu(M)$ uniformly in $N$, hence
\[
 \varrho_*\ge\underline\varrho_\nu(M):=e^{-C_\nu\Phi_\nu(M)}>0\qquad\text{uniformly in }N.
\]
Since $P_{\le N}\cS u_0\to\cS u_0$ in $\Hs1$, we can fix $N$ such that
\begin{equation}\label{eq:initial-error}
 Y_0=\norm{\nabla(\cS u_0-P_{\le N}\cS u_0)}_2^2\le\tfrac12\underline\varrho_\nu(M)Y_* .
\end{equation}
Then $\varrho_*Y_*-Y_0\ge\frac12\underline\varrho_\nu(M)Y_*>0$, whatever the size of~$M$. The higher norms of the chosen reference, such as $\norm{\Delta b_0}_2$, may be large, but they are fixed once $N$ is fixed. They only enter the constant~$C_{\nu,k,b}$ in~\eqref{eq:beta-lambda}. In particular, this is not a smallness assumption on the data: The reference may be large, and only the initial error is small. With this choice, \autoref{prop:exit}(c) gives
\begin{equation}\label{eq:exit-final}
 \bP(\sigma<\infty)\le\sqrt{1+r_0}\,e^{C_\nu\Phi_\nu(M)/2}\exp\Bigl[-\frac{\underline\varrho_\nu(M)Y_*}{4\lambda}\Bigr],
\end{equation}
for all $0<\lambda\le\lambda_0$ and all $\beta$ satisfying~\eqref{eq:beta-lambda}.

\subsection{The component integral and the proof of the theorem}\label{ss:component}
From \autoref{prop:exit}(b), $\varrho_b\ge\varrho_*$ and the bound
$a(r)D_v\ge a(r)r^2=r/(1+r)$, which follows from the Poincar\'e inequality $D_v\ge r^2$ as in
the proof of \autoref{cor:closed}, we obtain
\begin{equation}\label{eq:two-integrals}
 \E\int_0^{\sigma\wedge\tau_*}\frac r{1+r}\rd t\le\frac{4\ell_\lambda}{\varrho_*\lambda\nu},\qquad\E\int_0^{\sigma\wedge\tau_*}h_k(r)\rd t\le\frac{4\ell_\lambda}{\varrho_*\lambda\beta^2}.
\end{equation}
For $1<m\le2k$, $0<\delta\le1$, and $r\ge0$, we have the scalar splitting
\begin{equation}\label{eq:scalar-split}
 r^m\le2\delta^{m-1}\frac r{1+r}+4\delta^{-(2k+2-m)}h_k(r).
\end{equation}
Indeed, for $r\le\delta$, we have $r^m\le\delta^{m-1}r\le2\delta^{m-1}r/(1+r)$, while for
$r>\delta$, we have $r^m/h_k(r)=r^{m-2k-2}(1+r)^2\le4\delta^{-(2k+2-m)}$, as one sees by
treating $r\le1$ and $r>1$ separately and using $m\le2k$.
Consequently,
\begin{equation}\label{eq:moment}
 \E\int_0^{\sigma\wedge\tau_*}r^m\rd t\le\frac{4\ell_\lambda}{\varrho_*\lambda}
 \biggl(\frac{2\delta^{m-1}}\nu+\frac{4\delta^{-(2k+2-m)}}{\beta^2}\biggr).
\end{equation}

\begin{proof}[Proof of \autoref{thm:nonlinear}]
Fix $\eta,\eps>0$ and choose the reference flow as in~\autoref{ss:choice}. Choose $\lambda\in(0,\lambda_0]$ so small that the right-hand side of \eqref{eq:exit-final} is at most $\eps/2$, and keep $\lambda$ fixed. Since $u_3=v_3$ has mean zero, $\int_{\T^3}|u_3|^6\rd x\le C_Sr^6$ with a numerical constant~$C_S$. As $k\ge4$, we may take $m=6$ in \eqref{eq:moment}: We first choose $\delta\in(0,1]$ such that $C_S\frac{4\ell_\lambda}{\varrho_*\lambda}\frac{2\delta^5}\nu\le\frac{\eps\eta}4$, and then $\beta_0$ such that, for $\beta\ge\beta_0$, \eqref{eq:beta-lambda} holds and $C_S\frac{4\ell_\lambda}{\varrho_*\lambda}\frac{4\delta^{-(2k-4)}}{\beta^2}\le\frac{\eps\eta}4$. 
Combining these choices with the bound $\int_{\T^3}|u_3|^6\rd x\le C_Sr^6$ and
\eqref{eq:moment}, we obtain
$\E\int_0^{\sigma\wedge\tau_*}\!\int_{\T^3}|u_3|^6\rd x\rd t\le\eps\eta/2$, and Markov's
inequality yields
\[
 \bP\Bigl(\int_0^{\sigma\wedge\tau_*}\!\!\int_{\T^3}|u_3|^6\rd x\rd t>\eta\Bigr)\le\frac\eps2 .
\]
On $\{\sigma=\infty\}$, we have $\tau_*=\infty$ by \autoref{prop:exit}(a), and the stopped integral is the full integral. 
Since $\bP(\sigma<\infty)\le\eps/2$ by \eqref{eq:exit-final} and the choice of $\lambda$
(recall that \eqref{eq:beta-lambda} holds for $\beta\ge\beta_0$), a union bound shows that
the event in \eqref{eq:main-nl} has probability at least $1-\eps$ for all $\beta\ge\beta_0$.

For the rate, keep the reference fixed, put $\gamma_k=2/(2k+1)$ and $\lambda=\Lambda_b\beta^{-\gamma_k}$ with $\Lambda_b^{2k+1}=C_{\nu,k,b}$. Then $\beta^2\lambda^{2k+1}=C_{\nu,k,b}$ for every $\beta$, and $\lambda\le\lambda_0$ for $\beta\ge\beta_1$. By \autoref{prop:exit}(a) and \eqref{eq:exit-final},
\[
 \bP(\tau_*<\infty)\le\bP(\sigma<\infty)\le\sqrt{1+r_0}\,e^{C_\nu\Phi_\nu(M)/2}\exp\Bigl[-\frac{\underline\varrho_\nu(M)Y_*}{4\Lambda_b}\,\beta^{\gamma_k}\Bigr],\qquad\beta\ge\beta_1,
\]
and enlarging $C$ covers $\beta\in[1,\beta_1]$.
\end{proof}

\begin{rem}\label{rem:rates}
With the same choice of $\lambda$ and with $\delta=\beta^{-\gamma_k}$, \eqref{eq:moment} gives, for $2<m\le2k$, the estimate
\[
 \E\int_0^{\sigma\wedge\tau_*}r^m\rd t\le C_{\nu,k,b,m}\Bigl[Y_0\,\beta^{-\gamma_k(m-2)}+\bigl(\log(1+r_0)+B_b\bigr)\beta^{-\gamma_k(m-1)}\Bigr],
\]
since $\gamma_k(2k+1)=2$. For $k=4$ and $m=6$ the powers are $\beta^{-8/9}$ and~$\beta^{-10/9}$. These are stopped moments.
We make no claim about expectations on the exceptional event. For $k\ge8$, one may take~$m=16$, which controls $\int\norm{u_3}_6^{16}\rd t$ up to $\sigma$, the mixed norm appearing in one of the criteria of \cite{KZ:06}. The proof above does not use such a criterion and works for $k=4$.
\end{rem}

\begin{rem}\label{rem:zero-reference}
Choosing $b\equiv0$ gives $\varrho_*=1$, $B_b=0$ and $Y_0=\norm{\nabla\cS u_0}_2^2$. If $Y_0<Y_*$, \autoref{prop:exit}(c) reads $\bP(\sigma<\infty)\le\min\{1,\sqrt{1+r_0}\,e^{-(Y_*-Y_0)/(2\lambda)}\}$. Thus the case of a small unforced part needs no separate argument.
\end{rem}

\section{One-component damping and a blow-up criterion for the linear noise}\label{sec:tamed}
We return to the deterministic model \eqref{eq:intro-tamed} and fix $\nu,\alpha>0$. Recall that our Leray projection~$\PP$ removes the spatial mean (\autoref{ss:spaces}). This matters here: With the projection which keeps constant fields, the mean of $u$ would evolve according to $\frac{\rd}{\rd t}\int_{\T^3}u\rd x=-\alpha\int_{\T^3}u_3^5\rd x\,\bfe_3$, and the mean-zero framework would not be invariant. With our convention,
\begin{equation}\label{eq:tamed}
 \del_tu+\cB(u)+\nu Au+\alpha\PP(u_3^5\bfe_3)=0,\qquad u(0)=u_0\in\Hs1,
\end{equation}
preserves the class of mean-zero divergence-free fields. We refer to \eqref{eq:tamed} as the damped system. Before the projection, its damping term acts only on the vertical component. By contrast, the taming term in the tamed Navier--Stokes equations of R\"ockner and Zhang \cite{RZ:09} acts on the whole velocity field.

Global regularity under one-component damping is not new in spirit. Yamazaki \cite[Theorem~2.2]{Yam:18} proved global well-posedness for the four-dimensional Navier--Stokes equations with damping $|u_j|^8u_j$ in two components, by combining the energy identity with a two-component regularity criterion, and remarked that an analogous result holds in three dimensions with damping in $u_3$ only, see \cite[pp.~451--452]{Yam:18}.
For the quintic power, the energy identity yields $u_3\in \rL^6_{t,x}$, which is exactly the endpoint pair~$(6,6)$ of the one-component criterion of Zhou and Pokorn\'y~\cite{ZP:10}. The periodic, mean-zero, $\rH^1$ setting of \eqref{eq:tamed} is not literally covered by these results, and in \autoref{ss:tamed-global}, we give a self-contained proof of global well-posedness, as a corollary of an abstract a priori bound (\autoref{thm:tamed}).

We include the proof because it can be formulated as a pathwise statement: An a priori bound for functions satisfying two integrated inequalities with an arbitrary constant $A\ge1$ (\autoref{prop:abstract}). With $A=1$, this gives~\autoref{thm:tamed}. With a random constant $A$, obtained from random integrating factors, it gives a one-component blow-up criterion for the linear noise \eqref{eq:lin}, valid for every $\beta>0$ (\autoref{prop:lin-criterion}). No existing criterion can be quoted for the latter, since the transformed equation is not the Navier--Stokes system and the constants are random.

\subsection{Pointwise estimates}\label{ss:tamed-convection}
The following two lemmas contain the structure exploited in one-component regularity criteria~\cite{KZ:06,CT:08,ZP:10}. They do not involve any equation.

\begin{lem}\label{lem:tamed-convection}
For $u\in\Hs2$,
\begin{align}
 |(\cB(u),\Delta u)_2|&\le C\int_{\T^3}|\nabla_hu||\nabla u|^2\rd x,\label{eq:full-gradient}\\
 |(\cB(u),\Delta_hu)_2|&\le C\int_{\T^3}|u_3||\nabla u||\nabla\nabla_hu|\rd x.\label{eq:horizontal-gradient}
\end{align}
\end{lem}
\begin{proof}
Both estimates are integrations by parts. The estimate \eqref{eq:horizontal-gradient} is contained in \cite[\mbox{(2.7)--(2.8)}]{CT:08}, where the periodic setting is treated. For \eqref{eq:full-gradient}, see \cite{KZ:06,ZP:10}. We only indicate the structure.
By~\eqref{eq:conv-commutator}, $(\cB(u),\Delta u)_2$ equals $-\sum_{i,j,k=1}^3\int\del_ju_i\,\del_iu_k\,\del_ju_k$. In each term with $(i,j)\ne(3,3)$, one of the factors~$\del_iu_k$ and $\del_ju_i$ is a horizontal derivative, and for $i=j=3$, we use $\del_3u_3=-\mdiv_hu_h$. This gives~\eqref{eq:full-gradient}. For \eqref{eq:horizontal-gradient}, only $j\in\{1,2\}$ occurs. By the Cayley--Hamilton identity from the proof of \autoref{lem:reference}, the terms with $i,k\le2$ add up to $-\int\del_3u_3\bigl(|\nabla_hu_h|^2-\det\nabla_hu_h\bigr)$, and an integration by parts in $x_3$ moves the derivative from $u_3$ to the other factors. In all remaining terms, $u_3$ is differentiated in a horizontal direction, and an integration by parts in this direction has the same effect. In all cases, one obtains products of $u_3$, one factor from $\nabla u$, and one factor from~$\nabla\nabla_hu$.
\end{proof}

For $u\in\Hs2$, we use the notation
\begin{equation}\label{eq:tamed-notation}
 E=\norm{\nabla u}_2^2,\qquad D=\norm{\Delta u}_2^2,\qquad E_h=\norm{\nabla_hu}_2^2,\qquad D_h=\norm{\nabla\nabla_hu}_2^2,\qquad U=\norm{u_3}_6 .
\end{equation}
Note that $E_h\le E$, $D_h\le D$, and $E\le D$ by the Poincar\'e inequality. 
The following anisotropic Sobolev inequality is classical, see \cite[p.~512]{TZ:96}, where it is attributed to Solonnikov \cite{Sol:72}. It is the case $r=6$ of \cite[(1.10)]{CT:08}, where it is recalled for the periodic cube and used, as below, to estimate derivatives of $u$ with only a small amount of vertical regularity. 

\begin{lem}[An anisotropic Sobolev inequality]\label{lem:anisotropic}
There is a constant $C$ such that for all $f\in \rH^1(\T^3;\R^m)$,
\begin{equation}\label{eq:anisotropic}
 \norm f_6\le C\prod_{j=1}^3\bigl(\norm{\del_jf}_2+\norm f_2\bigr)^{1/3}.
\end{equation}
Consequently, for $u\in\Hs2$, we have
\begin{equation}\label{eq:anisotropic-grad}
 \norm{\nabla u}_6\le C(E+D_h)^{1/3}D^{1/6},
\end{equation}
and, by interpolation between $\rL^2$ and $\rL^6$,
\begin{equation}\label{eq:anisotropic-interpolation}
 \norm{\nabla u}_3\le CE^{1/4}(E+D_h)^{1/6}D^{1/12},\qquad\norm{\nabla u}_4^2\le CE^{1/4}(E+D_h)^{1/2}D^{1/4}.
\end{equation}
\end{lem}
\begin{proof}
Inequality \eqref{eq:anisotropic} is the case $r=6$ of \cite[(1.10)]{CT:08}, which is stated there for scalar functions on the periodic cube $(0,1)^3$. It carries over to $\T^3$ by rescaling, and to vector-valued $f$ by applying it to each component.
Applied to $f=\nabla u$, with $\norm{\del_j\nabla u}_2\le D_h^{1/2}$ for $j=1,2$ and $\norm{\del_3\nabla u}_2+\norm{\nabla u}_2\le2D^{1/2}$, it gives \eqref{eq:anisotropic-grad}, and \eqref{eq:anisotropic-interpolation} follows from $\norm{\nabla u}_3\le\norm{\nabla u}_2^{1/2}\norm{\nabla u}_6^{1/2}$ and $\norm{\nabla u}_4\le\norm{\nabla u}_2^{1/4}\norm{\nabla u}_6^{3/4}$.
\end{proof}

\begin{lem}\label{lem:pointwise}
There is a constant $C_\nu$, depending only on $\nu$, such that for $u\in\Hs2$
\begin{align}
 2|(\cB(u),\Delta_hu)_2|&\le\nu D_h+C_\nu\bigl(U^3E^{3/4}D^{1/4}+U^2E^{5/6}D^{1/6}\bigr),\label{eq:pointwise-hor}\\
 2|(\cB(u),\Delta u)_2|&\le\nu D+C_\nu E_h^{2/3}E^{1/3}(E+D_h)^{2/3}.\label{eq:pointwise-full}
\end{align}
\end{lem}
\begin{proof}
By \eqref{eq:horizontal-gradient}, H\"older's inequality with exponents $(6,3,2)$, \eqref{eq:anisotropic-interpolation} and $(E+D_h)^{1/6}\le E^{1/6}+D_h^{1/6}$,
\[
 |(\cB(u),\Delta_hu)_2|\le CU\norm{\nabla u}_3D_h^{1/2}\le CUE^{1/4}D^{1/12}D_h^{2/3}+CUE^{5/12}D^{1/12}D_h^{1/2}.
\]
Young's inequality with exponents $(\frac32,3)$ for the first term and $(2,2)$ for the second term, each time allocating $\frac\nu4D_h$, gives~\eqref{eq:pointwise-hor}. By \eqref{eq:full-gradient}, H\"older's inequality with exponents $(2,4,4)$ and \eqref{eq:anisotropic-interpolation},
\[
 |(\cB(u),\Delta u)_2|\le CE_h^{1/2}\norm{\nabla u}_4^2\le CE_h^{1/2}E^{1/4}(E+D_h)^{1/2}D^{1/4}\le\frac\nu2D+C_\nu E_h^{2/3}E^{1/3}(E+D_h)^{2/3},
\]
where we used Young's inequality with exponents $(4,\frac43)$. This gives~\eqref{eq:pointwise-full}.
\end{proof}

Note that the right-hand side of \eqref{eq:pointwise-full} contains only horizontal second derivatives, while the right-hand side of \eqref{eq:pointwise-hor} contains $D$ with small exponents, multiplied by powers of $U=\norm{u_3}_6$.

\subsection{An abstract a priori bound}\label{ss:abstract}
Throughout this subsection, let $0<T_0\le\infty$, let $u\in \rC([0,T_0);\Hs1)\cap \rL^2_{\rm loc}([0,T_0);\Hs2)$, and let $E$, $D$, $E_h$, $D_h$ and $U$ be as in \eqref{eq:tamed-notation}, evaluated at~$u(t)$. No equation is assumed. For a constant $A\ge1$, we consider the condition
\begin{equation}\tag{$\mathrm H_A$}\label{eq:HA}
\left\{\begin{aligned}
 E_h(t)+\nu\int_a^tD_h\rd s&\le A\Bigl[E_h(a)+C_\nu\int_a^t\bigl(U^3E^{3/4}D^{1/4}+U^2E^{5/6}D^{1/6}\bigr)\rd s\Bigr],\\
 E(t)+\nu\int_a^tD\rd s&\le A\Bigl[E(a)+C_\nu\int_a^tE_h^{2/3}E^{1/3}(E+D_h)^{2/3}\rd s\Bigr],
\end{aligned}\right.
\qquad0\le a\le t<T_0,
\end{equation}
with the constant $C_\nu$ of~\autoref{lem:pointwise}. 
The integrals are finite: On $[a,t]$, $E$ is bounded and $D$ is integrable by the assumption on $u$, and $U\le CE^{1/2}$ by the Sobolev embedding, since $u_3$ has mean zero. Hence the integrands are bounded by $C(\sup_{[a,t]}E)^{9/4}D^{1/4}$, $C(\sup_{[a,t]}E)^{11/6}D^{1/6}$ and $(\sup_{[a,t]}E)(E+D)^{2/3}$, which are integrable on bounded intervals.
For $0\le a<b<T_0$ and $I=[a,b]$, we put
\begin{equation}\label{eq:tamed-interval-quantities}
 e=\int_IE\rd t,\qquad h=\int_IU^6\rd t,\qquad R=\sup_IE_h+\nu\int_ID_h\rd t,\qquad S=\sup_IE+\nu\int_ID\rd t .
\end{equation}
These are finite. In the applications, $e$ and $h$ are controlled by an energy identity, while $R$ and $S$ are the quantities to be estimated.

\begin{lem}\label{lem:tamed-interval}
There is a constant $C_*\ge1$, depending only on $\nu$, such that the following holds. If~\eqref{eq:HA} is satisfied for some $A\ge1$, then for every interval $I=[a,b]$ with $0\le a<b<T_0$, we have, with $C_A=C_*A^{7/3}$, the estimate
\begin{equation}\label{eq:tamed-interval}
 S\le C_A\Bigl[E(a)+\bigl(E_h(a)+e\bigr)^{4/3}e^{1/3}\Bigr]+C_Ah^{2/3}e^{2/3}S+C_Ah^{4/9}e\,S^{2/3}.
\end{equation}
\end{lem}
\begin{proof}
All suprema and integrals below are taken over $I$, and $C_\nu$ denotes constants depending only on~$\nu$.

\emph{Horizontal part.} The right-hand side of the first inequality of \eqref{eq:HA} is nondecreasing in $t$. Bounding~$\sup_IE_h$ and $\nu\int_ID_h\rd t$ separately by its value at $t=b$, we obtain
\[
 R\le2A\Bigl[E_h(a)+C_\nu\int_I\bigl(U^3E^{3/4}D^{1/4}+U^2E^{5/6}D^{1/6}\bigr)\rd t\Bigr].
\]
For the first term we extract $E^{1/2}$ and apply H\"older's inequality in time with exponents $(2,4,4)$, and for the second term we extract $E^{1/3}$ and use the exponents~$(3,2,6)$. Since $\sup E\le S$ and $\int D\le S/\nu$, this gives
\begin{align*}
 \int_IU^3E^{3/4}D^{1/4}\rd t&\le\Bigl(\sup_IE\Bigr)^{1/2}h^{1/2}e^{1/4}\Bigl(\int_ID\rd t\Bigr)^{1/4}\le\nu^{-1/4}h^{1/2}e^{1/4}S^{3/4},\\
 \int_IU^2E^{5/6}D^{1/6}\rd t&\le\Bigl(\sup_IE\Bigr)^{1/3}h^{1/3}e^{1/2}\Bigl(\int_ID\rd t\Bigr)^{1/6}\le\nu^{-1/6}h^{1/3}e^{1/2}S^{1/2},
\end{align*}
and therefore
\begin{equation}\label{eq:tamed-R}
 R\le C_\nu A\bigl(E_h(a)+h^{1/2}e^{1/4}S^{3/4}+h^{1/3}e^{1/2}S^{1/2}\bigr).
\end{equation}

\emph{Full gradient.} In the same way, the second inequality of \eqref{eq:HA}, H\"older's inequality in time with exponents $(3,\frac32)$, $\sup E_h\le R$ and $\int D_h\le R/\nu$ give
\[
 S\le2A\Bigl[E(a)+C_\nu\Bigl(\sup_IE_h\Bigr)^{2/3}e^{1/3}\Bigl(\int_I(E+D_h)\rd t\Bigr)^{2/3}\Bigr]\le2A\Bigl[E(a)+C_\nu R^{2/3}e^{1/3}\bigl(e+\nu^{-1}R\bigr)^{2/3}\Bigr],
\]
and hence
\begin{equation}\label{eq:tamed-S}
 S\le C_\nu A\bigl(E(a)+(R+e)^{4/3}e^{1/3}\bigr).
\end{equation}

\emph{Combination.} By \eqref{eq:tamed-R} and $A\ge1$, $R+e\le C_\nu A\bigl(E_h(a)+e+h^{1/2}e^{1/4}S^{3/4}+h^{1/3}e^{1/2}S^{1/2}\bigr)$. We insert this into \eqref{eq:tamed-S}, use $(x_1+x_2+x_3)^{4/3}\le3^{1/3}(x_1^{4/3}+x_2^{4/3}+x_3^{4/3})$ for $x_j\ge0$, and note that
\[
 \bigl(h^{1/2}e^{1/4}S^{3/4}\bigr)^{4/3}e^{1/3}=h^{2/3}e^{2/3}S,\qquad\bigl(h^{1/3}e^{1/2}S^{1/2}\bigr)^{4/3}e^{1/3}=h^{4/9}e\,S^{2/3}.
\]
The factors of $A$ combine to $A\cdot A^{4/3}=A^{7/3}$, and \eqref{eq:tamed-interval} follows.
\end{proof}

The term $C_Ah^{2/3}e^{2/3}S$ in \eqref{eq:tamed-interval} is linear in~$S$. It can be absorbed into the left-hand side if~$h^{2/3}e^{2/3}$ is small, which is possible on short intervals as soon as $h$ is bounded independently of~$S$.

\begin{prop}[Abstract a priori bound]\label{prop:abstract}
Suppose that \eqref{eq:HA} holds for some $A\ge1$, and that
\begin{equation}\label{eq:abstract-budget}
 \int_0^{T_0}E\rd t\le M_e,\qquad\int_0^{T_0}U^6\rd t\le M_h
\end{equation}
for some $M_e\ge0$ and $M_h>0$. Then there is a constant $K$, depending only on $\nu$, $A$, $M_e$, $M_h$ and~$E(0)$, such that
\[
 \sup_{t<T_0}E(t)+\nu\int_0^{T_0}D\rd t\le K .
\]
\end{prop}
\begin{proof}
\emph{Choice of the threshold.} Let $C_A$ be as in \autoref{lem:tamed-interval}, and put
\[
 \eps_*=\bigl(8C_A^{3/2}M_h\bigr)^{-1}.
\]
Let $I=[a,b]$ with $0\le a<b<T_0$ and $e\le\eps_*$. Since $h\le M_h$, we have $C_Ah^{2/3}e^{2/3}\le C_AM_h^{2/3}\eps_*^{2/3}=\frac14$, and Young's inequality with exponents $(\frac32,3)$ gives
\[
 C_Ah^{4/9}e\,S^{2/3}\le\frac14S+\frac{64}{27}C_A^3h^{4/3}e^3\le\frac14S+\frac{64}{27}C_A^3M_h^{4/3}\eps_*^3 .
\]
Since $S<\infty$, we may absorb both terms into the left-hand side of~\eqref{eq:tamed-interval}. Using also $E_h(a)\le E(a)$, we obtain
\begin{equation}\label{eq:tamed-F}
 S\le\Upsilon(E(a)),\qquad\Upsilon(y):=2C_A\Bigl[y+(y+\eps_*)^{4/3}\eps_*^{1/3}\Bigr]+\frac{128}{27}C_A^3M_h^{4/3}\eps_*^3 .
\end{equation}
The function $\Upsilon$ is increasing on $[0,\infty)$, satisfies $\Upsilon(y)\ge y$, and depends only on $\nu$, $A$ and~$M_h$.

\emph{A partition adapted to the dissipation.} Let $t_0=0$, and define recursively
\[
 t_{j+1}=\inf\Bigl\{t\in(t_j,T_0):\int_{t_j}^tE\rd s\ge\eps_*\Bigr\},\qquad\inf\emptyset=T_0 .
\]
As long as $t_{j+1}<T_0$, the interval $[t_j,t_{j+1}]$ consumes exactly $\eps_*$ of the budget $\int_0^{T_0}E\rd t\le M_e$. Hence,~$t_j=T_0$ for some $j\le J:=\lfloor M_e/\eps_*\rfloor+1$. The partition depends on $u$, but the bound $J$ on the number of intervals depends only on $\nu$, $A$, $M_e$ and~$M_h$.

\emph{Iteration.} Let $K_0=E(0)$ and $K_{j+1}=\Upsilon(K_j)$. This sequence is nondecreasing. Suppose that~$t_j<T_0$ and $E(t_j)\le K_j$, which holds for $j=0$. For every $b<t_{j+1}$, and also for $b=t_{j+1}$ if $t_{j+1}<T_0$, the interval $[t_j,b]$ satisfies $e\le\eps_*$, so \eqref{eq:tamed-F} and the monotonicity of $\Upsilon$ give
\[
 \sup_{[t_j,b]}E+\nu\int_{t_j}^bD\rd t\le\Upsilon(K_j)=K_{j+1}.
\]
In particular $E(t_{j+1})\le K_{j+1}$ if $t_{j+1}<T_0$. By induction, the same bound holds on all at most $J$ intervals, and therefore
\[
 \sup_{t<T_0}E(t)\le K_J,\qquad\nu\int_0^{T_0}D\rd t\le JK_J .
\]
This proves the assertion with $K=(J+1)K_J$.
\end{proof}

\subsection{Global regularity for the damped system}\label{ss:tamed-global}
\begin{lem}\label{lem:tamed-local}
For every $u_0\in\Hs1$, \eqref{eq:tamed} has a unique maximal solution $u$, which belongs to $\rC([0,T_*);\Hs1)\cap \rL^2_{\rm loc}([0,T_*);\Hs2)$, and if $T_*<\infty$, then
\[
 \sup_{t<T_*}\norm{\nabla u(t)}_2+\int_0^{T_*}\norm{\Delta u(t)}_2^2\rd t=\infty .
\]
\end{lem}
\begin{proof}
We apply \autoref{thm:AV} with $G=0$ and $F(u)=-\cB(u)-\alpha\PP(u_3^5\bfe_3)$. The convection term is treated as in~\autoref{prop:local}. By $|a^5-c^5|\le5(|a|^4+|c|^4)|a-c|$, H\"older's inequality and $\rH^{6/5}\hookrightarrow \rL^{10}$, the damping satisfies (H2) with $(\beta_1,\rho_1)=(\frac35,4)$. Since $(2\beta_1-1)(\rho_1+1)=1$, it is critical, reflecting the fact that $\rH^1$ is the scaling-critical space for the quintic heat nonlinearity in three dimensions. Hence, only the blow-up criterion of \autoref{thm:AV}(a) is available, which is the assertion.
\end{proof}

\begin{lem}[Energy balances]\label{lem:tamed-energy}
For $t<T_*$,
\begin{align}
 \norm{u(t)}_2^2+2\nu\int_0^t\norm{\nabla u}_2^2\rd s+2\alpha\int_0^t\norm{u_3}_6^6\rd s&=\norm{u_0}_2^2,\label{eq:tamed-L2}\\
 \norm{\nabla u(t)}_2^2+2\nu\int_0^t\norm{\Delta u}_2^2\rd s+10\alpha\int_0^t\norm{u_3^2\nabla u_3}_2^2\rd s&=\norm{\nabla u_0}_2^2+2\int_0^t(\cB(u),\Delta u)_2\rd s,\label{eq:tamed-H1}\\
 \norm{\nabla_hu(t)}_2^2+2\nu\int_0^t\norm{\nabla\nabla_hu}_2^2\rd s+10\alpha\int_0^t\norm{u_3^2\nabla_hu_3}_2^2\rd s&=\norm{\nabla_hu_0}_2^2+2\int_0^t(\cB(u),\Delta_hu)_2\rd s.\label{eq:tamed-hor}
\end{align}
In particular,
\begin{equation}\label{eq:tamed-component}
 \int_0^{T_*}\norm{u_3(t)}_6^6\rd t\le\frac{\norm{u_0}_2^2}{2\alpha}.
\end{equation}
\end{lem}
\begin{proof}
For $T<T_*$ we have $\del_tu\in \rL^2(0,T;\rL^2_\sigma)$, so $\norm u_2^2$, $\norm{\nabla u}_2^2$ and $\norm{\nabla_hu}_2^2$ are absolutely continuous, with derivatives $2(\del_tu,u)_2$, $2(\del_tu,-\Delta u)_2$ and $2(\del_tu,-\Delta_hu)_2$. Since $u$, $\Delta u$ and~$\Delta_hu$ are mean-zero and divergence-free, and $\PP$ is self-adjoint, we have $(\PP(u_3^5\bfe_3),u)_2=\int u_3^6$, $(\PP(u_3^5\bfe_3),-\Delta u)_2=(u_3^5,-\Delta u_3)_2=5\norm{u_3^2\nabla u_3}_2^2$, and the analogous identity with $\Delta_h$. For a.e.~$t$, $u_3\in \rH^2\subset \rL^\infty$, which justifies the integration by parts. The convection term vanishes in the first identity.
\end{proof}

\begin{cor}[Global well-posedness of the damped system]\label{thm:tamed}
For every $u_0\in\Hs1$, the maximal solution of \eqref{eq:tamed} is global, i.e., $T_*=\infty$, and
\[
 \sup_{t\ge0}\norm{\nabla u(t)}_2^2+\nu\int_0^\infty\norm{\Delta u(t)}_2^2\rd t\le K,\qquad\int_0^\infty\norm{u_3(t)}_6^6\rd t\le\frac{\norm{u_0}_2^2}{2\alpha},
\]
with $K$ depending only on $\nu$, $\alpha$, $\norm{u_0}_2$ and $\norm{\nabla u_0}_2$.
\end{cor}
\begin{proof}
If $u_0=0$, then $u\equiv0$, so let $u_0\ne0$. Subtracting \eqref{eq:tamed-hor} at times $a$ and $t$, dropping the nonnegative damping term and using \eqref{eq:pointwise-hor}, for $0\le a\le t<T_*$, we obtain 
\[
 E_h(t)+2\nu\int_a^tD_h\rd s\le E_h(a)+\int_a^t\bigl[\nu D_h+C_\nu\bigl(U^3E^{3/4}D^{1/4}+U^2E^{5/6}D^{1/6}\bigr)\bigr]\rd s,
\]
and absorbing $\nu\int_a^tD_h\rd s$, which is finite, into the left-hand side gives the first inequality of \eqref{eq:HA} with $A=1$. In the same way, \eqref{eq:tamed-H1} and \eqref{eq:pointwise-full} give the second one. By \eqref{eq:tamed-L2} and \eqref{eq:tamed-component}, \eqref{eq:abstract-budget} holds with
\[
 M_e=\frac{\norm{u_0}_2^2}{2\nu},\qquad M_h=\frac{\norm{u_0}_2^2}{2\alpha}.
\]
The second bound is where the damping enters; it does not depend on any higher norm of~$u$. \autoref{prop:abstract} with $T_0=T_*$ and $A=1$ gives $\sup_{t<T_*}\norm{\nabla u}_2^2+\nu\int_0^{T_*}\norm{\Delta u}_2^2\rd t\le K$, with $K$ depending only on $\nu$, $\alpha$, $\norm{u_0}_2$ and $\norm{\nabla u_0}_2$. If $T_*<\infty$, this contradicts the blow-up criterion of~\autoref{lem:tamed-local}. Hence, $T_*=\infty$. Uniqueness is part of \autoref{lem:tamed-local}, and the last assertion is~\eqref{eq:tamed-component}.
\end{proof}

\begin{rem}\label{rem:tamed-damping}
(1) The damping was used only through the bound \eqref{eq:tamed-component} and through the sign of the damping terms in \eqref{eq:tamed-H1} and \eqref{eq:tamed-hor}, which were dropped. The bound \eqref{eq:tamed-component} is what allows us to fix the threshold $\eps_*$ in terms of the data alone. Without the damping, the Sobolev inequality only gives $h\le C\int_IE^3\rd t\le C\,eS^2$, and the critical term in \eqref{eq:tamed-interval} becomes $C_*h^{2/3}e^{2/3}S\le Ce^{4/3}S^{7/3}$, which is superlinear in $S$ and can no longer be absorbed by a threshold depending only on the data.

(2) For the same reason, \autoref{prop:abstract} with $A=1$ applies to the Navier--Stokes equations ($\alpha=0$) if the bound on $\int\norm{u_3}_6^6\rd t$ is assumed instead of derived: If $u$ is a strong solution of the Navier--Stokes equations on $[0,T)$ with $T<\infty$ and $\int_0^T\norm{u_3}_6^6\rd t<\infty$, then $\sup_{t<T}\norm{\nabla u(t)}_2+\int_0^T\norm{\Delta u}_2^2\rd t<\infty$, and $u$ can be continued beyond~$T$. This recovers, on $\T^3$, the endpoint case $(p,q)=(6,6)$ of the one-component criterion $u_3\in \rL^p(0,T;\rL^q)$, $\frac2p+\frac3q\le\frac34+\frac1{2q}$, $q>\frac{10}3$, of Zhou and Pokorn\'y \cite[Theorem~2.1]{ZP:10}.
The strategy of combining the horizontal and the full enstrophy balance goes back to \cite{KZ:06,CT:08}, see also \cite{AGR:23}, and the idea of feeding the integrability provided by a partial damping into a component-reduction criterion is due to Yamazaki~\cite{Yam:18}.
\end{rem}

\begin{rem}\label{rem:higher-powers}
For the damping $\alpha\PP(|u_3|^{2m-2}u_3\bfe_3)$, the energy identity gives $u_3\in \rL^{2m}_{t,x}$, and the damping again contributes nonnegative terms to the enstrophy balances. For $m\ge4$, H\"older's inequality on $[0,T]\times\T^3$ turns this into a bound on $\int_0^T\norm{u_3}_6^6\rd t$, so \autoref{prop:abstract} applies on bounded time intervals. The damping is then supercritical in the sense of \autoref{thm:AV}, however, and a local theory would have to rely on its monotonicity. For $m=2$, the bound $u_3\in \rL^4_{t,x}$ is not sufficient for the argument. Thus, the quintic damping is the smallest odd power covered by \autoref{prop:abstract}, and the one adapted to $\rH^1$ data, since $\rH^1\hookrightarrow \rL^6$.
\end{rem}

\subsection{A blow-up criterion for the linear noise}\label{ss:lin-criterion}
For the linear noise \eqref{eq:lin}, the stochastic terms in the energy and enstrophy balances are bounded multiples of the balanced quantities, and they can be removed by random integrating factors, as in the proof of~\autoref{lem:continuation}. This leads to \eqref{eq:HA} with a random constant~$A$. The horizontal balance requires an It\^o formula for the seminorm $\norm{\nabla_h\cdot}_2$, which we derive from \cite[Lemma~A.5]{AV:25} by means of a Fourier multiplier.

\begin{lem}\label{lem:lin-balances}
Let $\beta>0$, $u_0\in\Hs1$, and let $(u,\tau_*)$ be the maximal solution of~\eqref{eq:lin}. Then a.s., for all $t<\tau_*$, we have
\begin{align}
 \norm{u(t)}_2^2+2\nu\int_0^t\norm{\nabla u}_2^2\rd s&=\norm{u_0}_2^2+\beta^2\int_0^t\norm{Ku}_2^2\rd s+2\beta\int_0^t\norm{u_3}_2^2\rd W_s,\label{eq:lin-L2}\\
 \norm{\nabla_hu(t)}_2^2+2\nu\int_0^t\norm{\nabla\nabla_hu}_2^2\rd s&=\norm{\nabla_hu_0}_2^2+\int_0^t\bigl[2(\cB(u),\Delta_hu)_2+\beta^2\norm{\nabla_hKu}_2^2\bigr]\rd s\notag\\
 &\quad+2\beta\int_0^t\norm{\nabla_hu_3}_2^2\rd W_s.\label{eq:lin-hor}
\end{align}
\end{lem}
\begin{proof}
Let $M$ be the Fourier multiplier with symbol $M_\xi=|\xi|^{-1}$ for \eqref{eq:lin-L2}, and $M_\xi=|\xi_h|/|\xi|$ for~\eqref{eq:lin-hor}. Since $0\le M_\xi\le1$, the operator $M$ is bounded on every $\Hs s$, and it commutes with~$\Delta$, $\PP$ and~$K$. Applying $M$ to the equation in \autoref{def:solution}, we see that $z=Mu$ belongs to the same class as $u$ and satisfies
\[
 z(t)=Mu_0+\int_0^t\bigl[\nu\Delta z-M\cB(u)\bigr]\rd s+\beta\int_0^tMKu\rd W_s .
\]
The It\^o formula \cite[Lemma~A.5]{AV:25}, localized as in the proof of \autoref{lem:energy}, therefore gives
\begin{multline*}
 \norm{\nabla z(t)}_2^2+2\nu\int_0^t\norm{\Delta z}_2^2\rd s=\norm{\nabla Mu_0}_2^2+\int_0^t\bigl[2(\Delta M^2u,\cB(u))_2+\beta^2\norm{\nabla MKu}_2^2\bigr]\rd s\\
 +2\beta\int_0^t(\nabla MKu,\nabla Mu)_2\rd W_s ,
\end{multline*}
where we used $(-\Delta z,M\cB(u))_2=-(\Delta M^2u,\cB(u))_2$. By Plancherel's theorem, for $M_\xi=|\xi|^{-1}$, we have $\norm{\nabla Mf}_2=\norm f_2$, $\norm{\Delta Mf}_2=\norm{\nabla f}_2$ and $\Delta M^2=-\Id$, so $(\Delta M^2u,\cB(u))_2=-(u,\cB(u))_2=0$. For $M_\xi=|\xi_h|/|\xi|$, we have $\norm{\nabla Mf}_2=\norm{\nabla_hf}_2$, $\norm{\Delta Mf}_2=\norm{\nabla\nabla_hf}_2$ and $\Delta M^2=\Delta_h$. Finally, $(\nabla MKu,\nabla Mu)_2$ equals $(Ku,u)_2=\norm{u_3}_2^2$ in the first case and $(\nabla_hKu,\nabla_hu)_2=\norm{\nabla_hu_3}_2^2$ in the second case, by \autoref{lem:geometry}(iv).
\end{proof}

\begin{prop}[A one-component blow-up criterion for the linear noise]\label{prop:lin-criterion}
Let $\beta>0$, $u_0\in\Hs1$, and let $(u,\tau_*)$ be the maximal solution of~\eqref{eq:lin}. Then
\[
 \bP\Bigl(\tau_*<\infty,\ \int_0^{\tau_*}\norm{u_3(t)}_6^6\rd t<\infty\Bigr)=0 .
\]
\end{prop}
\begin{proof}
Fix $T\in\N$, and let $G_T=\{\tau_*\le T,\ \int_0^{\tau_*}\norm{u_3}_6^6\rd t<\infty\}$. We show that $\sup_{t<\tau_*}\norm{\nabla u(t)}_2<\infty$ a.s.\ on~$G_T$. By \eqref{eq:blowup}, this implies $\bP(G_T)=0$, and the assertion follows by taking the union over~$T\in\N$. We use the notation \eqref{eq:tamed-notation} and write $E_0=\norm u_2^2$.

\emph{Step 1: Random integrating factors.} Let $X\in\{E_0,E,E_h\}$. By \autoref{lem:lin-balances} and \autoref{lem:energy}, we have, a.s.\ for all $t<\tau_*$,
\begin{equation}\label{eq:lin-X0}
 X(t)=X(0)+\int_0^t\bigl(-2\nu D_X+N_X+\Phi_X\bigr)\rd s+\int_0^t\Psi_X\rd W_s,
\end{equation}
where $(D_X,N_X)$ equals $(E,0)$, $(D,2(\cB(u),\Delta u)_2)$ and $(D_h,2(\cB(u),\Delta_hu)_2)$, the It\^o corrections are
\[
 \Phi_X=\beta^2\norm{Ku}_2^2,\qquad\beta^2\norm{\nabla Ku}_2^2,\qquad\beta^2\norm{\nabla_hKu}_2^2,
\]
and the noise integrands are
\[
 \Psi_X=2\beta\norm{u_3}_2^2,\qquad2\beta\norm{\nabla u_3}_2^2,\qquad2\beta\norm{\nabla_hu_3}_2^2,
\]
respectively. Since $K$ commutes with derivatives and its symbol has norm at most one, and since $|\nabla u_3|\le|\nabla u|$ pointwise,
\[
 0\le\Phi_X\le\beta^2X,\qquad0\le\Psi_X\le2\beta X .
\]
We put $c_X=\Phi_X/X$ and $b_X=\Psi_X/X$ on $\{X>0\}$, and $c_X=b_X=0$ on $\{X=0\}$. These are progressively measurable functions of $u$ with $0\le c_X\le\beta^2$ and $0\le b_X\le2\beta$, and $\Phi_X=c_XX$, $\Psi_X=b_XX$ also on $\{X=0\}$, because $\Phi_X$ and $\Psi_X$ vanish there by the bounds above. Hence, \eqref{eq:lin-X0} can be written as
\begin{equation}\label{eq:lin-X}
 \rd X=\bigl(-2\nu D_X+N_X+c_XX\bigr)\rd t+b_XX\rd W\qquad\text{on }[0,\tau_*).
\end{equation}
Let $I_X(t)=\int_0^t\mathbf 1_{[0,\tau_*)}b_X\rd W_s$, which is a continuous martingale on $[0,T]$, and put
\[
 \Gamma_X(t)=\exp\Bigl(-\int_0^tc_X\rd s+\frac12\int_0^tb_X^2\rd s-I_X(t)\Bigr),\qquad\Lambda_T=\exp\Bigl(3\beta^2T+\sum_X\sup_{t\le T}|I_X(t)|\Bigr)<\infty\quad\text{a.s.}
\]
Since $0\le c_X\le\beta^2$ and $0\le b_X\le2\beta$, we have $|\log\Gamma_X(t)|\le\beta^2T+2\beta^2T+\sup_{s\le T}|I_X(s)|\le\log\Lambda_T$, i.e., $\Lambda_T^{-1}\le\Gamma_X\le\Lambda_T$ on $[0,\tau_*\wedge T)$.
Moreover, $\rd\Gamma_X=\Gamma_X[(b_X^2-c_X)\rd t-b_X\rd W]$, and It\^o's product rule gives
\[
 \rd(\Gamma_XX)=\Gamma_X\rd X+X\rd\Gamma_X-\Gamma_Xb_X^2X\rd t=\Gamma_X(-2\nu D_X+N_X)\rd t
\]
on $[0,\tau_*)$ after localization with $\tau_n$. Unlike in the proof of \autoref{lem:continuation}, the convection term $N_X$ is not included in the integrating factor. This keeps $\Gamma_X$ bounded by $\Lambda_T$, independently of higher norms of $u$.
Hence, for $0\le a\le t<\tau_*\wedge T$, we get
\begin{equation}\label{eq:lin-weighted}
 \Gamma_X(t)X(t)+2\nu\int_a^t\Gamma_XD_X\rd s=\Gamma_X(a)X(a)+\int_a^t\Gamma_XN_X\rd s .
\end{equation}

\emph{Step 2: Verification of~\eqref{eq:HA}.} For $X=E_0$, we have $D_{E_0}=E$, $N_{E_0}=0$ and $\Gamma_{E_0}(0)=1$, so \eqref{eq:lin-weighted} with $a=0$, $\Gamma_{E_0}\ge\Lambda_T^{-1}$ and monotone convergence give
\[
 \int_0^{\tau_*\wedge T}E\rd t\le\frac{\Lambda_T\norm{u_0}_2^2}{2\nu}.
\]
For $X=E_h$, we insert \eqref{eq:pointwise-hor} into \eqref{eq:lin-weighted} and absorb $\nu\int_a^t\Gamma_{E_h}D_h\rd s$, which is finite, into the left-hand side. This gives, for $0\le a\le t<\tau_*\wedge T$,
\[
 \Gamma_{E_h}(t)E_h(t)+\nu\int_a^t\Gamma_{E_h}D_h\rd s\le\Gamma_{E_h}(a)E_h(a)+C_\nu\int_a^t\Gamma_{E_h}\bigl(U^3E^{3/4}D^{1/4}+U^2E^{5/6}D^{1/6}\bigr)\rd s .
\]
Bounding $\Gamma_{E_h}$ from below by $\Lambda_T^{-1}$ on the left-hand side and from above by $\Lambda_T$ on the right-hand side, where all terms are nonnegative, and multiplying by $\Lambda_T$, we obtain the first inequality of \eqref{eq:HA} with $A=\Lambda_T^2$ and $T_0=\tau_*\wedge T$. In the same way, $X=E$ and \eqref{eq:pointwise-full} give the second one. Note that~$\Lambda_T$ does not depend on $a$ and $t$.

\emph{Step 3: Conclusion.} Fix $\omega\in G_T$ outside a null set. Then $\tau_*\le T$, and by Step~2 the path $u(\omega)\in \rC([0,\tau_*);\Hs1)\cap \rL^2_{\rm loc}([0,\tau_*);\Hs2)$ satisfies \eqref{eq:HA} on $[0,\tau_*)$ with $A=\Lambda_T(\omega)^2$, as well as \eqref{eq:abstract-budget} with
\[
 M_e=\frac{\Lambda_T(\omega)\norm{u_0}_2^2}{2\nu},\qquad M_h=1+\int_0^{\tau_*}\norm{u_3}_6^6\rd t<\infty .
\]
\autoref{prop:abstract} with $T_0=\tau_*(\omega)$ gives $\sup_{t<\tau_*}E(t)+\nu\int_0^{\tau_*}D\rd t<\infty$, which proves the claim.
\end{proof}

\begin{rem}\label{rem:lin-criterion}
(1) In \autoref{thm:tamed}, the bound \eqref{eq:tamed-component} comes first and implies global existence. In \autoref{thm:nonlinear} and \autoref{thm:linear}(c), the integral $\int_0^\infty\!\int_{\T^3}|u_3|^6\rd x\rd t$ is controlled only a posteriori, on an event of high probability on which global existence has been obtained by other means.

(2) \autoref{prop:lin-criterion} holds for every $\beta>0$ and every $u_0\in\Hs1$, without any condition on $\cS u_0$, and it is the stochastic analogue of the criterion in \autoref{rem:tamed-damping}(2). Since its constants are random, it does not provide moment bounds. The proof uses that the stochastic terms in the balances are bounded multiples of the balanced quantities, i.e., that the noise is linear. For the nonlinear noise \eqref{eq:nl}, the corresponding coefficients are of order $\beta^2\norm{\nabla\cR u}_2^{2k}$, but there \autoref{thm:nonlinear} already provides global existence with high probability.
\end{rem}

\subsection*{Acknowledgments}
I.K.\ was supported in part by NSF grant DMS-2205493 and Simons grant SFI-MPS-TSM-00025796.

Parts of this work were carried out with the assistance of GPT-6 Astra and Claude Opus~5.5. They were used to draft and edit parts of the manuscript, to propose and check arguments, to support the literature search, and for typesetting. The main ideas and the overall approach are due to the authors, who have checked all arguments and take full responsibility for the content of this paper.

\appendix
\section{The logarithm at zero}\label{app:log}
We prove~\autoref{lem:log}. Fix $n$ and work on $[0,\tau_n]$, where $\tau_n$ is as in the proof of \autoref{lem:Xminus}. Thus, $r\le\norm{\nabla u}_2\le n$ and $\int_0^{\tau_n}\norm{\Delta u}_2^2\le n$. By \eqref{eq:energy-v}, $E=r^2$ is a real continuous semimartingale. For $\eps\in(0,1]$, let
\[
 \psi_\eps(E)=\log\bigl(1+\sqrt{E+\eps^2}\bigr),\qquad r_\eps=\sqrt{r^2+\eps^2},\qquad a_\eps=\frac1{r_\eps(1+r_\eps)}.
\]
Then $\psi_\eps'(E)=\frac12a_\eps$ and $\psi_\eps''(E)=-\frac{1+2r_\eps}{4r_\eps^3(1+r_\eps)^2}$, and It\^o's formula gives
\begin{multline*}
 \psi_\eps(E(t))=\psi_\eps(E(0))-\int_0^ta_\eps\bigl(\nu D_v+T_v\bigr)\rd s+\beta\int_0^t\frac{r^{k+2}}{r_\eps(1+r_\eps)}\rd W_s\\
 +\int_0^t\Bigl[-\frac{\beta^2r^{2k+4}}{2r_\eps^2(1+r_\eps)^2}+\frac{\beta^2\eps^2r^{2k+2}}{2r_\eps^3(1+r_\eps)}\Bigr]\rd s,
\end{multline*}
where the integrand is the exact noise contribution $\psi_\eps'\beta^2r^{2k+2}+\frac12\psi_\eps''\cdot4\beta^2r^{2k+4}$. We let $\eps\downarrow0$ term by term.
Then:
\begin{enumerate}[(i)]
 \item $\psi_\eps(E)\to\psi(r)$ uniformly.
 \item Since $\eps^2r^2/(r^2+\eps^2)^{3/2}\le C\eps$, which can be shown by substituting $r=\eps s$, the positive part of the bracket is at most $C\beta^2n^{2k}\eps$, and its integral vanishes in the limit.
 \item The negative part of the integrand increases to $-\frac{\beta^2}2h_k(r)$ in absolute value as $\eps\downarrow0$; after moving it to the left-hand side, monotone convergence applies.
 \item The viscous term $a_\eps D_v$ increases to $a(r)D_v$ on $\{v\ne0\}$, and vanishes for a.e.\ $t$ in $\{v=0\}$. Monotone convergence applies.
 \item For $r>0$ we have $a_\eps\le a(r)$, and by \eqref{eq:Tv} the function $a(r)|T_v|$ is bounded by a constant multiple of $r^{1/2}D_v^{3/4}+(D_e^{1/2}+D_b^{1/2})r^{1/2}D_v^{1/4}+D_e+D_b$, which is integrable on $[0,\tau_n]$. At times with $v=0$ we have $T_v=0$. Dominated convergence applies.
 \item The martingale integrands converge to $\beta r^{k+1}/(1+r)$ and are bounded by $\beta n^{k+1}$. The stochastic integrals converge uniformly in $\rL^2(\Omega)$ by Doob's inequality.
\end{enumerate}
This gives \eqref{eq:log} on $[0,\tau_n]$, without any local-time contribution at $r=0$, and $n\to\infty$ completes the proof. Alternatively, the Young inequalities of \autoref{lem:young} can be applied before the limit $\eps\downarrow0$, with~$r(1+r)$ replaced by $r_\eps(1+r_\eps)$. The regularization errors vanish on the same localizations at fixed $\beta$ and~$\lambda$.

\end{document}